\documentclass[dvipdfmx,a4paper]{article}
\usepackage[margin=25mm]{geometry}
\usepackage{amsmath,amssymb}
\usepackage{bm}
\usepackage{ulem}
\usepackage{graphicx}
\usepackage{ascmac}
\usepackage{amsthm}
\usepackage{authblk}
\usepackage {algorithmic}
\usepackage{algorithm}
\usepackage{mathrsfs}
\usepackage{comment}
\usepackage[T1]{fontenc}
\usepackage{newpxtext,newtxmath}
\newtheorem{thm}{Theorem}[section]
\newtheorem{main_theorem}[thm]{Main threorem}
\newtheorem{dfn}[thm]{Definition}
\newtheorem{prp}[thm]{Proposition}
\newtheorem{lmm}[thm]{Lemma}
\newtheorem{cor}[thm]{Corollary}
\newtheorem{alg}[thm]{Algorithm}
\newtheorem{xmp}[thm]{Example}

\newcommand{\CC}{\mathbb C}
\newcommand{\EE}{\mathbb E}
\newcommand{\NN}{\mathbb N}
\newcommand{\PP}{\mathbb P}
\newcommand{\QQ}{\mathbb Q}
\newcommand{\RR}{\mathbb R}
\newcommand{\ZZ}{\mathbb Z}

\newcommand{\calNC}{\mathcal{NC}}
\newcommand{\calLNC}{\mathcal{LNC}}
\newcommand{\calSP}{\mathcal{SP}}
\newcommand{\calC}{\mathcal{C}}
\newcommand{\calD}{\mathcal{D}}
\newcommand{\calK}{\mathcal{K}}
\newcommand{\calL}{\mathcal{L}}
\newcommand{\calM}{\mathcal{M}}
\newcommand{\calS}{\mathcal{S}}
\newcommand{\calT}{\mathcal{T}}

\newcommand{\Loco}{\mathrm{Lo.co.}}
\newcommand{\IR}{\mathrm{IR}}
\newcommand{\DIR}{\mathrm{DIR}}
\newcommand{\sign}{\mathrm{sign}}
\newcommand{\smplx}{\lambda_\mathrm{smplx}}

\newcommand{\bmat}[1]{\begin{bmatrix} #1 \end{bmatrix}}

\newcommand{\mat}[1]{\begin{matrix} #1 \end{matrix}}
  \makeatletter
    
    \@addtoreset{equation}{section}
  \makeatother
  
    \makeatletter
    
    \@addtoreset{figure}{section}
  \makeatother
  
 \makeatletter
    
    \@addtoreset{thm}{section}
    \renewcommand{\ALG@name}{Algorithm(pseudo code)}
  \makeatother

\title{Blow-up Algorithm for Sum-of-Products Polynomials and Real Log Canonical Thresholds}
\author[1]{Joe Hirose}
\affil[1]{Department of Mathematical and Computing Science, Tokyo Institute of Technology, 2-12-1, 
Oookayama, Meguro-ku, Tokyo, 152-8552, Japan}
\date{\today}
\begin{document}
\maketitle
\begin{abstract}
When considering a real log canonical threshold (RLCT) that gives a Bayesian generalization error, 
in general, papers replace a mean error function with a relatively simple polynomial 
whose RLCT corresponds to that of the mean error function, 
and obtain its RLCT by resolving its singularities through an algebraic operation called blow-up.
Though it is known that the singularities of any polynomial can be resolved 
by a finite number of blow-up iterations, 
it is not clarified whether or not it is possible to resolve singularities of a specific polynomial 
by applying a specific blow-up algorithm.
Therefore this paper considers the blow-up algorithm for the polynomials 
called sum-of-products (sop) polynomials and its RLCT.

The main theorem can be divided into three parts. 

The first is the case of a bivariate sop binomial.
In this case, it is possible to construct the manifold whose any local coordinate 
the bivariate sop binomial is normal crossing on 
by repeatedly blow-up centered with the origin.
A normal crossing makes the calculation of RLCT easier.
The RLCT of a bivariate sop binomial can be obtained 
by considering the invariant through the blow-up.

The second is the case of a general sop binomial.
In this case, it is possible to construct the manifold whose any local coordinate 
the bivariate sop binomial is normal crossing on  
by repeatedly blow-up centered with the variables whose 
difference of the degree between terms are maximum and minimum.
The RLCT of a sop binomial can be obtained by using the result for the bivariate case.

The last is the case of a general sop polynomial.
Because it is difficult to judge whether or not 
a given sop polynomial is normal crossing in this case, 
this paper defines local normal crossing which is a weaker condition than normal crossing.
it is possible to construct the manifold whose any local coordinate 
the sop polynomial is local normal crossing on.
The manifold can be obtained by swapping terms and applying the blow-up algorithm 
for sop binomial repeatedly.
Because the algorithm does not necessarily make the sop polynomial 
normal crossing on all local coordinates, 
it is derived that not the exact value of RLCT but the upper bound of that.
This upper bound (called simplex upper bound in this paper) of RLCT corresponds to 
the inverse of the optimal value of the linear programming problem 
for multi-indexes of the sop polynomial, 
and is derived by constructing a local coordinate system with the pole of the zeta function 
corresponding to the upper bound by using the algorithm.

The main theorem can be applied to a general polynomial as well.
This paper compares the simplex upper bound with RLCT of the statistical model and 
the specific polynomial revealed in previous studies.
Furthermore, this paper reveals that the optimal solution of the linear programming problem 
for multi-indexes of polynomial derives the optimal weight of weighted blow-up.
\end{abstract}

\setcounter{tocdepth}{3}
\tableofcontents

\section{Introduction}
When considering an invariant that gives a Bayesian generalization error, 
that is a real log canonical threshold, 
in general, papers replace a mean error function with a relatively simple polynomial 
whose real log canonical threshold corresponds to that of the mean error function, 
and obtain its real log canonical threshold by resolving 
its singularities through an algebraic operation called blow-up.
Though it is known that the singularities of any polynomial can be resolved 
by a finite number of blow-up iterations, 
it is not clarified well whether or not it is possible to resolve singularities of a specific polynomial 
by applying a specific blow-up algorithm.
Therefore this paper proposes a blow-up algorithm that can be applied to the polynomials 
called sum-of-products polynomials and proves that it halts.
Furthermore, this paper considers real log canonical thresholds of sum-of-products polynomials 
by using the algorithm.

First, this section explains the foundation of Bayesian learning theory and 
details the relation to a real log canonical threshold and blow-up.
Then this section defines exclusive sum-of-products polynomials 
which is subject to previous studies and 
explains the novelty and utility of this paper.

\subsection{Bayesian learning theory}
This subsection explains Bayesian learning theory which is the main study background of this paper.
First, this paper explains how to construct Bayesian predictive distribution and 
that its generalization loss is given by real log canonical threshold.
Second, this paper introduces how to derive a real log canonical threshold.
This paper considers the pole of the zeta function in order to derive a real log canonical threshold.
it needs to resolve the singularities in order to derive the poles of the zeta function, 
so blow-up is introduced as the method of resolution of singularities.

\subsubsection{framework of Bayesian inference}\label{sec:framework_of_Bayes_method}
This subsubsection explains the framework of Bayesian statistics.
There are three distributions needed to make a predictive distribution.
It consists of \textbf{true distribution} $q(x)$ which generate data, 
\textbf{statistical model} $p(x|w)$ which is the set of distribution decided by parameter $w$ 
and \textbf{prior distribution} $\varphi(w)$ of parameter $w$.
The accuracy of Bayesian inference is also decided by the triplet $(q(x),p(x|w),\varphi(w))$, 
it is explained later.

Let $X^n=(X_1,X_2,\ldots,X_n)$ be a set of data generated from 
true distribution $q(x)$ independently and identically.

True distribution $q(x)$ is not given in practice.
In the framework of Bayesian inference, it is supposed that 
the data $X^n$ is generated from the statistical model $p(x|w)$ for some $w\in W$ instead and that 
the parameter $w$ follows some distribution $\varphi(w)$ which is called prior distribution.
The joint probability distribution between the data and the parameter can be defined 
under the assumption mentioned above.
The conditional probability distribution $p(w|X^n)$ given the data $X^n$ is written by the following.

\[p(w|X^n)=\frac{1}{Z_n}\varphi(w)\prod_{i=1}^np(X_i|w)\]
$p(w|X^n)$ is called \textbf{posterior distribution}.
Where $Z_n$ is the normalization constant to let the integral be 1 and called \textbf{marginal likelihood}.

\[Z_n=\int\varphi(w)\prod_{i=1}^np(X_i|w)dw\]
It is possible to consider that the posterior distribution extracts the information from the data $X^n$.
\textbf{Predictive distribution} $p^*(x)$ is the statistical model $p(x|w)$ given weight 
by the posterior distribution $p(w|X^n)$, that is, 

\[p^*(x)=\int p(x|w)p(w|X^n)dw\]
Assuming that the predictive distribution $p^*(x)$ infers the true distribution $q(x)$ is Bayesian inference.

\subsubsection{accuracy of Bayesian inference}\label{sec:accuracy_evaluation}
This subsubsection assesses objectively how closely the predictive distribution $p^*(x)$ 
approximates the true distribution $q(x)$.
Although a variety of evaluation indexes are sometimes used in practice depending on circumstances, 
this paper introduces the generalization error, 
the generalization loss, and the free energy as general evaluation indexes.

The \textbf{generalization error} is defined as the following.

\[K_n=\mathrm{KL}(q\|p^*)=\int q(x)\log\frac{q(x)}{p^*(x)}dx\]
Where $\mathrm{KL}(q\|p^*)$ is the \textbf{Kullback-Leibler divergence} 
from the true distribution $q(x)$ to the predictive distribution $p^*(x)$.
The following properties hold for the Kullback-Leibler divergence.
$P$ and $Q$ are any probability distribution in the following.

\begin{enumerate}
\item $\mathrm{KL}(P\|Q)\geq0\ ({}^{\forall}P,Q)$
\item\label{kl_prp:definiteness} $\mathrm{KL}(P\|Q)=0\Longleftrightarrow P=Q$
\end{enumerate}
Though the Kullback-Leibler divergence is also called \textbf{Kullback-Leibler distance}, 
neither the symmetry nor the triangle inequality that is the axiom of distance holds in general, that is, 
neither $\mathrm{KL}(P\|Q)=\mathrm{KL}(Q\|P)$ nor 
$\mathrm{KL}(P\|R)\leq\mathrm{KL}(P\|Q)+\mathrm{KL}(Q\|R)$ holds for any probability distribution $R$.
However, it is said that the smaller the generalization error is, 
the more the predictive distribution $p^*(x)$ approximates the true distribution $q(x)$ 
because the Kullback-Leibler divergence is $0$ if and only if two probability distributions are the same.

The generalization error can be decomposed into the following as the definition shows.

\begin{equation}\label{eq:decomp_Kn}
K_n=-\left(-\int q(x)\log q(x)dx\right)+\left(-\int q(x)\log p^*(x)dx\right)
\end{equation}
The left bracket of the right side of the equation \eqref{eq:decomp_Kn} is the constant in the inference 
because it is derived only from the true distribution $q(x)$.
This is the \textbf{entropy} being denoted as $S$ hereafter.

\[S=-\int q(x)\log q(x)dx\]
It is said that entropy indicates the randomness of the probability distribution.
For example, consider the probability mass function on $\{0,1,\ldots,N-1\}$.
Let $q_k(k=0,1,\ldots,N-1)$ be the following probability mass function.

\[q_k(i)=\left\{\begin{array}{lr}1/(k+1)&(i=0,\ldots,k)\\0&(\mathrm{otherwise})\end{array}\right.\]
Let $S_k$ be the entropy of $q_k$. $S_k$ is calculated as the following.

\[S_k=-\sum_{i=0}^{N-1}q_k(i)\log q_k(i)=\sum_{i=0}^k\frac{1}{k+1}\log(k+1)=\log(k+1)\]
This increases monotonically with $k$.
Thus the more biased the probability is, the less the entropy is, and vice versa.

The right bracket of the right side of the equation \eqref{eq:decomp_Kn} is called 
the \textbf{generalization loss} being denoted as $G_n$ hereafter.

\[G_n=-\int q(x)\log p^*(x)dx\]
Since the entropy $S$ is the constant in the inference as mentioned above, 
the less the generalization loss is, the less the generalization error is, and vice versa.
The fact enables using the generalization loss instead of the generalization error.

Finally, this paper introduces the \textbf{free energy}.
The free energy is defined as the following.

\[F_n=-\log Z_n=-\log\int\varphi(w)\prod_{i=1}^np(X_i|w)dw\]
In other words, the free energy $F_n$ is the negative log marginal likelihood.
The following equation holds between the free energy $F_n$ and the generalization loss $G_n$.

\[\EE[G_n]=\EE[F_{n+1}]-\EE[F_n]\]
The mean of the increment of the free energy coincides with the mean of the generalization loss.
It is said that the small generalization loss has a relation to the small free energy.

The following accuracy is gained for the evaluation indexes introduced above.
If the true distribution $q(x)$ is feasible by the statistical model $p(x|w)$, that is, 
there exists $w^*\in W$ such that $q(x)=p(x|w^*)$ for any $x$, 
the following asymptotic behaviors hold for the generalization loss and the free energy
\cite{AAnonregularLM,AAnonidentifiableLM}.

\[\EE[G_n]=S+\frac{\lambda}{n}-\frac{m-1}{n\log n}+o\left(\frac{1}{n\log n}\right)\]
\[\EE[F_n]=nS+\lambda\log n-(m-1)\log\log n+O(1)\]
where $\lambda>0,m\in\NN$ are the \textbf{real log canonical threshold} (\textbf{RLCT} for short) and 
the \textbf{multiplicity} respectively.
The real log canonical threshold is also called the \textbf{learning coefficient} in the learning theory field.
If the statistical model satisfied the regular condition 
where the posterior distribution approximates the normal distribution, 
$\lambda$ is equal to $d/2$ letting $d$ be the dimension of the parameter space\cite{AAnonregularLM}.
The statistical model that doesn't satisfy the regular condition is 
called the \textbf{singular model} on the other hand.
Deriving the RLCT and the multiplicity of a specific singular model is 
one of the main study objects in Bayesian learning theory 
because it is more difficult to derive the formula of RLCT.
The singular models whose RLCT is revealed are the following: 
Latent Dirichlet Allocation\cite{LDA}, Non-negative Matrix Factorization\cite{NNMF}, 
restricted Boltzmann Machine\cite{restricted_BM}, 
three-layered neural networks whose activation functions are 
identity function\cite{rrr} (also called reduced rank regression), 
hyperbolic tangent function\cite{vandermonde} and Swish function\cite{swish_tanaka,swish_baba}, 
and mixture models whose components are normal distribution\cite{vandermonde}, 
Poisson distribution\cite{pm_arxiv,pm_master} and multinomial distribution\cite{mm_arxiv}.
The Vandermonde matrix-type singularities are one of a few polynomials whose RLCT is revealed.

It is known that the following propositions hold for any statistical model.

\begin{prp}[parameter upper bound\cite{AGaSLT}]\label{thm:parameter_upper_bound}
Let $q(x),p(x|w),\varphi(w)$ be the true distribution, the statistical model, 
and the prior distribution respectively and $W\subset\RR^d$ be the parameter space.
If there exists an open set $U\subset W$ such that

\[\left\{w\in U\middle| \int q(x)\log\frac{q(x)}{p(x|w)}dx=0\ ,\varphi(w)>0\right\}\neq\varnothing,\]
the following inequality holds for the real log canonical threshold $\lambda$.

\[\lambda\leq\frac{d}{2}\]
$d/2$ is called the \textbf{parameter upper bound} of RLCT in this paper.
\end{prp}

\begin{prp}[Jeffreys prior\cite{AGaSLT}]\label{thm:Jeffreys_prior}
Let $q(x), p(x|w)$ be the true distribution and the statistical model$(w\in\RR^d)$.
The \textbf{Fisher information matrix} $I(w)$ is defined by the following.

\[I_{ij}(w)=\int\frac{\partial f(x,w)}{\partial w_i}\frac{\partial f(x,w)}{\partial w_j}p(x|w)dx\]
where

\[f(x,w)=\log\left(\frac{q(x)}{p(x|w)}\right).\]
The \textbf{Jeffrey prior} is defined by the following and is used as the prior distribution.

\[\varphi(w)=\left\{\begin{array}{lr}\frac{1}{Z}\sqrt{\det I(w)}&(w\in W)\\
0&(\mathrm{otherwise})\end{array}\right.\]
Where $Z$ is the normalization constant such that $\int\varphi(w)dw=1$.
The following inequality holds for the real log canonical threshold $\lambda$.

\[\lambda\geq\frac{d}{2}\]
\end{prp}

\subsubsection{zeta function and real log canonical threshold}\label{sec:zeta_func_and_rlct}
While the previous subsubsection \ref{sec:accuracy_evaluation} explains that 
the RLCT and multiplicity have a relation to the accuracy of Bayesian inference, 
this subsubsection and the next subsubsection \ref{sec:reso_and_blowup} explain how to derive them.
To derive the RLCT and the multiplicity, it is effective to consider the following \textbf{zeta function}.
($z\in\CC$)

\begin{equation}\label{eq:def_zeta}
\zeta(z)=\int_WK(w)^z\varphi(w)dw
\end{equation}
Where the integral region $W$ is the support of the prior distribution $\varphi(w)$ and $K(w)$ is 
the \textbf{mean error function} defined by the following.

\begin{equation}\label{eq:MEF}
K(w)=\int q(x)\log\frac{q(x)}{p(x|w)}dx
\end{equation}
Suppose that $K(w)$ is a real analytic function and $\varphi(w)$ is infinitely differentiable 
and its support $W$ is compact for the zeta function defined by the equation \eqref{eq:def_zeta}.
If there exists the parameter $w\in W$ such that $K(w)\varphi(w)>0$, 
this zeta function is regular on $\Re(z)>0$ and can be analytically continued uniquely 
on the entire complex plane as the meromorphic function\cite{RSDD}.
Furthermore, the analytically continued zeta function has the poles on the negative real number, 
and the maximum pole and its order coincide with the negative RLCT $-\lambda$ and 
the multiplicity $m$ respectively.
Because it is difficult to derive the mean error function of the statistical model analytically directly, 
it had better replace the mean error function with a relatively simple polynomial 
by using the following properties before deriving the RLCT.

\begin{prp}[property of inequality\cite{AGaSLT}]\label{eq:rlct_prp:inequality}%要確認
Let $i=1,2$.
The following zeta function is defined for the pair $(K_i(w),\varphi_i(w))$.

\[\zeta_i(z)=\int_W K_i(w)^z\varphi_i(w)dw\]
Let the maximum pole and its order of $\zeta_i$ be $-\lambda_i,m_i$ respectively.
If

\[K_1(w)\leq K_2(w)\ ,\ \varphi_1(w)\geq\varphi_2(w)\]
holds for any $w\in W$,

\[\lambda_1<\lambda_2\]
or

\[\lambda_1=\lambda_2\ ,\ m_1\geq m_2\]
holds.

\end{prp}
The property of inequality also holds for the inclusive relation of the integral region.
In fact, it is derived by letting $K_1(w)=K_2(w)$ and 
substituting $\varphi(w)\mathbb{I}[w\in W_i]$ for $\varphi_i(w)$ ($W_2\subset W_1$).
Where $\mathbb{I}[\ ]$ is the \textbf{indicator function}, that is, 
returns $1$ if the condition written in $[\ ]$ is satisfied otherwise $0$.

\begin{prp}[property of summation \cite{AGaSLT}]\label{eq:rlct_prp:sum}
Let the parameter $w$ be divided into $w=(w_a,w_b)$.
If $K(w),\varphi(w)$ are also divided into

\[K(w)=K_a(w_a)+K_b(w_b)\]
\[\varphi(w)=\varphi_a(w_a)\varphi_b(w_b),\]
the following equations hold for the RLCT and the multiplicity.

\[\lambda=\lambda_a+\lambda_b\]
\[m=m_a+m_b-1\]
\end{prp}
The following properties also hold.

\begin{prp}[property of product]\label{eq:rlct_prp:prd}
Let the parameter $w$ be divided into $w=(w_a,w_b)$.
If $K(w),\varphi(w)$ are also divided into

\[K(w)=K_a(w_a)K_b(w_b)\]
\[\varphi(w)=\varphi_a(w_a)\varphi_b(w_b),\]
the following equations hold for the RLCT and the multiplicity.

\[\lambda=\min(\lambda_a,\lambda_b)\]
\[m=\left\{\begin{array}{cc}
m_a&(\lambda=\lambda_a\neq\lambda_b)\\
m_b&(\lambda=\lambda_b\neq\lambda_a)\\
m_a+m_b&(\lambda=\lambda_a=\lambda_b)\end{array}\right.\]
\begin{proof}
Transform the zeta function into the following.

\[\zeta(z)=\int K(w)^z\varphi(w)dw=\left(\int K(w_a)^z\varphi_a(w_a)dw_a\right)
\left(\int K_b(w_b)^z\varphi_b(w_b)dw_b\right)=\zeta_a(z)\zeta_b(z)\]
Since the negative RLCT is the maximum pole of the zeta function, 
the RLCT $\lambda$ is the minimum of the RLCT of $\zeta_a$ and $\zeta_b$.
Since the multiplicity is its order, it is the sum of both if $\lambda_a=\lambda_b$, 
otherwise the multiplicity of the minimum one.

\end{proof}
\end{prp}

\begin{prp}[property of equality]\label{eq:rlct_prp:equality}
Define the following zeta function for $K(w),H(w)$.

\[\zeta_K(z)=\int K(w)^z\varphi(w)dw\]
\[\zeta_H(z)=\int H(w)^z\varphi(w)dw\]
If there exists $c_1,c_2>0$ such that 

\[c_1H(w)\leq K(w)\leq c_2H(w),\]
the RLCT and the multiplicity of $\zeta_K$ and $\zeta_H$ are equal respectively.

\begin{proof}
It is obvious by the property of inequality(Proposition \ref{eq:rlct_prp:inequality}).
\end{proof}
\end{prp}
The most useful proposition among the above is the property of equality(Proposition\ref{eq:rlct_prp:equality}).
Letting $K(w)$ a mean error function, 
if the polynomial (or analytic function) $H(w)$ whose RLCT is equal to $K(w)$ is ready, 
it is substantially sufficient to consider the RLCT and the multiplicity of $H(w)$.
The property of inequality (Proposition\ref{eq:rlct_prp:inequality}) 
is also useful to derive the upper bound of RLCT.

\subsubsection{resolution of singularities and blow-up}\label{sec:reso_and_blowup}

Though the previous subsubsection \ref{sec:zeta_func_and_rlct} explains how to 
replace the mean error function $K(w)$ with the polynomial $H(w)$ whose RLCT coincides with that of $K(w)$, 
it is difficult to calculate the zeta function \eqref{eq:def_zeta} 
for a polynomial or an analytic function $H(w)$.
In order to resolve this problem, the resolution of singularities theorem is used.

\begin{thm}[Atiyah's resolution of singularities theorem \cite{RSDD}]\label{thm:atiyah_reso_of_sing}
Let $f(x)$ be an analytic function defined on a neighborhood of the origin of $\RR^d$ and not identically zero.
There exist an open set $U$, a real analytic manifold $W$, and a proper analytic map $g:W\rightarrow U$ such that

\begin{enumerate}
\item $g:W\rightarrow U$ is a real analytic isomorphism of $W\backslash W_0$ and $U\backslash U_0$, 
where $U_0=f^{-1}(0)$ and $W_0=g^{-1}(U_0)$.

\item For an arbitrary $P\in W_0$, there is a local coordinate $w=(w_1,\ldots,w_d)$ whose origin is $P$ 
and that satisfies the following equations.

\begin{equation}\label{eq:nc_form1}f(g(w))=a(w)w_1^{k_1}w_2^{k_2}\cdots w_d^{k_d}\end{equation}
\begin{equation}\label{eq:nc_form2}|g'(w)|=|b(w)w_1^{h_1}w_2^{h_2}\cdots w_d^{h_d}|\end{equation}

Where $k_i,h_i(i=1,\ldots,d)$ are non-negative integers and 
$a(w),b(w)\neq0(^\forall w\in W)$ are analytic functions.

\end{enumerate}
\end{thm}
As the equation \eqref{eq:nc_form1} and \eqref{eq:nc_form2} in the theorem, 
the analytic function represented by an identically positive or negative analytic function $h(w)$ 
and non-negative integers $n_1,\ldots,n_d\in\NN_0$ is called \textbf{normal crossing}.

\[h(w)w_1^{n_1}w_2^{n_2}\cdots w_d^{n_d}\]
The set of normal crossing is denoted as $\calNC$ hereafter.
Using the word normal crossing, theorem \ref{thm:atiyah_reso_of_sing} insists on the following: 
for any real analytic function $f$ that is not identically $0$, 
there is a real analytic manifold $W$ and a proper analytic map $g$ such that 
for any point on $W_0$ there is a local coordinate whose origin is the point and 
$f\circ g$ and $|g'|$ can be represented as a normal crossing.

A normal crossing singularity makes the calculation of RLCT easier.
In fact, fixed a local coordinate, the integral of each variable is elementary 
and the poles of the zeta function on the local coordinate are

\[\left\{-\frac{h_1+1}{k_1},-\frac{h_2+1}{k_2},\ldots,-\frac{h_d+1}{k_d}\right\}.\]
Thus, considering all local coordinates, the RLCT and the multiplicity are

\[\lambda=\underset{u:\Loco}{\min}\left\{\underset{i=1,\ldots,d}{\min}\left(
\frac{h_i^{(u)}+1}{k_i^{(u)}}\right)\right\}\]
\[m=\underset{u:\Loco}{\max}\left|\left\{\ i\ \middle|\lambda=\frac{h_i^{(u)}+1}{k_i^{(u)}}\right\}\right|\]
respectively. Where $\underset{u:\Loco}{\min}$($\underset{u:\Loco}{\max}$) means 
minimizing(maximizing) through any local coordinate $u$ and 
$h_i^{(u)},k_i^{(u)}$ are negative integers on local coordinate $u$.

As mentioned above, constructing a real analytic manifold $W$ and a proper analytic map $g$ that 
realize the resolution of singularities in the context of theorem \ref{thm:atiyah_reso_of_sing} 
derives the RLCT and the multiplicity.
The next description is how to construct such $W$ and $g$.
In learning theory, the operation called blow-up is used in general\cite{AGaSLT}.
The definition of \textbf{blow-up} is the following.

\begin{dfn}[blow-up on Euclid space]
Let $r\in\NN$ satisfy $2\leq r\leq d$. The blow-up centered with $d-r$ dimensional non-singular algebraic set

\[V=\{x\in\RR^d|x_1=x_2=\cdots=x_r=0\}\]
is the subset of direct product $\RR^d\times\PP^{r-1}$ defined by the following.

\[B_V(\RR^d)=\overline{\{(x,(x_1:x_2:\cdots:x_r))|x\in\RR^d\backslash V\}}\]
Where $\overline{A}$ is closure of $A$.

\end{dfn}
Next, the blow-up defined above is considered from an algebraic perspective.
Letting $U_i=\{(x_1:x_2:\cdots:x_r|x_i\neq0)\}(i=1,\ldots,r)$, 
the projective space can be represented as $\PP^{r-1}=\bigcup_{i=1}^rU_i$.
For any element $(x,(x_1:x_2:\cdots:x_r))\in\RR^d\times U_i$, the equation

\[(x_1:\cdots:x_i:\cdots:x_r)=\left(\frac{x_1}{x_i}:\cdots:1:\cdots:\frac{x_r}{x_i}\right)\]
holds, so by using the notation $x^{(i)}_j=x_j/x_i(j=1,\ldots,r,j\neq i)$, 
the following variable transformation $g_i:\RR^d\rightarrow\RR^d$ is gained.

\begin{equation}\label{eq:blowup_var_trans}
g_i=\left\{\begin{array}{rcl}
x_1&=&x_1^{(i)}x_i\\
x_2&=&x_2^{(i)}x_i\\
&\vdots&\\
x_{i-1}&=&x_{i-1}^{(i)}x_i\\
x_i&=&x_i\\
x_{i+1}&=&x_{i+1}^{(i)}x_i\\
&\vdots&\\
x_r&=&x_r^{(i)}x_i\\
x_{r+1}&=&x_{r+1}\\
&\vdots&\\
x_d&=&x_d
\end{array}
\right.
\end{equation}
and the coordinate of $B_V(\RR^d)$ whose $x_i$ is not $0$ is isomorphic to $\RR^d$, 
which means that the point on the coordinate is identified with that of $\RR^d$.
The gained local coordinate is called \textbf{$x_i$-chart}.
The blow-up is gained by gluing the points of these $r$ local coordinates together that go to the same point.

Blowing-up can be also applied for a local coordinate gained by another blow-up.
That means blow-up constructs and grows a tree structure $\calT=(V,E)$ in a sense.
The idea that identifies the iteration of blow-up with the tree structure of graph theory 
is utilized by previous studies\cite{3esop,2esop}.
The vertex set $V$ of this tree $\calT$ is the blown-up local coordinates system or 
the polynomial set on their local coordinates and 
the edge $e\in E$ is considered to correspond to the variable transformation with blow-up.
How to construct the blow-up tree is specifically written in the following.
In the following, the vertex set $V$ is the polynomial set on local coordinates.

\begin{enumerate}
\item\label{initialize_tree} Let the polynomial subject to blow-up be $f$ and initialize $\calT=(V,E)$,
$V=\{f\},E=\varnothing$.
\item\label{add_tree} Let $h\in V\backslash\pi_1(E)$.
Where $\pi_1$ is the projection that extracts the first element, that is $\pi_1:(x,y)\mapsto x$.
If needed, blow up against $h$.
For example, blowing up centered with non-singular set $\{x\in\RR^d|x_1=x_2=\cdots=x_r=0\}$,
polynomial $h\circ g_i$ is gained on $x_i$ chart where $g_i$ is defined by \eqref{eq:blowup_var_trans}.
The gained polynomials are added to $V$, then the pairs of them and the original polynomial $h$ are added to $E$,
that is, $V$ and $E$ is renewed as $V\cup\{h\circ g_i|i=1,2,\ldots,r\}$ and 
$E\cup\{(h,h\circ g_i)|i=1,2,\ldots,r\}$ respectively.
\item The algorithm finishes and returns to \ref{add_tree}. The specific halting judgment differs by algorithms.
\end{enumerate}

The tree $\calT=(V,E)$ constructed by the above procedure is called a \textbf{blowup tree}.
In the context of tree components, the elements of $V$ are called \textbf{nodes}.
The density $|V|$ of a set of nodes is called the \textbf{size} of the tree.
If there exists an edge $(g,h)\in E$, 
$h$ is called a \textbf{child node} of $g$ and $g$ is done the \textbf{parent node} of $h$.
Also, node $g$ is a \textbf{descendant} of node $g'$ or node $g'$ is an \textbf{ancestor} of node $g$ 
if there is an edge sequence $(g_i, g_{i+1})_{i=0}^{n-1}\subset E$ such that $g_0=g,g_n=g'$.
For convenience, node $g$ itself is assumed to be both an ancestor and a descendant of node $g$.
The only node without a parent node is $f$, which is called a \textbf{root} node.
Conversely, a node with no children is called a \textbf{leaf} node.
The set of leaf nodes $V\backslash\pi_1(E)$ of a tree $\calT=(V,E)$ is denoted by $\calL(\calT)$ or simply $\calL$.
For each node $h$, the root node $f$ can be reached by traversing its edges, 
and the number of edges passed is denoted by $d(h)$, which is the \textbf{depth} of node $h$.
The depth of the root node $f$ is set as $d(f)=0$.
Let $\max d(h)$ be the depth of tree $\calT$ and denote it by $d(\calT)$.

In the following, we discuss restricting the type of nonsingular sets that 
are central to the blowup. For $d\in\NN,d\geq2$ for $X_d=\{1,2,\ldots,d\}$ 
We define the following nonsingular sets for subsets $T\subset X_d$ of concentration $2$ or more.

\[C_T=\left\{x\in\RR^d\middle|x_t=0(^\forall t\in T)\right\}\]
This family of nonsingular sets is also defined below.

\[\calC=\left\{C_T\middle|T\subset\{1,2,\ldots,d\},\ |T|\geq2\right\}\]
$\calC$ is called a family of \textbf{nonsingular sets around the coordinate axes}.
This study mainly considers blowups around elements of $\calC$.
The variable transformation by blow-up around the elements of $\calC$ is 
a linear transformation of the logarithm of the variables.
A concrete example is shown below.
Let $x\in\RR^d$ be a variable before variable transformation and $x'\in\RR^d$ be a variable 
after variable transformation.
Also, let $\log x=(\log x_1,\log x_2,\ldots,\log x_d)^\top\in\RR^d$.
The variable transformation $g_i$ expressed in equation \eqref{eq:blowup_var_trans} 
can be expressed as a matrix related to the exponents of the variables as follows.

\[\log x=\bmat{I_{i-1}&1_{i-1}&&\\
&1&&\\
&1_{r-i}&I_{r-i}&\\
&&&I_{d-r}}\log x'=\left(\prod_{j=1}^r R_{ji}\right)\log x'\]
where $R_{ji}\in\RR^{d\times d}$ is the basic transformation matrix with diagonal and 
$(j,i)$ components $1$ and other components $0$.
As described above, the variable transformation by blow-up around the elements of $\calC$ is 
a linear mapping between the logarithms of the variables.
The matrix can be expressed as a product of the basic transformation matrix $R_{ji}$.
This is a linear map between the logarithms of the variables.
This also holds in the case of a chain of blowups around the elements of $\calC$.
Let us consider the blowup tree $\calT=(V,E)$ that consists of blowups centered on elements of $\calC$.
Let $g,g'\in V$ be the nodes of the blowup tree and there exists a variable transformation $h$ such that 
$g'=g\circ h$.
In this case, if let $x=h(x')$, there exists a regular matrix $B\in\NN^{d\times d}$ such that 
$\log x = B\log x'$.
Thus, the matrix $B$ associated with the blowups centered with the elements of $\calC$ is called 
\textbf{blow matrix}.
The blow matrix $B$ is an element of the following set by its definition.

\begin{equation}\label{eq:def_set_of_basic_transform_matrix}
\calK=\left\{B\in\NN_0^{d\times d}\middle|^\exists(R_{i_k,j_k})_{k=1}^n\ s.t.\ 
B=\prod_{k=1}^nR_{i_k,j_k}\right\}
\end{equation}

$\calK$ is the set of matrices consisting of the product of the basic transformation matrices $R_{ij}$.
Therefore, $\det(B)=1$ for any $B\in\calK$.

Whether the blowup repeatedly results in a normal crossing at all local coordinates or not 
depends on the choice of the nonsingular set that is the center of the blowup and other factors.
This nonsingular set is not necessarily an element of $\calC$.
It is known that the singularity can be resolved by repeating appropriate blowups.

\begin{thm}[Hironaka's resolution of singularities theorem\cite{RSAV1,RSAV2,AGaSLT}]\label{thm:hironaka_reso_of_sing}
Let $f\in\RR[x_1,x_2,\ldots,x_d]$ be any polynomial.
Then the sequence of pairs of real algebraic sets $\{(V_i,W_i)\}_{i=0}^n$ 
satisfy the following conditions.

\begin{enumerate}
\item $V_i\subset W_i(i=1,2,\ldots,n)$
\item $V_0=\{x\in\RR^d|f(x)=0\},W_0=\RR^d$
\item $W_i(i=0,1,\ldots,n)$ are non-singular algebraic sets.
\item $V_n$ is represented by a normal crossing in any local coordinate of $W_n$.
\item For $i=1,2,\ldots,n$,$W_i=B_{C_{i-1}}(W_{i-1})$, that is, $W_i$ is a blowup centered on $C_i$ in $W_{i-1}$.
\item Let $\pi_i:W_i\rightarrow W_{i-1}$ be the projection defined on the blowup $B_{C_{i-1}}(W_{i-1})$.
Then $V_i$ is the total pullback of $\pi_i$, that is, $V_i=\pi^{-1}(V_{i-1})$.
\item Each blowup centered with $C_i$ is a non-singular real algebraic set 
contained in the critical point set of $f\circ\pi_1\circ\pi_2\circ\cdots\circ\pi_i$.
\end{enumerate}
\end{thm}
Theorem \ref{thm:hironaka_reso_of_sing} shows that the real log canonical threshold and multiplicity of 
the learning model can be obtained by a finite number of blowups.
However, it has not been clarified enough whether or not the resolution of singularities is 
possible by using a specific blowup algorithm when a specific polynomial is given.
This paper considers the real log canonical threshold through a blowup algorithm applicable to 
a polynomial called the sum-of-products polynomial.

In addition to the above-mentioned blowups,
there are other methods for resolving singularities called toric modification \cite{AAND,SCND} 
and weighted blowup \cite{weighted}.
For a study of real log canonical thresholds using toric modification, see \cite{AMEIBS}.
As for weighted blowups, 
the definitions and applications of our results are given in section \ref{sec:apply_for_weighted_blowup}.

\subsection{previous studies}

This section discusses previous studies. First, we define the sum-of-products polynomial, 
which is the subject of this paper, and the exclusive sum-of-products polynomial, 
which is the subject of previous studies, and then present the results of the previous studies.
Finally, the novelty and usefulness of this paper are reviewed.

\subsubsection{definition of sum-of-products polynomial}

Before introducing the previous studies, 
this section will give a definition of the sum-of-products polynomial, 
which is the subject of this paper.
It is known that the sum-of-products polynomial has a structure similar to 
the mean error function of a Boltzmann machine \cite{3esop}.
In general, in the study of real log canonical threshold in Bayesian learning theory, 
Though few mean error functions are equivalent to the sum-of-products polynomial, 
the coefficients of the polynomial are not important for the blowup of a nonsingular set $C_T\in\calC$ 
around a coordinate axis.
It is sufficient to consider the sum-of-products polynomial with coefficients fixed to $1$ for simplicity.

\begin{dfn}[sum-of-products polynomial]\label{dfn:sop_polynomial}
A real polynomial with $d$-variables $f$ is a \textbf{sum-of-products polynomial} 
if the coefficients of all terms are $1$. That is,

\[f(w)=\sum_{j=1}^n\prod_{i=1}^dw_i^{a_{ij}}\ 
\left(a_{ij}\in\NN_0,\ j\neq j'\Longrightarrow\sum_{i=1}^d|a_{ij}-a_{ij'}|\geq1\right) .\]
A \textbf{sop polynomial} is short for sum-of-products polynomial.

\end{dfn}

The set of sum-of-products polynomials with $d$-variables and $n$-terms is denoted by $\calSP^{d,n}$.
In particular, the set of non-negative sum-of-products polynomials with $d$-variables and $n$-terms
is denoted by $\calSP^{d,n}_{\geq0}=\{f\in\calSP^{d,n}|f\geq0\}$.
Since this paper assumes that the sum-of-products polynomial is used instead of the mean error function, 
the main object of consideration is the non-negative sum-of-products polynomial $\calSP^{d,n}_{\geq0}$.

\begin{xmp}
Examples of sum-of-products polynomials are given in the following.

\[f(w_1,w_2)=w_1^3w_2^5+w_1^6w_2+w_1^4w_2^2\in\calSP^{2,3}\]
\[g(w_1,w_2,w_3)=w_1^4w_2^3w_3^2+w_1^3w_2^2w_3^3+w_1^2w_2^3w_3^2+w_1^3w_2+w_3^2\in\calSP^{3,5}\]
\end{xmp}
A sum-of-products polynomial is characterized by the multi-indexes of each term.
Let $A=(a_{ij})\in\NN_0^{d\times n}$, a sum-of-products can be represented as

\[f(w)=1_n^\top\exp(A^\top\log w).\]
Where $1_n\in\RR^d$ is a vector whose components are all $1$s, 
and application of $\exp$ and $\log$ on a vector returns 
the vector of its components applied $\exp$ and $\log$ respectively.
That is, let $k\in\NN$ be the vector length, then $\exp:(a_1,a_2,\ldots,a_k)^\top\mapsto
(e^{a_1},e^{a_2},\ldots,e^{a_k})^\top$ and $\log:(a_1,a_2,\ldots,a_k)^\top\mapsto
(\log a_1,\log a_2,\ldots,\log a_k)^\top$.
$A\in\NN^{d\times n}$ is called a \textbf{multi-indexes matrix} and can be denoted by $A_f$ 
if the original sum-of-products polynomial $f$ needs to be clarified.
By its definition, let the multi-indexes matrix $A=\bmat{a_{*1}&a_{*2}&\cdots&a_{*n}}$, 
$a_{*j}\neq a_{*j'}$ holds for any different $j,j'$.
Consider applying an iterative blowup centered with a nonsingular set $C_T\in\calC$ to 
a sum-of-products polynomial $f$. Let $g$ be the variable transformation obtained by repeated blowups, 
and let $B$ be its blow matrix. Let $w\in\RR^d$ be the variable before variable transformation and 
$w'\in\RR^d$ be the variable after variable transformation. Then $\log w= B\log w'$ holds.
And the following polynomial is obtained.

\[f(g(w'))=1_n^\top\exp(A_f^\top B\log w')=1_n^\top\exp((B^\top A_f)^\top\log w')\]

Since $B^\top$ is a regular matrix, let $A=\bmat{a_{*1}&a_{*2}&\cdots&a_{*n}}\in\NN^{d\times n}$, 
then $B^\top a_{*j}\neq B^\top a_{*j'}$ for different $j\neq j'$.
Therefore,$f(g(w'))$ is also a product-sum polynomial.
From the above, the sum-of-products polynomial $\calSP^{d,n}$ is closed with respect to 
the operation of a blowup centered with a nonsingular set $C_T\in\calC$ around the coordinate axes.

Next, the exclusive sum-of-products polynomial which is the subject of previous studies 
is defined in the following.

\begin{dfn}[exclusive sum-of-products polynomial]\label{dfn:exclusive_sop_poly}
The real polynomial $f$ with $d$-variables is an \textbf{exclusive sum-of-products polynomial} 
if $f$ is a sum-of-products polynomial and each variable is contained in only one term. That is,

\[f(w)=\sum_{j=1}^n\prod_{i=1}^dw_i^{a_{ij}}\ 
\left(a_{ij}\in\NN_0,
\ j\neq j'\Longrightarrow\sum_{i=1}^da_{ij}a_{ij'}=0,\ \sum_{i=1}^d|a_{ij}-a_{ij'}|\geq1\right) .\]
\end{dfn}

\begin{xmp}
Examples of exclusive sum-of-products polynomials are given in the following.
\[f(w_1,w_2,w_3,w_4)=w_1^4w_2^3+w_3^5+w_4\]
\[f(w_1,w_2,w_3,w_4,w_5,w_6)=w_1^5w_2^4+w_3^6+w_4^3w_5^1+w_6^2\]
\end{xmp}

Because of the definition, $n-1\leq d$ holds for the exclusive sum-of-products polynomial.
Also, for a multi-indexes matrix $A\in\NN^{d\times n}$,
let $A=\bmat{a_{*1}&a_{*2}&\cdots&a_{*n}}$, $a_{*j}^\top a_{*j'}=0$ holds for different $j\neq j'$.

\subsubsection{RLCT and blowup algorithm for exclusive sop polynomial}\label{sec:ex_sop_rlct_alg}

The previous studies consider a special case of exclusive sum-of-products polynomials.
An exclusive sum-of-products polynomial with at most $m$ variables in each term is 
called a \textbf{$m$-multiple sum-of-products polynomial}.
The following results are obtained for the double sum-of-products polynomial.

\begin{thm}[RLCT of double exclusive sop polynomial \cite{2esop}]\label{thm:double_sop_rlct}
The double exclusive sop polynomial $K$ with $n$-terms is defined as follows.

\[K(x_1,x_2,\ldots,x_n,y_1,y_2,\ldots,y_n)=\sum_{i=1}^nx_i^{m_i}y_i^{n_i}\]
The zeta function of $K$ is defined as follows.

\[\zeta(z)=\iint_WK(x,y)^zdxdy\]
Then the maximum pole $-\lambda$ of $\zeta(z)$ is given by the appropriate choice of the integral domain $W$.

\begin{equation}\label{eq:double_sop_rlct}
\lambda=\sum_{i=1}^n\min\left(\frac{1}{m_i},\frac{1}{n_i}\right)
\end{equation}
\end{thm}

It is easy to show that the upper bound of the double exclusive sop polynomial is \eqref{eq:double_sop_rlct} 
by using the combination of the real log canonical threshold property 
(Propositions \ref{eq:rlct_prp:inequality},\ref{eq:rlct_prp:sum}, and \ref{eq:rlct_prp:prd}).
In the proof of Theorem \ref{thm:double_sop_rlct}, 
an appropriate choice of the integral domain $W$ of the zeta function and 
constituting blow-up in the sense of Theorem \ref{thm:hironaka_reso_of_sing} derive RLCT \cite{2esop}.

The following results are obtained for the triple sum-of-products polynomial.

\begin{thm}[RLCT and blowup algorithm for triple exclusive sop polynomial \cite{3esop}]\label{thm:triple_sop}
Two subsets of the triple exclusive sop polynomials $\mathcal{F},\mathcal{G}'$ and 
two algorithms are defined as follows.

\[\mathcal{F}=\left\{\sum_{i=1}^2x_i^{m_i}y_i^{n_i}z_i^{l_i}\middle|m_i,n_i,l_i\in\NN_0(^\forall i=1,2),
0\leq m_i,n_i,l_i\leq1(^\exists i=1,2),
\sum_{i=1}^2(m_i+n_i+l_i)\neq0\right\}\]
\[\mathcal{G}'=\left\{\sum_{i=1}^2x_i^{m_i}y_i^{n_i}z_i^{l_i}\middle|m_i,n_i,l_i\in\NN_0(^\forall i=1,2),
0\leq m_i,n_i,l_i\leq3(^\forall i=1,2),
1\in\{m_i,n_i,l_i,|i=1,2\}\right\}\]

\begin{alg}[\textbf{minimum degree selective blow-up}]
Let $K\in\mathcal{F}$.
If $K$ is a normal crossing, no more blow-up is performed in the branch.
If $K$ is not a normal crossing, we first factor out common factors.
Then there exist monomials $p_1,p_2$ which have no common factors and monomial $c$, and $K=c(p_1+p_2)$ holds.
Let $s_1,s_2$ be the variables with the minimum degree in $p_1$ and $p_2$ respectively.
Then blow-up centered with $\{s_1=s_2=0\}$ for $K$.
Repeat the above procedure for the polynomials $K_1$ and $K_2$ obtained by blow-up, respectively.

The pseudo-code for this algorithm is shown in \ref{alg:argmin_selected_blowup}.

%\stepcounter{thm}
\begin{algorithm}[htbp]
\caption{minimum degree selective blow-up}
\label{alg:argmin_selected_blowup}
\begin{algorithmic}[1]
\REQUIRE $K\in\mathcal{F}$
\ENSURE $\calT  =(V,E)$\COMMENT{$\calT  $ : tree}
%\FUNCTION{argmin\_selected\_blowup}{$K$}
\STATE $Q \gets [\ ]$\COMMENT{$Q$ : queue}
\STATE $V \gets \varnothing$\COMMENT{$V$ : polynomial (as vertex) set}
\STATE $E \gets \varnothing$\COMMENT{$E$ : edge (pair of vertex) set}
\STATE $Q.\mathrm{enqueue}(K)$
\STATE $V.\mathrm{add}(K)$\COMMENT{$V \gets V\cup\{K\}$}
\WHILE{$Q \neq [\ ]$}
	\STATE $f \gets Q.\mathrm{dequeue}()$
	\IF{$f \notin\calNC $}
			\STATE $c,p_1,p_2\gets f.\mathrm{factorize}()$
			\COMMENT{$c$ : common factor, $p_1,p_2$ : coprime, $f=c(p_1+p_2)$}
			\FOR{$i \in \{1,2\}$}
				\STATE $s_i \gets p_i.\mathrm{argmindeg}()$
				\COMMENT{$s_i$ : argmin of the degree of variables of $p_i$}
			\ENDFOR
			\COMMENT{blow-up with center $\{s_1=s_2=0\}$}
			\FOR{$i \in \{1,2\}$}
				\STATE $f_i \gets f(s_i=s_1s_2)$
				\COMMENT{substitute $s_1s_2$ for $s_i$}
				\STATE $Q.\mathrm{enqueue}(f_i)$
				\STATE $V.\mathrm{add}(f_i)$
				\STATE $E.\mathrm{add}((f,f_i))$
			\ENDFOR
	\ENDIF
\ENDWHILE
\STATE $\calT\gets (V,E)$
\STATE \textbf{Return} $\calT$
%\ENDFUNCTION
\end{algorithmic}
\end{algorithm}

\end{alg}

\begin{alg}[\textbf{maximum degree selective blow-up}]
Let $K\in\mathcal{G}'$.
If $K$ is a normal crossing, no more blow-up is performed in the branch.
If $K$ is not a normal crossing, we first factor out common factors.
Then there exist monomials $p_1,p_2$ which have no common factors and monomial $c$, and $K=c(p_1+p_2)$ holds.
Let $s_1,s_2$ be the variables with the maximum degree in $p_1$ and $p_2$ respectively.
Then blow-up centered with $\{s_1=s_2=0\}$ for $K$.
Repeat the above procedure for the polynomials $K_1$ and $K_2$ obtained by blow-up, respectively.

The pseudo-code for this algorithm is shown in \ref{alg:argmax_selected_blowup}.

%\stepcounter{thm}
\begin{algorithm}[htbp]
\caption{maximum degree selective blow-up}
\label{alg:argmax_selected_blowup}
\begin{algorithmic}[1]
\REQUIRE $K\in\mathcal{G}'\backslash\mathcal{F}$
\ENSURE $\calT  =(V,E)$
%\Function{argmax\_selected\_blowup}{$K$}
\STATE $Q \gets [\ ]$
\STATE $V \gets \varnothing$
\STATE $E \gets \varnothing$
\STATE $Q.\mathrm{enqueue}(K)$
\STATE $V.\mathrm{add}(K)$
\WHILE{$Q \neq [\ ]$}
	\STATE $f \gets Q.\mathrm{dequeue}()$
	\IF{$f \notin\calNC $}
			\STATE $c,p_1,p_2\gets f.\mathrm{factorize}()$
			\FOR{$i \in \{1,2\}$}
				\STATE $s_i \gets p_i.\mathrm{argmaxdeg}()$
				\COMMENT{$s_i$ : argmax of the degree of variables of $p_i$}
			\ENDFOR
			\FOR{$i \in \{1,2\}$}
				\STATE $f_i \gets f(s_i=s_1s_2)$
				\STATE $Q.\mathrm{enqueue}(f_i)$
				\STATE $V.\mathrm{add}(f_i)$
				\STATE $E.\mathrm{add}((f,f_i))$
			\ENDFOR
	\ENDIF
\ENDWHILE
\STATE $\calT   \gets (V,E)$
\STATE \textbf{Return} $\calT  $
%\EndFunction
\end{algorithmic}
\end{algorithm}

\end{alg}
In the pseudo-codes \ref{alg:argmin_selected_blowup} and \ref{alg:argmax_selected_blowup} of the algorithm,
$Q$ is a queue, that is, a first-in first-out queue data structure.
Let $f$ be a polynomial, $f$.factorize() denotes a method to factor out the common factors of $f$.
The return of the method are monomials $c,p_1$, and $p_2$.
Where $c$ is their common factor and $p_1$ and $p_2$ are prime to each other.
For monomial $p_i$, $p_i$.argmindeg() and $p_i$.argmaxdeg() are methods that return 
the variable $s_i$ with the minimum and maximum degree among the variables in $p_i$ respectively.

The minimum degree selective blow-up algorithm \ref{alg:argmin_selected_blowup} for any $K\in\mathcal{F}$ stops 
and the maximum pole $-\lambda$ of the zeta function $\zeta(z)=\int_{W_K}K(w)^zdw$ is obtained by

\[\lambda=1.\]
where the integral domain is $W_K=\{K(w)\geq0|w\in W\}$ and $W$ is a sufficiently large compact set.
Also the maximum degree selective blow-up algorithm \ref{alg:argmax_selected_blowup} for 
$K\in\mathcal{G'}\backslash\mathcal{F}$ stops and the maximum pole $-\lambda$ of 
the zeta function $\zeta(z)=\int_{W_K}K^z$ is obtained by

\[\lambda=\sum_{i=1}^2\min\left(\frac{1}{m_i},\frac{1}{n_i},\frac{1}{l_i}\right).\]
\end{thm}

Theorem \ref{thm:triple_sop} proposes two blow-up algorithms and shows the existence of certain 
triple exclusive sop polynomials such that each algorithm stops.
On the other hand, there is a counter-example for which the minimum degree selective 
algorithm \ref{alg:argmin_selected_blowup} does not stop \cite{3esop}.
And it has been confirmed that the maximum degree selective algorithm \ref{alg:argmax_selected_blowup} 
does stop even when the degree limit rises, but halting for any multiple exclusive sop polynomial 
has not been proved.

\subsubsection{novelty and utility of this study}

Compared with the previous studies, this study extends the scope of study from exclusive sop polynomials 
to general sop polynomials.
New propositions and theorems are shown for the blowup algorithm and 
the real log canonical threshold with the extension.

As for the blowup algorithm, as a novel point, 
this paper modifies the maximum degree selective blowup algorithm \ref{alg:argmax_selected_blowup} 
proposed in the previous study and proposes an algorithm for the resolution of the singularities 
in the sense of Theorem \ref{thm:hironaka_reso_of_sing} for nonnegative sop binomials, 
and proves halting (Section \ref{sec:blowup_algorithm_2sop}).
For general sop polynomials, by using the proposed algorithm as a subroutine and relaxing the conditions,
this paper proposes an algorithm that can derive an upper bound on the real log canonical threshold, 
and proves its halting (Section \ref{sec:dn_sop_alg}).

As for the RLCT, as a novel point, 
we obtain the exact value of the RCLT for nonnegative sop binomials (Section \ref{sec:d2_case}) 
and the upper bound of the RLCT for general nonnegative sop polynomials (Section \ref{sec:dn_case}) 
by using our proposed blow-up algorithm.
As a useful point, we confirm that the upper bound of the RLCT of general sop polynomials does not 
depend on the coefficients of the polynomials, that is, it can be applied to general polynomials.
If we can find a polynomial whose RLCT is equivalent to that of the mean error function of learning model 
whose RLCT has not known yet, we can obtain an upper bound of the RLCT by solving 
a linear programming problem on the multi-indexes of the polynomial.
We also show that the upper bound is tighter than the parameter upper bound, 
which is as widely applicable as this upper bound (Section \ref{sec:simplex_upper_bound}).
We also confirm that the optimal weights of the weighted blow-up can be computed from the optimal solution 
of the linear programming problem for multi-indexes of polynomials (Section \ref{sec:apply_for_weighted_blowup}).

\section{Proposed algorithm and main theorem}\label{sec:proposed_algorithm_and_main_theorem}

In this section, this paper describes the proposed algorithm and the main theorem. 
The main theorem is divided into a blow-up algorithm and a RLCT.
First, the concepts necessary for the statement of the main theorem are defined.
Let the set of $d$-variable monomials whose multi-indexes are a real number greater than $-1$ be denoted by

\[\calM^d=\left\{\prod_{i=1}^dw_i^{s_i-1}\middle|s_i>0(i=1,2,\ldots,d)\right\}.\]
If the analytic function $f(w)\in\RR^d\rightarrow\RR$ is represented by

\[f(w)=h(w)w_1^{k_1}w_2^{k_2}\cdots w_d^{k_d}\ (h(w)\neq0\ \mathrm{on}\ N)\]
where $h(w)$ is an analytic function that is not zero in a neighborhood $N\subset\RR^d$ of the origin 
and  $k_1,k_2,\ldots,k_d\in\NN_0$ are non-negative integers, $f$ is called a 
\textbf{local normal crossing near the origin}.
In the following, this paper will mainly consider the cases near the origin, 
so that when we just speak of local normal crossing, we mean local normal intersections near the origin, 
unless otherwise described.
The set of local normal crossings is denoted by $\calLNC$.

Next, we propose three blowup algorithms (Algorithms \ref{alg:between_two_vars_with_jacobian},
\ref{alg:between_terms}, and \ref{alg:local_nc_blowup}).
Algorithm \ref{alg:between_two_vars_with_jacobian} (blowup between variables with Jacobian) 
is a subroutine of Algorithm \ref{alg:between_terms} (blowup between terms), 
and Algorithm \ref{alg:between_terms} is a subroutine of Algorithm \ref{alg:local_nc_blowup} 
(local normal crossing blowup).

\begin{alg}[\textbf{blow-up between variables with Jacobian}]
Let $f$ be a non-negative bivariate binomial and $g\in\calM^2$.
If $f$ is a normal crossing, no more blow-up is performed in the branch.
If $f$ is not a normal crossing, blowup centered with the origin $\{(0,0)\in\RR^2|w_1=w_2=0\}$.
Note that, for $g$, the term due to the Jacobian with the blowup is also multiplied.
The above procedure is repeated for $(f_1,g_1)$ and $(f_2,g_2)$ obtained by blowup respectively.

The pseudo-code of this algorithm is shown in \ref{alg:between_two_vars_with_jacobian}.

%\stepcounter{thm}
\begin{algorithm}[htbp]
\caption{blow-up between variables with Jacobian}
\label{alg:between_two_vars_with_jacobian}
\begin{algorithmic}[1]
\REQUIRE $f(w_1,w_2)\in\calSP ^{2,2}_{\geq0},\ g(w_1,w_2)\in\calM^2$
\ENSURE $\calT  =(V,E)$
%\Function{blowup\_between\_two\_vars\_with\_jacobian}{$f,g,w_1,w_2$}
\STATE $Q \gets [\ ]$
\STATE $V \gets \varnothing$
\STATE $E \gets \varnothing$
\STATE $Q.\mathrm{enqueue}((f,g))$
\STATE $V.\mathrm{add}((f,g))$
\WHILE{$Q \neq [\ ]$}
	\STATE $(h,k) \gets Q.\mathrm{dequeue}()$
	\IF{$h \notin\calNC $}
		\FOR{$i \in \{1,2\}$}
			\STATE $h_i \gets h(w_i=w_1w_2)$
			\STATE $k_i \gets w_{3-i}k(w_i=w_1w_2)$
			\COMMENT{$w_{3-i}$ : term of Jacobian}
			\STATE $Q.\mathrm{enqueue}((h_i,k_i))$
			\STATE $V.\mathrm{add}((h_i,k_i))$
			\STATE $E.\mathrm{add}(((h,k),(h_i,k_i)))$
		\ENDFOR
	\ENDIF
\ENDWHILE
\STATE $\calT\gets(V,E)$
\STATE \textbf{Return} $\calT$
%\EndFunction
\end{algorithmic}
\end{algorithm}

\end{alg}

\begin{alg}[\textbf{blow-up between terms}]
Let $f$ be a non-negative $d$-variable binomial and $g\in\calM^d$.
If $f$ is a normal crossing, no more blow-up is performed in the branch.
If not, we factor out the common factors.
Then there exist monomials $c,p_1$ and $p_2$ where $c$ is the common factor and $p_1,p_2$ are prime each other,
and $f=c(p_1+p_2)$ holds.
Let the variables whose degrees are maximum in $p_1,p_2$ be $s_1,s_2$ respectively.
Then run blow-up between variables with Jacobian \ref{alg:between_two_vars_with_jacobian} 
for polynomial $f$, monomial $g$, and variables $s_1,s_2$.
The above procedure is repeated for each pair of polynomials $(f',g')$ in each local coordinate obtained by 
the blowup between variables with Jacobian.

The pseudo-code of the algorithm is shown in \ref{alg:between_terms}.

%\stepcounter{thm}
\begin{algorithm}[htbp]
\caption{blow-up between terms}
\label{alg:between_terms}
\begin{algorithmic}[1]
\REQUIRE $f\in\calSP^{d,2}_{\geq0},g\in\calM^d$
\ENSURE $\calT=(V,E)$
%\Function{blowup\_between\_terms}{$f,g$}
\STATE $Q \gets [\ ]$
\STATE $V \gets \varnothing$
\STATE $E \gets \varnothing$
\STATE $Q.\mathrm{enqueue}((f,g))$
\STATE $V.\mathrm{add}((f,g))$
\WHILE{$Q \neq [\ ]$}
	\STATE $(h,k) \gets Q.\mathrm{dequeue}()$
	\IF{$h \notin\calNC $}
			\STATE $c,p_1,p_2\gets h.\mathrm{factorize}()$
			\FOR{$i \in \{1,2\}$}
				\STATE $s_i \gets p_i.\mathrm{argmaxdeg}()$
			\ENDFOR
			\STATE $\calT' \gets \mathrm{blowup\_between\_two\_vars\_with\_jacobian}(h,k,s_1,s_2)$
			\FOR{$(h',k') \in V(\calT')$}
				\IF{$(h',k') \in \calL(\calT')$}
					\STATE $Q.\mathrm{enqueue}((h',k'))$
				\ENDIF
				\STATE $V.\mathrm{add}((h',k'))$
			\ENDFOR
			\FOR{$e \in E(\calT')$}
				\STATE $E.\mathrm{add}(e)$
			\ENDFOR
	\ENDIF
\ENDWHILE
\STATE $\calT   \gets (V,E)$
\STATE \textbf{Return} $\calT$
%\EndFunction
\end{algorithmic}
\end{algorithm}

\end{alg}

\begin{alg}[\textbf{local normal crossing blow-up}]
Let $f$ be a nonnegative sop polynomial with $d$-variables and $n$-terms and $g\in\calM^d$.
Let $n_1$ and $n_2$ be counter variables that indicate up to which term the blowup is done.
Initially, $n_1=1$ and $n_2=2$.
If $f$ is a local normal crossing, no more blowup is performed in the branch.
If not, expand $f$ as a sum of monomials. That is, let monomials be denoted by $p_1,\ldots,p_n$,
$f=p_1+\cdots+p_n$ holds.
Then focus on $(p_{n_1}+p_{n_2},g)$ among them and run blowup between terms \ref{alg:between_terms}.
Let $(p'_{n_1}+p'_{n_2},g')$ be the pair at any local coordinate obtained 
by the blowup between terms \ref{alg:between_terms}.
Next, update the counter variables at any local coordinate.
If $p'_{n_1}\leq p'_{n_2}$, update the counter $n_2=\max(n_1,n_2)+1$ in the local coordinate.
If not, update the counter $n_1=\max(n_1,n_2)+1$ in the local coordinate.
Let $f'$ be the result of the same variable transformation as 
$(p_{n_1}+p_{n_2},g)\rightarrow(p'_{n_1}+p'_{n_2},g')$ for $f$.
Then $f'=p'_1+\cdots+p'_n$ holds.
Repeat the above procedure for any $(f',g')$ respectively.

The pseudo-code of this algorithm is shown in \ref{alg:local_nc_blowup}.

%\stepcounter{thm}
\begin{algorithm}[htbp]
\caption{local normal crossing blow-up}
\label{alg:local_nc_blowup}
\begin{algorithmic}[1]
\REQUIRE $f\in\calSP^{d,n}_{\geq0},g\in\calM^d$
\ENSURE $\calT=(V,E)$
%\Function{local\_normal\_crossing\_blowup}{$f,g$}
\STATE $Q \gets [\ ]$
\STATE $V \gets \varnothing$
\STATE $E \gets \varnothing$
\STATE $n_1 \gets 1$
\STATE $n_2 \gets 2$
\STATE $Q.\mathrm{enqueue}((f,g,n_1,n_2))$
\STATE $V.\mathrm{add}((f,g))$
\WHILE{$Q \neq [\ ]$}
	\STATE $(h,k,n'_1,n'_2) \gets Q.\mathrm{dequeue}()$
	\IF{$h \notin\calLNC $}
			\STATE $p_1,p_2,\ldots,p_n\gets h.\mathrm{expand}()$
			\COMMENT{$p_1,p_2,\ldots,p_n$ : monomial, $h=p_1+p_2+\cdots+p_n$}
			\STATE $\calT' \gets \mathrm{blowup\_between\_terms}(p_{n'_1}+p_{n'_2},k)$
			\FOR{$(h',k') \in V(\calT')$}
				\STATE $p'_{n'_1},p'_{n'_2} \gets h'.\mathrm{expand}()$
				\COMMENT{$p'_1,p'_2$ : monomial, $h'=p'_{n'_1}+p'_{n'_2}$}
				\STATE $h'' \gets h'.\mathrm{decompose}(p_{n'_1}+p_{n'_2})$
				\COMMENT{$h'':\RR^d\rightarrow\RR^d$, $h'=(p_{n'_1}+p_{n'_2})\circ h''$}
				\IF{$(h',k') \in \calL(\calT')$}
					\IF{$p'_{n'_1}\leq p'_{n'_2}$}
						\STATE $Q.\mathrm{enqueue}((h\circ h'',k',n'_1,\max\{n'_1,n'_2\}+1))$
					\ELSE
						\STATE $Q.\mathrm{enqueue}((h\circ h'',k',\max\{n'_1,n'_2\}+1,n'_2))$
					\ENDIF
				\ENDIF
				\STATE $V.\mathrm{add}((h\circ h'',k'))$
			\ENDFOR
			\FOR{$e\in E(\calT')$}
				\STATE $E.\mathrm{add}(e)$
			\ENDFOR
	\ENDIF
\ENDWHILE
\STATE $\calT   \gets (V,E)$
\STATE \textbf{Return} $\calT$
%\EndFunction
\end{algorithmic}
\end{algorithm}

\end{alg}
In the pseudo-code \ref{alg:local_nc_blowup} of the proposed algorithm,
$h$.expand() is a method to decompose $h$ into a sum of monomials, where $h$ is a polynomial.
The returns $p_1,\ldots,p_n$ are monomials and $h=p_1+p_2+\cdots+p_n$ holds.
Also, $h'$.decompose(${p'}_{n_1}+{p'}_{n_2}$) is a method to decompose into the composition of the function.
There exist $h'':\RR^d\rightarrow\RR^d$ such that $h'=(p_{n'_1}+p_{n'_2})\circ h''$,
where $(h',k')$ is the leaf of the blowup tree $\calT'$ whose root node is $(p_{n'_1}+p_{n'_2},k)$.
The method returns $h''$.

The following theorem holds for the proposed algorithms 
\ref{alg:between_two_vars_with_jacobian},\ref{alg:between_terms}, and \ref{alg:local_nc_blowup}.

\begin{main_theorem}[the property of the proposed algorithm]\label{thm:main_theorem_blowup_algorithm}
Let $f(w)\in\calSP^{d,2}_{\geq0}$ be a non-negative sop binomial and 
$g(w)=w^{s-1_d}\in\calM^d$ be a monomial$(s\in\RR^d_{\geq0})$.
If blow-up between terms algorithm (Algorithm \ref{alg:between_terms}) is applied to $(f,g)$, 
the algorithm returns a finite blow-up tree $\calT=(V,E)$ such that 
the binomial is represented as a normal crossing for all leaves $(h,k)\in\calL(\calT)$.

Also, let $f(w)\in\calSP^{d,n}_{\geq0}$ be a non-negative sop polynomial and 
apply the local normal crossing blowup algorithm (Algorithm \ref{alg:local_nc_blowup}) for $(f,g)$, 
the algorithm returns a finite blow-up tree $\calT=(V,E)$ such that 
the polynomial is represented as a local normal crossing for all leaves $(h,k)\in\calL(\calT)$.
\end{main_theorem}

Also, the following holds for the RLCT of the sop polynomial.

\begin{main_theorem}[RLCT of sop polynomial]\label{thm:main_theorem_rlct}
Let $f(w)\in\calSP^{d,n}_{\geq0}$ be a non-negative sop polynomial and 
$g(w)=w^{s-1_d}\in\calM^d$ be a monomial $(s\in\RR^d_{\geq0})$.
Let $A_f=(a_{hj})\in\NN_0^{d\times n}$ be the multi-indexes matrix of $f$ and 
let $W$ be a sufficiently large compact set.
For the largest pole $-\lambda$ of the zeta function $\zeta(z)$ defined below

\[\zeta(z)=\int_Wf(w)^z|g(w)|dw,\]
the following holds.
\begin{equation}\label{eq:main_theorem}
\frac{1}{\lambda}\geq\underset{\alpha\in S_{d-1}}{\max}\left\{\underset{j=1,2,\ldots,n}{\min}
\left(\sum_{h=1}^d\alpha_h\frac{a_{hj}}{s_h}\right)\right\}
\end{equation}
where $S_{d-1}$ is \textbf{$d-1$-dimensional simplex} and represented as follows.

\[S_{d-1}=\left\{\alpha\in\RR^d\middle|\alpha_h\geq0(^\forall h=1,2,\ldots,d),\ \sum_{h=1}^d\alpha_h=1\right\}\]
For the case $n=2$, the equation of \eqref{eq:main_theorem} holds.
In particular, for $h=1,2,\ldots,d$, $\nu_h,\mu_h$ related to multi-indexes is defined below,

\[\nu_h=\frac{\max(a_{h1},a_{h2})}{s_h},\ \mu_h=\frac{\min(a_{h1},a_{h2})}{s_h}\]
and the index set $H$ is defined as follows.

\[H=\left\{(h,h')\in\{1,2,\ldots,d\}^2\middle|h\neq h',\nu_h-\mu_h,\nu_{h'}-\mu_{h'}>0\right\}\]
Then the following equation holds.

\[\frac{1}{\lambda}=\max\left\{
\underset{h}{\max}\left(\mu_h\right),\ 
\underset{(h,h')\in H}{\max}\left(
\frac{\nu_h\nu_{h'}-\mu_h\mu_{h'}}{\nu_h+\nu_{h'}-\mu_h-\mu_{h'}}\right)\right\}\]
\end{main_theorem}

\section{Proof of main theorem}\label{sec:proof_of_main_theorem}
In this section, this paper shows the proof of the main theorem.
We prove the theorem in the order of bivariate sop binomial, sop binomial, and general sop polynomial 
and for each of them, we show the exact value or upper bound of the RLCT and that the blowup algorithm halts.

\subsection{case of bivariate sop binomial}\label{sec:22_case}
Before the proof, we review the properties of the bivariate sop binomial $\calSP ^{2,2}_{\geq0}$.

\begin{prp}\label{prp:deg_is_even_22}
Let $f\in\calSP ^{2,2}_{\geq0}$. Let the multi-indexes matrix $A_f=\bmat{p&q\\r&s}\in\NN_0^{2\times2}$, 
then $p,q,r,s\in 2\NN_0$.

\begin{proof}
Let $f(w_1,w_2)=w_1^pw_2^r+w_1^qw_2^s$.
From the symmetry, it is sufficient to show that $p\in2\NN_0+1$ contradicts the assumption.
When $q=p$, $f(-1,1)=-2<0$ holds, which contradicts $f\geq0$.
When $p>q$, $f(-M,1)<0$ holds for sufficiently large number $M>0$, which contradicts $f\geq0$.
When $p<q$, $f(-\varepsilon,1)<0$ holds for sufficiently small numbers $1>\varepsilon>0$, 
which contradicts $f\geq0$. Therefore $p\in2\NN_0$.

\end{proof}
\end{prp}
This property makes it easier to know whether or not 
a given bivariate sop binomial $f\in\calSP^{2,2}_{\geq0}$ is a normal crossing.

\subsubsection{blowup algorithm for bivariate sop binomial}
Next, we show that the blowup algorithm for the bivariate sop binomial $\calSP^{2,2}_{\geq0}$ halts.

\begin{prp}\label{prp:alg_halt}
Applying the following algorithm \ref{alg:between_two_vars} to any 
$^\forall K(w_1,w_2)\in\calSP ^{2,2}_{\geq0}$, the algorithm halts.

\begin{alg}[\textbf{blow-up between variables}]
Let $f$ be a non-negative bivariate sop binomial.
If $f$ is a normal crossing, no more blowup is performed in the branch.
If not, blowup centered with the origin $\{(0,0)\in\RR^2|w_1=w_2=0\}$.
Repeat the above procedure for each of $f_1$ and $f_2$ obtained by the blowup above.

The pseudo-code of this algorithm is shown in \ref{alg:between_two_vars}.

%\stepcounter{thm}
\begin{algorithm}[htbp]
\caption{blow-up between variables}
\label{alg:between_two_vars}
\begin{algorithmic}[1]
\REQUIRE $K\in\calSP ^{2,2}_{\geq0}$
\ENSURE $\calT  =(V,E)$
%\Function{blowup\_between\_two\_vars}{$K,w_1,w_2$}
\STATE $Q \gets [\ ]$
\STATE $V \gets \varnothing$
\STATE $E \gets \varnothing$
\STATE $Q.\mathrm{enqueue}(K)$
\STATE $V.\mathrm{add}(K)$
\WHILE{$Q \neq [\ ]$}
	\STATE $f \gets Q.\mathrm{dequeue}()$
	\IF{$f \notin\calNC $}
		\FOR{$i \in \{1,2\}$}
			\STATE $f_i \gets f(w_i=w_1w_2)$
			\STATE $Q.\mathrm{enqueue}(f_i)$
			\STATE $V.\mathrm{add}(f_i)$
			\STATE $E.\mathrm{add}((f,f_i))$
		\ENDFOR
	\ENDIF
\ENDWHILE
\STATE $\calT   \gets (V,E)$
\STATE \textbf{Return} $\calT  $
%\EndFunction
\end{algorithmic}
\end{algorithm}

\end{alg}

\begin{proof}
Let $f\in\calSP ^{2,2}_{\geq0}$ and $A_f=\bmat{p&q\\r&s}\in\NN_0^{2\times 2}$.
Define $\calD :\calSP ^{2,2}_{\geq0}\rightarrow\ZZ$ as $\calD :f\mapsto(p-q)(r-s)$.

First, we show that $\calD (f)\geq0\Longleftrightarrow f\in\calNC$.
Suppose $\calD (f)\geq0$. If $p-q,r-s\geq0$, 
$f$ can be represented as $f(w_1,w_2)=w_1^qw_2^s(1+w_1^{p-q}w_2^{r-s})$.
Since $p-q$ and $r-s$ are even by Proposition \ref{prp:deg_is_even_22}, $f\in\calNC$ holds.
The case of $p-q$ and $r-s\leq 0$ can be shown in the same way.
Next, suppose the inverse $\calD (f)<0$.
If $p-q>0>r-s$, then $f(w_1,w_2)=w_1^qw_2^r(w_1^{p-q}+w_2^{s-r})$, which shows $f\notin\calNC$.
The same discussion holds for $r-s>0>p-q$.
From the above, it follows that $\calD (f)\geq0\Longleftrightarrow f\in\calNC$.

Next, we show that $\calD$ is monotonically increasing with respect to blow-up.
Since $\calD(f)\geq0$ is equivalent to $f\in\calNC$, no more blowup is needed for $f$.
Therefore, we suppose $\calD (f)<0$.
Since $f_1(w_1,w_2)=f(w_1w_2,w_2)=w_1^{p+r}w_2^r+w_1^{q+s}w_2^s$ holds, 
we get $\calD (f_1)=(p+r-q-s)(r-s)$. We get $\calD (f_2)=(p-q)(r+p-s-q)$ in the same way.
We also get $\calD (f_1)-\calD (f)=(r-s)^2>0$ and 
$\calD (f_2)-\calD (f)=(p-q)^2>0$ by calculating the difference.
Thus, since $\calD$ is monotonically increasing with respect to blowup and discrete, 
for any polynomial $h$ added to $Q$ after finite times of blowup, 
$\calD (h)\geq0$ holds, that is, $h\in\calNC$.
Therefore for any $K\in\calSP ^{2,2}_{\geq0}$, the algorithm \ref{alg:between_two_vars} halts.

\end{proof}
\end{prp}
Next, we describe a property of $\calD$, which is necessary for the discussion of RLCTs below.
Suppose $f\in\calSP ^{2,2}_{\geq0}$. If $\calD (f)<0$, we get $\calD (f_1)\calD (f_2)\leq0$, 
and if the equality holds in particular, then $\calD (f_1)=\calD (f_2)=0$ holds.
In fact, it is obvious because $\calD (f_1)\calD (f_2)=(p+r-q-s)^2\calD (f)$.
From this property, if $f\notin\calNC$ and $f$ is not homogeneous, 
then either $f_1\notin\calNC\ni f_2$ or $f_2\notin\calNC\ni f_1$ holds.
Also, if $K\notin\calNC$, a homogenous binomial $f$ which is not normal crossing is 
added to $Q$ in the finite times of blowup, 
and $f_1$ and $f_2$ obtained by the blow-up of $f$ belong to $\calNC$.

Based on this property, we consider the structure of the blowup tree $\calT=(V,E)$ obtained 
by Algorithm \ref{alg:between_two_vars}.
First, the tree $\calT$ is a binary tree and $h\in\calNC$ is equivalent to $h\in\calL(\calT)$.
If $h\notin\calL(\calT)$ and $d(h)<d(\calT)-1$, one of its child nodes $h_1$ and $h_2$ is a leaf node and 
the other is not.
If $h\notin\calL(\calT)$ and $d(h)=d(\calT)-1$, $h$ is homogeneous binomial and 
both child nodes are $h_1,h_2\in\calL(\calT)$.
Let one of two child leaves $h_1$ and $h_2$ be $g$, 
and consider the path from the root node $f$ to the leaf node $g$.
The set of nodes on the path is called the \textbf{stem} of the blow-up binary tree $\calT$ and 
is denoted by $\calS(\calT)$.
If the considered tree $\calT$ is obviously unique, we simply denote it by $\calS$.
The stem includes the deepest leaf node $g$ and the root node $f$, that is, $f,g\in\calS$.

\subsubsection{RLCT of bivariate sop binomial}
Next, we consider the RLCT for bivariate sop binomials.
First, we consider the case of normal crossing, whose proof requires no algorithm.
We define the quantities that are important to derive the RLCT.

\begin{dfn}[left-right index ratio]
Let $f(w_1,w_2)\in\calSP^{2,2}_{\geq0}\cap\calNC$ and the monomial $g(w_1,w_2)\in\calM^2$ is 
represented as follows.

\[f(w_1,w_2)=w_1^pw_2^q(1+w_1^rw_2^{r'})\ ,\ g(w_1,w_2)=w_1^{s-1}w_2^{t-1}\]
where $p,q,r,r'\in\NN_0,s,t>0$ and $p+r,q+r'\in\NN$.
The following quantities related to indexes

\[\rho^l(f,g)=\frac{p}{s},\ \rho^r(f,g)=\frac{q}{t}\]
are called \textbf{left index ratio} and \textbf{right index ratio} of $(f,g)$ respectively.
The set of index ratios is denoted by $\IR(f,g)$, that is,

\[\IR(f,g)=\{\rho^l(f,g),\rho^r(f,g)\}.\]
\end{dfn}

The following holds for the index ratio.

\begin{prp}[RLCT of bivariate normal crossing sop polynomial]\label{prp:NC_22sop_rlct}
Let $f(w_1,w_2)\in\calSP^{2,2}_{\geq0}\cap\calNC$ and $g(w_1,w_2)\in\calM^2$.
For the maximum pole $-\lambda$ of the zeta function defined below

\[\zeta(z)=\iint_W f(w_1,w_2)^z|g(w_1,w_2)|dw_1dw_2,\]
the following equation holds.

\[\frac{1}{\lambda}=\max\IR(f,g)=\max(\rho^l(f,g),\rho^r(f,g))\]
Where integral domain $W$ is a sufficiently large compact set.

\begin{proof}
Because of the equivalence of RLCT (Proposition \ref{eq:rlct_prp:equality}), 
let $M>0$ be a sufficiently large real number, it is ok to replace the integral domain $W$ with $[-M,M]^2$.
Let $f(w_1,w_2)=w_1^pw_2^q(1+w_1^rw_2^{r'})$ and $g(w_1,w_2)=w_1^{s-1}w_2^{t-1}$.
There exist positive constants $c_1$ and $c_2$ such that 
$c_1<1+w_1^rw_2^{r'}<c_2$ on the compact set $W$.
Since the equivalence of RLCT (Proposition \ref{eq:rlct_prp:equality}),
it is sufficient to consider the pole of the following zeta function.

\[\zeta'(z)=\iint_W(w_1^pw_2^q)^z|w_1^{s-1}w_2^{t-1}|dw_1dw_2\]
Furthermore, because of symmetry, it is sufficient to consider 
the integral domain $w_1\geq0$ and $w_2\geq0$.
We obtain the pole of $\zeta'(z)$ by calculating the following integral.

\[\zeta'(z)=\left(\int_0^Mw_1^{pz+s-1}dw_1\right)\left(\int_0^Mw_2^{qz+t-1}dw_2\right)
=\frac{M^{(p+q)z+(s+t)}}{(pz+s)(qz+t)}\]
Thus $-s/p,-t/q$ is the pole of $\zeta'$.
Where $-s/p,-t/q=-\infty$ if $p,q=0$, and let $\frac{1}{\infty}=0$.
Since $\lambda=\min\left(\frac{s}{p},\frac{t}{q}\right)$ from above discussion,

\[\frac{1}{\lambda}=\max\left\{\frac{p}{s},\frac{q}{t}\right\}=\max\IR(f,g)\]
Therefore it is shown that the statement is true.

\end{proof}
\end{prp}
Next, we consider the RLCT of a bivariate sop binomial that is not a normal crossing.
As in the case of normal crossing, we first define the quantities that are important 
to derive the RLCT.

\begin{dfn}[left-right index ratio and its dual]\label{dfn:prp_and_pdrp}
Let $f(w_1,w_2)\in\calSP^{2,2}_{\geq0}\backslash\calNC$ and 
monomial $g(w_1,w_2)\in\calM^2$ can be represented as follows.

\[f(w_1,w_2)=w_1^pw_2^q(w_1^r+w_2^{r'})\ ,\ g(w_1,w_2)=w_1^{s-1}w_2^{t-1}\]
Where $p,q\in\NN_0,r,r'\in\NN,s,t>0$.
The four following quantities related to the indexes

\[\rho^l(f,g)=\frac{p}{s},\ \rho^r(f,g)=\frac{q}{t},\ 
\nu^l(f,g)=\frac{p+r}{s},\ \nu^r(f,g)=\frac{q+r'}{t},\]
is called \textbf{left index ratio}, \textbf{right index ratio}, \textbf{left dual index ratio}, 
and \textbf{right dual index ratio} respectively.
The set of index ratios and dual index ratios are denoted by $\IR(f,g)$ and $\DIR(f,g)$ respectively, 
that is,

\[\IR(f,g)=\{\rho^l(f,g),\rho^r(f,g)\},\ \DIR(f,g)=\{\nu^l(f,g),\nu^r(f,g)\}.\]

\end{dfn}
We modify the blowup between variables (algorithm \ref{alg:between_two_vars}) 
to that with jacobians (algorithm \ref{alg:between_two_vars_with_jacobian}) 
because jacobians are needed to consider the RLCT.
In this algorithm, we use $V\subset\calSP^{2,2}_{\geq0}\times\calM^2$ for the vertex set 
instead of $V\subset\calSP^{2,2}_{\geq0}$.
It halts because of Proposition \ref{prp:alg_halt} as well.
Let $\calT=(V,E)$ be the blowup tree obtained by this algorithm.
$k\in\pi_2(V)$ is a monomial corresponding to the prior distribution and 
the Jacobian with the variable transformation.
Hereafter, this monomial will be referred to as the \textbf{outer monomial}.
$h\in\pi_1(V)$ is referred to as the \textbf{inner polynomial} in the same way.
Before the proof of the main theorem, we consider the change of the indexes of 
the outer monomial with the blowup of Algorithm \ref{alg:between_two_vars_with_jacobian}.
Let $k(w_1,w_2)=w_1^{s-1}w_2^{t-1}$.
The following monomials $k_1$ and $k_2$ are obtained in the local coordinates with the blowup.

\[k_1(w_1,w_2)=w_2k(w_1w_2,w_2)=w_2w_1^{s-1}w_2^{s-1}w_2^{t-1}=w_1^{s-1}w_2^{t+s-1}\]
\[k_2(w_1,w_2)=w_1k(w_1,w_1w_2)=w_1w_1^{s-1}w_1^{t-1}w_2^{t-1}=w_1^{s+t-1}w_2^{t-1}\]
This indicates that the behavior of the outer monomial with the blowup is 
equal to that of the inner polynomial, considering the apparent indexes $+1$.

\begin{thm}[RLCT of bivariate non-normal crossing sop binomial]\label{thm:nonNC_22sop_rlct}
Let $f(w_1,w_2)\in\calSP^{2,2}_{\geq0}\backslash\calNC$ and $g(w_1,w_2)\in\calM^2$.
For the maximum pole $-\lambda$ of the zeta function defined below

\[\zeta(z)=\iint_W f(w_1,w_2)^z|g(w_1,w_2)|dw_1dw_2,\]
The following equation holds.

\[\frac{1}{\lambda}=\max\left[\IR(f,g)\cup
\left\{\frac{\nu^l(f,g)\nu^r(f,g)-\rho^l(f,g)\rho^r(f,g)}{
\nu^l(f,g)+\nu^r(f,g)-\rho^l(f,g)-\rho^r(f,g)}\right\}\right]\]
Where the integral domain $W$ is a sufficiently large compact set.

\begin{proof}
The following holds for the maximum pole of the zeta function.

\[-\lambda=\underset{u:\Loco}{\max}(-\lambda_u)\]
where $-\lambda_u$ is the maximum pole of the zeta function in the local coordinate $u$.
Let $\calT=(V,E)$ be the blowup tree obtained by Algorithm \ref{alg:between_two_vars_with_jacobian}, 
the local coordinate system is the leaves set $\calL$, and $(h,k)\in\calL$ is equivalent to $h\in\calNC$.
Therefore, by Proposition \ref{prp:NC_22sop_rlct}

\[\frac{1}{\lambda}=\underset{(h,k)\in\calL}{\max}\left\{\max\IR(h,k)\right\}\]
holds. We now prove the following lemma.

\begin{lmm}\label{lmm:from_leaf_to_stem}
The following holds for the stem $\calS$ and leaves $\calL$ of the blowup binary tree $\calT$.

\[\bigcup_{(h,k)\in\calS}\IR(h,k)
=\bigcup_{(h',k')\in\calL}\IR(h',k')\]

\begin{proof}
First, we prove by induction that the following equation holds for any $n=1,2,\ldots,N$ 
if $A_i\cup A_{i+1}=A_{i+1}\cup B_{i+1}$ holds for 
$i=0,1,2,\ldots,N-1$ and the set sequences $(A_i)_{i=0}^N$,$(B_i)_{i=1}^N$.

\begin{equation}\label{eq:lmm_of_lmm_set}
\bigcup_{i=0}^nA_i=\left(\bigcup_{i=0}^{n-1}B_{i+1}\right)\cup A_n
\end{equation}
It is obvious the case of $n=1$.
Next, suppose it holds for the case of $n=k\leq N-1$ and consider the case of $n=k+1$.

\begin{align*}
\bigcup_{i=0}^{k+1}A_i&=\left(\bigcup_{i=0}^kA_i\right)\cup A_{k+1}\\
&=\left\{\left(\bigcup_{i=0}^{k-1}B_{i+1}\right)\cup A_k\right\}\cup A_{k+1}\\
&=\left(\bigcup_{i=0}^{k-1}B_{i+1}\right)\cup\left(A_{k+1}\cup B_{k+1}\right)\\
&=\left(\bigcup_{i=0}^{(k+1)-1}B_{i+1}\right)\cup A_{k+1}
\end{align*}
Thus the equation \eqref{eq:lmm_of_lmm_set} is proved.
Let the depth $d(\calT)=n$ of the blowup binary tree $\calT$.
Since the root node $(f,g)\in V(\calT)$ is $f\notin\calNC$, $n$ is greater than or equal to $1$.
Let the root node be $(h_0,k_0)$, the stem node of the depth $i=1,\ldots,n$ be $(h_i,k_i)$, 
and the other one (leaf node) be $(h'_i,k'_i)$.
For $i=0,\ldots,n-1$,

\[\IR(h_i,k_i)\cup\IR(h_{i+1},k_{i+1})=\IR(h_{i+1},k_{i+1})\cup\IR(h'_{i+1},k'_{i+1})\]
holds. In fact, suppose $r\geq r'>0$ for the blowup between variables and 
let $h_i,k_i$ be as follows$(h_i\notin\calL)$,

\[h_i(w_1,w_2)=w_1^pw_2^q(w_1^r+w_2^{r'})\ ,\ k_i(w_1,w_2)=w_1^{s-1}w_2^{t-1}\]
we get $h_{i+1}$ and $h'_{i+1}$ as below.

\[h_{i+1}(w_1,w_2)=w_1^{p+q+r'}w_2^q(w_1^{r-r'}+w_2^{r'})\ ,\ k_{i+1}(w_1,w_2)=w_1^{s+t-1}w_2^{t-1}\]
\[h'_{i+1}(w_1,w_2)=w_1^pw_2^{q+p+r'}(w_1^rw_2^{r-r'}+1)\ ,\ k'_{i+1}(w_1,w_2)=w_1^{s-1}w_2^{t+s-1}\]
Thus $(h'_{i+1},k'_{i+1})$ belongs to $\calL$ and 
the index ratio set $\IR(h'_{i+1},k'_{i+1})$ is $\{p/s,(q+p+r')/(t+s)\}$.
As for $(h_{i+1},k_{i+1})$, $(h_{i+1},k_{i+1})\notin\calL$ holds if $r>r'$ and 
$(h_{i+1},k_{i+1})\in\calL$ holds if $r=r'$.
And the index ratio set $\IR(h_{i+1},k_{i+1})$ is $\{(p+q+r')/(s+t),q/t\}$.
Thus, substitute $\IR(h_i,k_i)$ and $\IR(h'_{i+1},k'_{i+1})$ for $A_i$ and $B_i$ respectively, 
then the condition \eqref{eq:lmm_of_lmm_set} is satisfied.
Therefore by the equation \eqref{eq:lmm_of_lmm_set},

\[\bigcup_{i=0}^n\IR(h_i,k_i)=
\left(\bigcup_{i=0}^{n-1}\IR(h'_{i+1},k'_{i+1})\right)\cup\IR(h_{n},k_{n})\]
holds. where $(h_n,k_n)\in\calS\cap\calL$, $(h_i,k_i)\in\calS\backslash\calL(i=0,1,\ldots,n-1)$,
$(h'_{i+1},k'_{i+1})\in\calL(i=0,1,\ldots,n-1)$, so the lemma is proved.

\end{proof}
\end{lmm}
(The proof of the theorem \ref{thm:nonNC_22sop_rlct} continues below.)
From Lemma \ref{lmm:from_leaf_to_stem}, 
we can consider the index ratio of the stem set $\calS$ instead of the leaf node set $\calL$,
which is a normal crossing, that is, the following equation holds.

\[\frac{1}{\lambda}=\underset{(h,k)\in\calS}{\max}\left\{\max\IR(h,k)\right\}\]

We then consider the relationship between the index ratios and their duals of $\calS$.
The following lemma holds.

\begin{lmm}\label{lmm:rp_of_stem_property}
Let $\calS$ be the stem of a blowup binary tree $\calT$.
Let the left index ratio, the right index ratio, the left dual index ratio, 
and the right dual index ratio at the trunk nodes of the depth $i=0,1,\ldots,n$ 
be $\rho_i^l,\rho_i^r,\nu_i^l$, and $\nu_i^r$ respectively.
Then the following holds.

\begin{enumerate}
\item\label{lmm_pole_monotone} The sequences $(\rho_i^l)_{i=0}^n,(\rho_i^r)_{i=0}^n,(\nu_i^l)_{i=0}^n$, 
and $(\nu_i^r)_{i=0}^n$ are monotonical.
\item\label{lmm_pole_reverse1} The increase and decrease of $(\rho_i^l)_{i=0}^n$ and $(\nu_i^r)_{i=0}^n$ 
are regressive.
\item\label{lmm_pole_reverse2} The increase and decrease of $(\rho_i^r)_{i=0}^n$ and $(\nu_i^l)_{i=0}^n$ 
are regressive.
\item\label{lmm_pole_or_decrease} Either $(\rho_i^l)_{i=0}^n$ or $(\rho_i^r)_{i=0}^n$ increases.
\end{enumerate}

Where the monotonical sequences $(a_i)_{i=0}^n$ and $(b_i)_{i=0}^n$ are \textbf{regressive} means that
$(a_i)_{i=0}^n$ is monotonically increasing if and only if $(b_i)_{i=0}^n$ is monotonically decreasing.

\begin{proof}
Let the stem node $(h_i,k_i)\in\calS$ of depth $i=0,1,\ldots,n$ be

\[h_i(w_1,w_2)=w_1^pw_2^q(w_1^r+w_2^{r'}),\ k_i(w_1,w_2)=w_1^{s_i-1}w_2^{t_i-1}.\]
where $p_i,q_i,r_i,r'_i\in\NN_0$ and $s_i,t_i>0$.
Let $R_0=r_0=r,R_1=r'_0=r'$ and Euclid's Let $\{\alpha_j\}_{j=1}^N$ and $\{R_j\}_{j=0}^{N+1}$ 
be the quotient and remainder sequences generated by Euclid's algorithm, that is,

\[R_{j-1}=\alpha_jR_j+R_{j+1}\ \ \ (1\leq R_{j+1}< R_j(1\leq^\forall j\leq N-1),R_{N+1}=0).\]
Note that the following changes may be applied for the only last step $R_{N-1}=\alpha_NR_N+R_{N+1}$ and 
one more step is allowed.

\[\left\{\begin{array}{cc}R_{N-1}=(\alpha_N-1)R_N+R_{N+1}&(R_{N+1}=R_N)\\
R_N=1R_{N+1}+R_{N+2}&(R_{N+2}=0)\end{array}\right.\]
This change results in two ways of quotient and remainder sequences 
but can be made to correspond to the deepest leaf selection of the stem in the two ways.
Let $i_{2j}$ and $i_{2j+1}$ be the indexes such that $r_i=R_{2j}$ and $r'_i=R_{2j+1}$ hold for the first time, 
then we get

\[(r_i,r'_i)=\left\{\begin{array}{cc}(R_{2j},R_{2j-1}-(i-i_{2j})R_{2j})&(i_{2j}\leq i\leq i_{2j+1})\\
(R_{2j}-(i-i_{2j+1})R_{2j+1},R_{2j+1})&(i_{2j+1}\leq i\leq i_{2j+2})\end{array}\right..\]
where $i_{j+1}-i_j=\alpha_j$. In the same way,

\[(p_i,q_i)=\left\{
\begin{array}{cc}(p_{i_{2j}},q_{i_{2j-1}}+(i-i_{2j})(p_{i_{2j}}+R_{2j})&(i_{2j}\leq i\leq i_{2j+1})\\
(p_{i_{2j}}+(i-i_{2j+1})(q_{i_{2j+1}}+R_{2j+1}),q_{2j+1})&(i_{2j+1}\leq i\leq i_{2j+2})\end{array}\right.
\]
\[(s_i,t_i)=\left\{
\begin{array}{cc}(s_{i_{2j}},t_{i_{2j-1}}+(i-i_{2j})(s_{i_{2j}}+1)&(i_{2j}\leq i\leq i_{2j+1})\\
(s_{i_{2j}}+(i-i_{2j+1})(t_{i_{2j+1}}+1),t_{i_{2j+1}})&(i_{2j+1}\leq i\leq i_{2j+2})\end{array}\right.
\]
holds. Using the equation above, left index ratios can be represented as

\[\rho_i^l=\frac{p_i}{s_i}=\left\{\begin{array}{cc}
\frac{p_{i_{2j}+1}}{s_{i_{2j}}}&(i_{2j}\leq i\leq i_{2j+1})\\
\frac{p_{i_{2j}}+(i-i_{2j+1})(q_{i_{2j+1}}+r_{i_{2j+1}})}
{s_{i_{2j}}+(i-i_{2j+1})t_{i_{2j+1}}}&(i_{2j+1}\leq i\leq i_{2j+2})
\end{array}\right..\]

Remind that either $\frac{a}{b}\geq\frac{a+c}{b+d}\geq\frac{c}{d}$ or 
$\frac{a}{b}\leq\frac{a+c}{b+d}\leq\frac{c}{d}$ holds for $a,c\geq0$ and $b,d>0$ and 
both hold if and only if $\frac{a}{b}=\frac{c}{d}$.
Because of this fact, either $\rho_{i_{2j}}^l\geq\rho_i^l\geq\nu_{i_{2j+1}}^r$ or 
$\rho_{i_{2j}}^l\leq\rho_i^l\leq\nu_{i_{2j+1}}^r$ holds for $i_{2j+1}\leq i\leq i_{2j+2}$, 
and the former case is monotonically decreasing and the latter case is monotonically increasing 
with respect to the increase of $i$ in this interval.
A similar fact holds for right index ratios.

\[\rho_i^r=\frac{q_i}{t_i}=\left\{\begin{array}{cc}
\frac{q_{i_{2j-1}}+(i-i_{2j})(p_{i_{2j}}+r_{i_{2j}})}
{t_{i_{2j-1}}+(i-i_{2j})s_{i_{2j}}}&(i_{2j}\leq i\leq i_{2j+1})\\
\frac{q_{i_{2j+1}}}{t_{i_{2j+1}}}&(i_{2j+1}\leq i\leq i_{2j+2})
\end{array}\right.\]
Thus either $\rho_{i_{2j-1}}^r\geq\rho_i^r\geq\nu_{i_{2j}}^l$ or 
$\rho_{i_{2j-1}}^r\leq\rho_i^r\leq\nu_{i_{2j}}^l$ holds for $i_{2j}\leq i\leq i_{2j+1}$, 
and the former case is monotonically decreasing and the latter case is monotonically increasing 
with respect to the increase of $i$ in this interval.
In the case of left dual index ratios,

\[\nu_i^l=\frac{p_i+r_i}{s_i}=\left\{\begin{array}{cc}
\frac{p_{i_{2j}}+r_{i_{2j}}}{s_{i_{2j}}}&(i_{2j}\leq i\leq i_{2j+1})\\
\frac{(p_{i_{2j}}+r_{i_{2j}})+(i-i_{2j+1})q_{i_{2j+1}}}
{s_{i_{2j}}+(i-i_{2j+1})t_{i_{2j+1}}}&(i_{2j+1}\leq i\leq i_{2j+2})
\end{array}\right.\]
either $\nu_{i_{2j}}^l\geq\nu_i^l\geq\rho_{i_{2j+1}}^r$ or 
$\nu_{i_{2j}}^l\leq\nu_i^l\leq\rho_{i_{2j+1}}^r$ holds for $i_{2j+1}\leq i\leq i_{2j+2}$,
and the former case is monotonically decreasing and the latter case is monotonically increasing 
with respect to the increase of $i$ in this interval.
In the case of right dual index ratios,

\[\nu_i^r=\frac{q_i+r'_i}{t_i}=\left\{\begin{array}{cc}
\frac{(q_{i_{2j-1}}+r'_{i_{2j-1}})+(i-i_{2j})p_{i_{2j}}}
{t_{i_{2j-1}}+(i-i_{2j})s_{i_{2j}}}&(i_{2j}\leq i\leq i_{2j+1})\\
\frac{q_{i_{2j+1}}+r_{i_{2j+1}}}{t_{i_{2j+1}}}&(i_{2j+1}\leq i\leq i_{2j+2})
\end{array}\right.\]
either $\nu_{i_{2j-1}}^r\geq\nu_i^r\geq\rho_{i_{2j}}^l$ or 
$\nu_{i_{2j-1}}^r\leq\nu_i^r\leq\rho_{i_{2j}}^l$ holds for $i_{2j}\leq i\leq i_{2j+1}$,
and the former case is monotonically decreasing and the latter case is monotonically increasing 
with respect to the increase of $i$ in this interval.

Suppose $\rho_0^l\geq\nu_0^r$. From the definition of indexes, 
$i_0=i_1=0$ holds so $\rho_{i_0}^l\geq\nu_{i_1}^r$.
From the previous discussion of left index ratios, 
we obtain $\rho_{i_0}^l\geq\rho_{i_2}^l\geq\nu_{i_1}^r$.
Then from the previous discussion of right index ratios, 
we obtain $\rho_{i_2}^l\geq\nu_{i_3}^r\geq\nu_{i_1}^r$.
Here we again obtained $\rho_{i_2}^l\geq\nu_{i_3}^r$, 
so the following inequation is obtained by repeatedly applying the discussion of 
left index ratios and right index ratios. 

\[\rho_{i_0}^l\geq\rho_{i_2}^l\geq\cdots\geq\rho_{i_{2\lfloor N/2\rfloor}}^l\geq\nu_{i_{2\lfloor N/2\rfloor\pm1}}^r\geq\cdots\geq\nu_{i_3}^r\geq\nu_{i_1}^r\]
Furthermore, since $\rho_i^l$ is monotonically decreasing and 
$\nu_i^r$ is monotonically increasing in each interval, the inequation

\[\rho_0^l\geq\rho_1^l\geq\cdots\geq\rho_n^l\geq\nu_n^r\geq\cdots\geq\nu_1^r\geq\nu_0^r\]
is obtained. the reverse inequation is obtained in the same way for the case of $\rho_0^l\leq\nu_0^r$.
The same discussion holds for $\rho^r,\nu^l$.
Thus \ref{lmm_pole_monotone},\ref{lmm_pole_reverse1}, and \ref{lmm_pole_reverse2} of the lemma are proved.

Finally, we prove \ref{lmm_pole_or_decrease} of the lemma.
Suppose $\rho^l$ and $\rho^r$. Because of the previous discussion, 
this is equivalent to $\rho_0^l\geq\nu_0^r$ and $\rho_0^r\geq\nu_0^l$.
writing them with $p_0,q_0,r_0r'_0,s_0,t_0$, we get

\[\frac{p_0}{s_0}\geq\frac{q_0+r'_0}{t_0}\ \land\ \frac{q_0}{t_0}\geq\frac{p_0+r_0}{s_0},\]
but this contradicts $\frac{p_0+r_0}{s_0}>\frac{p_0}{s_0}$ by the assumption $r_0>0$.
Therefore \ref{lmm_pole_or_decrease} of the lemma is proved.

\end{proof}
\end{lmm}
(The proof of the theorem \ref{thm:nonNC_22sop_rlct} continues below.)
Let $\calS^l$ be the left-selected stem at the deepest leaves and $\calS^r$ be the right-selected stem.
Where selecting the left means letting $r_n=0$ and selecting the right means letting $r'_n=0$.
Note that $\rho_n^l$ in $\calS^l$ is equal to $\rho_n^r$ in $\calS^r$.
This value is denoted by $\rho_n^*$ in the following.
In the following discussion, the stem is assumed to be $\calS^l$.
The same discussion holds for $\calS^r$.
One of the following three inequations holds for $(\rho^l_i)_{i=0}^n$ and $(\rho^r_i)_{i=0}^n$ 
because of the Lemma \ref{lmm:rp_of_stem_property}.

\begin{enumerate}
\item $\rho_0^l\leq\rho_i^l\leq\rho_n^*\geq\rho_i^r\geq\rho_0^r\ (^\forall i=0,1,\ldots,n)$
\item $\rho_0^l\leq\rho_i^l\leq\rho_n^*\leq\rho_i^r\leq\rho_0^r\ (^\forall i=0,1,\ldots,n)$
\item $\rho_0^l\geq\rho_i^l\geq\rho_n^*\geq\rho_i^r\geq\rho_0^r\ (^\forall i=0,1,\ldots,n)$
\end{enumerate}
Thus whichever inequation holds,

\[
\frac{1}{\lambda}=
\underset{i=0,1,\ldots,n}{\max}\left\{\max\IR(h_i,k_i)\right\}=
\max\left\{\rho_0^l,\rho_0^r,\rho_n^*\right\}=
\max\IR(f,g)\cup\{\rho_n^*\}\]

holds. Next, we show $\rho_n^*=\frac{p_0r'_0+r_0q_0+r_0r'_0}{s_0r'_0+r_0t_0}$.
We use the following lemma.

\begin{lmm}
Let the stem nodes of the blowup binary tree $\calT$ be $\calS=\{(h_i,k_i)\}_{i=0}^n$.
Where index $i$ means the depth of the stem node.
$(h_i,k_i)$ is represented as the following.

\[h_i(w_1,w_2)=w_1^{p_i}w_2^{q_i}(w_1^{r_i}+w_2^{r'_i}),\ k_i(w_1,w_2)=w_1^{s_i-1}w_2^{t_i-1}\]
where $p_i,q_i,r_i,r'_i\in\NN_0$ and $s_i,t_i>0$,$(i=0,1,\ldots,n)$.
For any $a,b,c\in\RR$, there exists a constant $C\in\RR$ such that 
the following equation holds for any $i=0,1,\ldots,n$.

\begin{equation}\label{eq:blowup_invariant}
a\left\{p_ir'_i+r_iq_i+r_ir'_i\right\}+b\left\{s_iq_i-(p_i+r_i)t_i\right\}+c\left\{s_ir'_i+r_it_i\right\}=C
\end{equation}
The equation \eqref{eq:blowup_invariant} is referred to as 
\textbf{invariant equation} with respect to the stem of the blowup binary tree.

\begin{proof}
One of the two following relations holds for blowup between adjacent stems.

\[(p_{i+1},r_{i+1},s_{i+1},q_{i+1},r'_{i+1},t_{i+1})=\left\{\begin{array}{cc}
(p_i+q_i+r'_i,r_i-r'_i,s_i+t_i,q_i,r'_i,t_i)&(r_i\geq r'_i)\\
(p_i,r_i,s_i,q_i+p_i+r_i,r'_i-r_i,t_i+s_i)&(r'_i\geq r_i)
\end{array}\right.\]
We substitute in each case and show that the lemma holds.
It is sufficient to confirm that each coefficient of $a,b$, and $c$ does not change.
Suppose $r_i\geq r'_i$. The coefficient of $a$ is

\[
p_{i+1}r'_{i+1}+r_{i+1}q_{i+1}+r_{i+1}r'_{i+1}=
(p_i+q_i+r'_i)r'_i+(r_i-r'_i)q_i+(r_i-r'_i)r'_i
=p_ir'_i+r_iq_i+r_ir'_i.\]
Since there is the symmetry of the equation and transformation between $(p,r)$ and $(q,r')$, 
the same discussion holds for the case of $r_i\leq r'_i$.
Thus the coefficient of $a$ is invariant for the blowup.
For the case of $r_i\geq r'_i$, the coefficient of $b$ is

\[
s_{i+1}q_{i+1}-(p_{i+1}+r_{i+1})t_{i+1}
=(s_i+t_i)q_i-(p_i+q_i+r'_i+r_i-r'_i)t_i
=s_iq_i-(p_i+r_i)t_i.
\]
For the case of $r_i\leq r'_i$, it is
\[
s_{i+1}q_{i+1}-(p_{i+1}+r_{i+1})t_{i+1}
=s_i(q_i+p_i+r_i)-(p_i+r_i)(t_i+s_i)
=s_iq_i-(p_i+r_i)t_i.
\]
Thus it is invariant for the blowup.
Finally, we calculate the coefficient of $c$. For the case of $r_i\geq r'_i$, it is
\[
r_{i+1}t_{i+1}+s_{i+1}r'_{i+1}
=(r_i-r'_i)t_i+(s_i+t_i)r'_i
=r_it_i+s_ir'_i.
\]
Since there is the symmetry of the equation and transformation between $(r,s)$ and $(r',t)$, 
the same discussion holds for the case of $r_i\leq r'_i$.
Thus it is invariant for the blowup.
Since the coefficient of all $a,b$, and $c$ is invariant for the blowup among stems, 
the summation of them are invariant for fixed $a,b$, and $c$, that is, the lemma is proved.
\end{proof}
(Note)
The proof is given above but in the following, 
we will show how to find blow-up invariants such as those obtained in the lemma.
Let $v=(p,r,s,q,r',t)^\top\in\RR^6$ be a vector of variables.
We derive a quadratic invariant of the variable vector $v$.
Let $A\in\RR^{6\times6}$ be a symmetric matrix, and the quadratic invariant of the blowup is 
expected to be represented in the quadratic form $v^\top Av=C$, where $C$ is a constant.
There are two blowups between the stem and we expect them to be invariant with respect to each other.
Let the linear transformation matrices of the exponents with blowups be $P,Q\in\RR^{6\times6}$, 
the objective is to design $A$ so that the following equation holds.

\begin{equation}\label{invariant_form}
(Pv)^\top A(Pv)=(Qv)^\top A(Qv)=v^\top Av=C
\end{equation}
where

\[P=\bmat{\mat{&&\\&I_3&\\&&}&\mat{1&1&0\\0&-1&0\\0&0&1}\\
\mat{&&\\&O_{3,3}&\\&&}&\mat{&&\\&I_3&\\&&}}\ ,\ 
Q=\bmat{\mat{&&\\&I_3&\\&&}&\mat{&&\\&O_{3,3}&\\&&}\\
\mat{1&1&0\\0&-1&0\\0&0&1}&\mat{&&\\&I_3&\\&&}}.\]
It is sufficient to satisfy $P^\top AP-A=O$ and $Q^\top AQ-A=O$ in order to hold \eqref{invariant_form}, 
so we solve this equation for each component of $A$.
Note that $A$ is assumed to be a symmetric matrix, so the number of variables is $21$, 
which corresponds to the number of components in the upper triangle.
The solution to this equation is as follows.

\[A=\bmat{O_{3,3}&A'\\A'^\top&O_{3,3}}\ ,\ 
A'=\bmat{0&a_{25}&-a_{34}\\a_{25}&a_{25}&-a_{34}+a_{35}\\a_{34}&a_{35}&0}
\]
Thus the quadratic invariant with blowup is obtained by

\[v^\top Av=a_{25}\{pr'+rq+rr'\}+a_{34}\left\{sq-(p+r)t\right\}+a_{35}\left\{sr'+rt\right\}.\]
Replacing $a_{25},a_{34}$, and $a_{35}$ with $a,b$, and $c$ respectively, we obtain

\[a\{pr'+rq+rr'\}+b\left\{sq-(p+r)t\right\}+c\left\{sr'+rt\right\}=C.\]
The linear and cubic invariants are obtained in the same way.

\end{lmm}
(The proof of the theorem \ref{thm:nonNC_22sop_rlct} continues below.)
Note that $r_n=0$ and $r'_n>0$ since we are now considering $\calS^l$.
Substituting $(a,b,c)=(1,0,0)$ in \eqref{eq:blowup_invariant} and comparing $i=0$ and $i=n$, 
we obtain

\[p_0r'_0+r_0q_0+r_0r'_0=p_nr'_n.\]
Also Substituting $(a,b,c)=(0,0,1)$ in \eqref{eq:blowup_invariant} and comparing $i=0$ and $i=n$, 
we obtain

\[s_0r'_0+r_0t_0=s_nr'_n.\]
Thus we obtain

\[\rho_n^*=\rho_n^l=\frac{p_n}{s_n}=\frac{p_nr'_n}{s_nr'_n}=\frac{p_0r'_0+r_0q_0+r_0r'_0}{s_0r'_0+r_0t_0}.\]
Finally, we transform $\rho_n^*$.

\[
\frac{p_0r'_0+r_0q_0+r_0r'_0}{s_0r'_0+r_0t_0}=
\frac{\frac{(p_0+r_0)}{s_0}\frac{(q_0+r'_0)}{t_0}-\frac{p_0}{s_0}\frac{q_0}{t_0}}{
\frac{s_0r'_0}{s_0t_0}+\frac{r_0t_0}{s_0t_0}}
=\frac{\nu^l(f,g)\nu^r(f,g)-\rho^l(f,g)\rho^r(f,g)}{\nu^l(f,g)+\nu^r(f,g)-\rho^l(f,g)-\rho^r(f,g)}.
\]
The theorem is proved.
\end{proof}
\end{thm}

\subsection{case of general sop binomial}\label{sec:d2_case}
Next, we consider general sop binomials $\calSP^{d,2}$ by extending the number of variables.
The algorithm for the general sop binomials and the RLCT are based on 
the results obtained in the case of bivariate sop binomials.
First, we state the properties of$\calSP^{d,2}_{\geq0}$ as in the bivariate case.

\begin{prp}\label{prp:deg_is_even_d2}
Let $f\in\calSP ^{d,2}_{\geq0}$. Any component $a_{ij}$ of multi-indexes matrix 
$A_f=(a_{ij})\in\NN_0^{d\times2}$ is even.

\begin{proof}
It is proved in the same way as Proposition \ref{prp:deg_is_even_22}.

\end{proof}
\end{prp}
From Proposition \ref{prp:deg_is_even_d2}, in the local coordinates obtained by blowup, 
the binomial is a normal crossing if it contains a constant term when factoring out a common factor.

\subsubsection{blowup algorithm for sop binomial}\label{sec:blowup_algorithm_2sop}
In this section, we consider the construction of the blow-up algorithm and whether it halts or not.
In the bivariate case, There is a only blowup centered at the origin $\{w=(w_1,w_2)\in\RR^2|w_1=w_2=0\}$.
When generalizing the number of variables, 
we can consider multiple non-singular sets as the center of the blow-up. 
However, as mentioned in Section \ref{sec:ex_sop_rlct_alg}, 
even if such a set is chosen randomly or based on certain rules, 
the blowup algorithm does not necessarily halt.
Therefore, we show that the proposed algorithm (Algorithm \ref{alg:between_terms}, blowup between terms) halts, 
which is a modification of the maximum degree selective algorithm \ref{alg:argmax_selected_blowup} 
in the previous study\cite{3esop} which is expected to halt.

\begin{thm}\label{thm:terms_blowup_halt}
The algorithm \ref{alg:between_terms} halts for any $f(w)\in\calSP^{d,2}$ and $g(w)\in\calM^d$.

\begin{proof}
It is sufficient to focus on the inner polynomial $h$ in discussing that the algorithm halts.
We define $\calD:\calSP^{d,2}_{\geq0}\rightarrow\ZZ^d$ as $\calD:h\mapsto A_h^\top\bmat{1\\-1}$.

First, we show that $\|h\|_\infty\geq\|h'\|_\infty$ holds for the root node $(h,k)$ 
and each leaf node $(h',k')\in\calL(\calT')$ of a tree $\calT'$ obtained 
by a blowup between variables (Algorithm \ref{alg:between_two_vars_with_jacobian}).
Where the norm $\|\cdot\|_\infty:\RR^d\rightarrow\RR_{\geq0}$ is the maximum norm 
defined by $\|\cdot\|_\infty:x\mapsto\underset{i=1,2,\ldots,d}{\max}|x_i|$.
Let $h=c(p_1+p_2)$ be the inner polynomial of the stem node in the blow-up tree $\calT'$ 
except for the deepest leaf.
Where $c,p_1$, and $p_2$ are monomials, $c$ is a common factor, and  $p_1$ and $p_2$ are prime to each other.
Let $s_1$ be the variables whose degree is the maximum of the $p_1$, and 
let $s_2$ be the variables whose degree is the maximum of the $p_2$.
Then we can represent $p_1=p'_1s_1^r$ and $p_2=p'_2s_2^{r'}$.
(The monomials $p'_1$ and $p'_2$ do not contain the variables $s_1$ and $s_2$ respectively.)
In the following, we assume that $r\geq r'$ from symmetry.
Note that $\|\calD(h)\|_{\infty}=r$ by the definition.
Blow-up centered with $\{w\in\RR^d|s_1=s_2=0\}$ for $h$, 
then we obtain the following polynomial $h_a$ and $h_b$ in two local coordinates.

\begin{eqnarray}
h_a&=&cs_1^{r'}(p'_1s_1^{r-r'}+p'_2s_2^{r'}),\label{eq:blownup_poly1}\\
h_b&=&cs_2^{r'}(p'_1s_1^rs_2^{r-r'}+p'_2)\label{eq:blownup_poly2}
\end{eqnarray}
$h_a$ is a leaf node if and only if $r-r'=0$, otherwise it is a stem node.
(both leaf and stem are possible.)
By focusing on the order of each variable, we obtain the following.

\begin{eqnarray}
\|\calD(h_a)\|_\infty&=&\max\{r-r',\mathrm{maxdeg}(p'_1),r',\mathrm{maxdeg}(p'_2)\}\label{eq:blownup_deg1}\\
\|\calD(h_b)\|_\infty&=&\max\{r,r-r',\mathrm{maxdeg}(p'_1),\mathrm{maxdeg}(p'_2)\}\label{eq:blownup_deg2}
\end{eqnarray}
where $\mathrm{maxdeg}$ means the maximum degree of the monomial.
Since $\mathrm{maxdeg}(p'_1)\leq r$ and $\mathrm{maxdeg}(p'_2)\leq r'\leq r$, 
we obtain $\|\calD(h_a)\|_\infty\leq\|\calD(h_b)\|_\infty=\|\calD(h)\|_\infty$.
Since this holds for all stems of a binary tree $\calT'$ except for the deepest leaf, 
by repeated application, $\|\calD(h')\|_\infty\leq\|\calD(h)\|_\infty$ holds 
for any leaf node $h'\in\calL(\calT')$.

We then show that $V(\calT)$ is finite.
Note that $V(\calT)$ is finite if and only if the blow-up algorithm halts, 
and that $\calD(h)\geq0\lor\calD(h)\leq0\Longleftrightarrow h\in\calNC$.
We show that by deriving a contradiction.
Suppose there are the sequence of the blowup trees $(\calT'_i)_{i=1}^\infty$ obtained by 
blowup between variables(Algorithm \ref{alg:between_two_vars_with_jacobian}) and 
the sequence of the inner polynomials $(h_i)_{i=1}^\infty$.
Where $h_i\in\calL(\calT'_i)$ and $\calT'_{i+1}$ is a root node of $\calT'_{i+1}$.
Since the sequence $(\|\calD(h_i)\|_\infty)_{i=1}^\infty$ is 
a monotonically non-increasing sequence of integers, 
there exist we have natural numbers $k,M\in\NN$ such that 
$i\geq k\Longrightarrow\|\calD(h_i)\|_\infty=M$ holds.
We now show the following lemma.

\begin{lmm}\label{lmm:equal_condition}
Let $h_i=c_i(p_i^{(1)}+p_i^{(2)})$.
Where $c_i,p_i^{(1)}$, and $p_i^{(2)}$ are monomials such that 
$p_i^{(1)}$ and $p_i^{(2)}$ have no common factor.
Furthermore, let $p_i^{(1)}=(s_1^{(i)})^{r_i}\times{p'_i}^{(1)}$ and 
$p_i^{(2)}=(s_2^{(i)})^{r'_i}\times{p'_i}^{(2)}$.
Where $r_i$ and $r'_i$ are maximum degree of $p_i^{(1)}$ and $p_i^{(2)}$ respectively, and 
$s_1^{(i)}$ and $s_2^{(i)}$ are the variables whose degrees are them.
Then, the necessary and sufficient condition to hold $\|\calD(h_i)\|_\infty=\|\calD(h_{i+1})\|_\infty$ is 
satisfying one of the following conditions.
Where they include the case where $(r_i,{p'_i}^{(1)})$ and $(r'_i,{p'_i}^{(2)})$ are reversed.

\begin{enumerate}
\item\label{lmm:equal_cond1} $\|\calD(h_i)\|_\infty=r_i=r'_i$ holds.
\item\label{lmm:equal_cond2} $\|\calD(h_i)\|_\infty=r_i=\mathrm{maxdeg}({p'_i}^{(1)})>r'_i$ holds.
\item\label{lmm:equal_cond3} $\|\calD(h_i)\|_\infty=r_i>\mathrm{maxdeg}({p'_i}^{(1)}),r'_i$ and 
$d(h_{i+1})-d(h_i)=1$ hold.
\end{enumerate}

\begin{proof}
First, we show that $\|\calD(h_i)\|_\infty=\|\calD(h_{i+1})\|_\infty$ when either of the above holds.
If condition \ref{lmm:equal_cond1} holds, $|\calL(\calT_{i+1})|$ is $2$, 
and from equations \eqref{eq:blownup_deg1} and \eqref{eq:blownup_deg2}, 
$\|\calD(h_{i+1})\|_\infty=\|\calD(h_i)\|_\infty$ for $h_{i+1}\in\calL(\calT_{i+1})$.
If condition \ref{lmm:equal_cond2} holds, the variables in ${p'_i}^{(1)}$ are not 
involved in blow-up in $\calT_{i+1}$, so the degree of the variables in ${p'_i}^{(1)}$ are 
equal for any inner polynomial in $V(\calT_{i+1})$.
Thus $\|\calD(h_{i+1})\|_\infty=\|\calD(h_i)\|_\infty$ holds.
Finally, if condition \ref{lmm:equal_cond3} holds, 
there is only one $h_{i+1}\in\calL(\calT_{i+1})$ such that $d(h_{i+1})-d(h_i)=1$ holds, 
and satisfies equation \eqref{eq:blownup_deg2}.
Therefore, $\|\calD(h_{i+1})\|_\infty=\|\calD(h_i)\|_\infty$ holds.

Then, we show that the converse is true, that is, 
$\|\calD(h_{i+1})\|_\infty>\|\calD(h_i)\|_\infty$ holds if the following condition holds.

\begin{enumerate}
\setcounter{enumi}{3}
\item \label{lmm:not_equal_cond4} $\|\calD(h_i)\|_\infty=r_i>\mathrm{maxdeg}({p'_i}^{(1)}),r'_i$ and 
$d(h_{i+1})-d(h_i)\geq2$ hold.
\end{enumerate}

If the condition holds, because of $d(h_{i+1})-d(h_i)\geq2$, 
$\|\calD(h_{i+1})\|_\infty\leq\|\calD(h')\|_\infty$ hold, 
where $h'$ is the stem node of the depth $1$ in the blowup binary tree $\calT_{i+1}$.
Since the equation \eqref{eq:blownup_deg1} holds for the norm of the stem node, 
$\|\calD(h')\|_\infty<\|\calD(h_i)\|_\infty$ holds.
Therefore the lemma is proved.

\end{proof}
\end{lmm}
(The proof of the theorem \ref{thm:terms_blowup_halt} continues below.)
From Lemma \ref{lmm:equal_condition}, 
one of the conditions \ref{lmm:equal_cond1},\ref{lmm:equal_cond2}, and \ref{lmm:equal_cond3} holds 
for the $k$-th and following numbers in the sequence.
We then consider the order in which these conditions hold.
Once conditions \ref{lmm:equal_cond2} or \ref{lmm:equal_cond3} is satisfied for some $i(\geq k)$-th, 
conditions \ref{lmm:equal_cond1} are never satisfied after $i$-th.
Because the difference between each component of $\calD(h)$ and $0$ does not become large among blowups.
This fact can be confirmed from equations \eqref{eq:blownup_poly1} and \eqref{eq:blownup_poly2}.
Similarly, once condition \ref{lmm:equal_cond3} is satisfied for some $i(\geq k)$-th, 
condition \ref{lmm:equal_cond2} is never satisfied after the $i$-th.
Because $r_i>r_i-r'_i,\mathrm{maxdeg}({p'_i}^{(1)})$ holds 
although $h_{i+1}$ can be applied to equation \eqref{eq:blownup_poly2}.
From the above, the order in which the conditions hold is 
\ref{lmm:equal_cond1}$\rightarrow$\ref{lmm:equal_cond2}$\rightarrow$\ref{lmm:equal_cond3}.
Where the number of times each condition holds can be zero.

Finally, we show that the number of times each condition holds is finite, which leads to a contradiction.
First, we consider the number of times condition \ref{lmm:equal_cond1} continually holds.
If $r_i=r'_i(i\geq k)$ holds, $\calD(h_j)_{s_t^{(j)}}=0(j>i,t=1\lor 2)$ holds 
for one of two variables $s_1^{(i)}$ and $s_2^{(i)}$ which is center of blowup.
Where $\calD(h_j)_{s_t^{(j)}}$ is the component of the vector $\calD(h_j)$ 
with respect to the variable $s_t^{(j)}$.
Therefore, the number of times condition \ref{lmm:equal_cond1} continually holds after the $i$-th is 
not greater than the number of variables whose component of $\calD(h_i)$ is $\pm r_i$.
Thus the number of times condition \ref{lmm:equal_cond1} continually holds is finite.

We consider the number of times condition \ref{lmm:equal_cond2} continually holds.
We consider the upper bound of the times for the case of $d(h_{i+1})-d(h_i)=1$ and 
$d(h_{i+1})-d(h_i)\geq2$ separately.
The number of times that $d(h_{i+1})-d(h_i)=1$ holds is not greater than 
the dimension of the variable ($d$-times).
Because it is needed not to satisfy the normal crossing condition $\calD(h_j)\geq0\lor\calD(h_j)\leq0$.
The number of times that $d(h_{i+1})-d(h_i)\geq2$ is not greater than 
the number of variables whose components of $\calD(h_j)$ are $\pm r_j$.
Therefore, the number of times where condition \ref{lmm:equal_cond2} continually holds is finite.

Finally, we consider the number of times that condition \ref{lmm:equal_cond3} continually holds.
This number of times is not greater than the dimension of the variable ($d$-times), 
since it is needed not to satisfy $\calD(h_j)\geq0\lor\calD(h_j)\leq0$.
Therefore, the number of times that condition \ref{lmm:equal_cond3} continually holds.

Thus all condition does not allow the iteration of the infinite times, 
which contradicts the assumption that the sequences are infinite.
From the above, the theorem is proved.

\end{proof}
\end{thm}

\subsubsection{RLCT of sop binomial}\label{sec:d2_case_rlct}

In this subsubsection, we consider the RLCT of the sop binomial.
First, we check the case of a normal crossing.
The method is the same as for the bivariate case.
In the following, for $w\in\RR^d$ and $a\in\RR^d_{\geq0}$, we define

\[w^a=\prod_{i=1}^dw_i^{a_i}.\]
$a$ is referred to as \textbf{multi-index}.

\begin{dfn}[index ratio]
Let $f(w)\in\calSP^{d,2}_{\geq0}\cap\calNC$ and $g(w)\in\calM^d$, and they are represented as follows.

\[f(w)=w^p(1+w^r)\ ,\ g(w)=w^{s-1}\]
where $p,r\in\NN_0^d,s\in\RR^d_{>0}$ and $p+r\in\NN^d$.
We define the following vector for multi-index.

\[\rho(f,g)=\frac{p}{s}=\left(\frac{p_i}{s_i}\right)_{i=1}^d\]
The components of the vector are referred to as \textbf{index ratios} of $(f,g)$.

\end{dfn}

The following proposition holds for index ratios.

\begin{prp}[RLCT of normal crossing sop binomial]\label{prp:NC_d2sop_rlct}
Let $f(w)\in\calSP^{d,2}_{\geq0}\cap\calNC$ and $g(w)\in\calM^d$.
For the maximum pole $-\lambda$ of the zeta function defined below

\[\zeta(z)=\int_W f(w)^z|g(w)|dw,\]
the following equation holds.

\[\frac{1}{\lambda}=\max\rho(f,g)\]
Where the integral domain $W$ is a sufficiently large compact set.

\begin{proof}
The method of proof is the same as in the case of Proposition \ref{prp:NC_22sop_rlct}.

\end{proof}
\end{prp}

We discuss the case of non-normal crossing sop binomials.

\begin{dfn}[index ratio and potential index ratio]\label{dfn:prp_and_pdrp_on_d2_sop}
Let $f(w)\in\calSP^{d,2}_{\geq0}\backslash\calNC$ and $g(w)\in\calM^d$, 
and they are represented as follows.

\[f(w)=w^p(w^r+w^{r'})\ ,\ g(w)=w^{s-1}\]
Where $p,r,r'\in\NN^d_0,s\in\RR^d_{\geq0}$, $r,r'\neq0$ and $r^\top r'=0$.
We define the following vector for multi-index.

\[\rho(f,g)=\frac{p}{s}=\left(\frac{p_i}{s_i}\right)_{i=1}^d\]
The components of the vector are referred to as \textbf{index ratios} of $(f,g)$.
We define the following index sets: $I_{(f,g)}=\{i\in\{1,2,\ldots,d\}|r_i>0\}$,
$J_{(f,g)}=\{j\in\{1,2,\ldots,d\}|r'_j>0\}$, 
and $H_{(f,g)}=\{1,2,\ldots,d\}\backslash(I_{(f,g)}\cup J_{(f,g)})$.
We define the following matrix for multi-indexes.

\[\pi(f,g)=\left(\frac{p_ir'_j+r_ip_j+r_ir'_j}{s_ir'_j+r_is_j}\right)_{i\in I_{(f,g)},j\in J_{(f,g)}}
\in\RR^{|I_{(f,g)}|\times|J_{(f,g)}|}\]
The components of the matrix are referred to as \textbf{potential index ratios}.

\end{dfn}
Note that neither $I_{(f,g)}$ nor $J_{(f,g)}$ 
defined for non-normal crossing $f\in\calSP^{d,2}_{\geq0}\backslash\calNC$ is an empty set.
With the above preparation, we prove the following theorem.

\begin{thm}[RLCT of non-normal crossing sop binomial]\label{thm:nonNC_d2sop_rlct}
Let $f(w)\in\calSP^{d,2}_{\geq0}\backslash\calNC$ and $g(w)\in\calM^d$.
For the maximum pole $-\lambda$ of the zeta function defined below

\[\zeta(z)=\int_Wf(w)^z|g(w)|dw,\]
The following equation holds.

\[\frac{1}{\lambda}=\max\left\{\max\rho(f,g),\max\pi(f,g)\right\}\]
Where the integral domain $W$ is a sufficiently large compact set.

\begin{proof}
In preparation for the proof, we prove the following lemma at first.

\begin{lmm}\label{lmm_all_monotonous}
Let $f(w)\in\calSP^{d,2}_{\geq0}\backslash\calNC$, $g(w)\in\calM^d$, 
$i\in I_{(f,g)}$, and $j\in J_{(f,g)}$.
Also, let the blowup binary tree obtained by 
blowup between variables centered with $\{w\in\RR^d|w_i=w_j=0\}$ 
(Algorithm \ref{alg:between_two_vars_with_jacobian}) be $\calT$ and its stem be $\calS$.
$(f_d,g_d)$ is defined as the stem node of depth $d=0,1,\ldots,D$.
Then let $M_d=\max\{\max\rho(f_d,g_d),\max\pi(f_d,g_d)\}$ and 
$m_d=\min\{\min\rho(f_d,g_d),\min\pi(f_d,g_d)\}$,$(d=0,1,\ldots,D)$, 
the sequence $(M_d)_{d=0}^D$ is monotonically non-increasing and 
the sequence $(m_d)_{d=0}^D$ is monotonically non-decreasing.

\begin{proof}
The components $w_h,w_{h'}(h,h'\neq i,j)$ are not involved in the blowup, 
so $\rho(f_d,g_d)_h,\pi(f_d,g_d)_{h,h'}(d=0,1,\ldots,D)$ are constants with respect to depth $d$.
From the invariant equation \eqref{eq:blowup_invariant} for blowups, 
$\pi(f_d,g_d)_{i,j}(d=0,1,\ldots,D)$ are also constants with respect to depth $d$.
Thus it is sufficient to focus on the change of 
$\rho(f_d,g_d)_i,\rho(f_d,g_d)_j,\pi(f_d,g_d)_{i,h'},\pi(f_d,g_d)_{h,j}(d=0,1,\ldots,D)$ 
with respect to depth $d$.

First, we consider $(\rho(f_d,g_d)_i)_{d=0}^D,(\rho(f_d,g_d)_j)_{d=0}^D$.
From Lemma \ref{lmm:rp_of_stem_property}, $(\rho(f_d,g_d)_i)_{d=0}^D,(\rho(f_d,g_d)_j)_{d=0}^D$ are 
monotonical sequences and one of the following three holds for $\pi(f_0,g_0)_{i,j}$.

\begin{enumerate}
\item $\rho(f_0,g_0)_i\leq\rho(f_d,g_d)_i\leq\pi(f_0,g_0)_{i,j}\geq
\rho(f_d,g_d)_j\geq\rho(f_0,g_0)_j\ (^\forall d=0,1,\ldots,D)$
\item $\rho(f_0,g_0)_i\leq\rho(f_d,g_d)_i\leq\pi(f_0,g_0)_{i,j}\leq
\rho(f_d,g_d)_j\leq\rho(f_0,g_0)_j\ (^\forall d=0,1,\ldots,D)$
\item $\rho(f_0,g_0)_i\geq\rho(f_d,g_d)_i\geq\pi(f_0,g_0)_{i,j}\geq
\rho(f_d,g_d)_j\geq\rho(f_0,g_0)_j\ (^\forall d=0,1,\ldots,D)$
\end{enumerate}
Thus $(\max\{\max\rho(f_d,g_d),\pi(f_d,g_d)_{i,j}\})_{d=0}^D$ and 
$(\min\{\min\rho(f_d,g_d),\pi(f_d,g_d)_{i,j}\})_{d=0}^D$ are monotonically non-increasing and
non-decreasing respectively.

Next, we consider $(\pi(f_d,g_d)_{i,h'})_{d=0}^D$,$(\pi(f_d,g_d)_{h,j})_{d=0}^D$ 
for the case of $h\neq i,h\in I_{(f,g)},h'\neq j,h'\in J_{(f,g)}$.
Because of symmetry, it is sufficient to consider $(\pi(f_d,g_d)_{h,j})_{d=0}^D$.
Let $f_d,g_d$ be represented as follows.

\[f_d(w)=w_i^{p^{(d)}_i}w_h^{p^{(d)}_h}w_j^{p^{(d)}_j}
(w_{\backslash i,h,j}^*w_i^{r^{(d)}_i}w_h^{r^{(d)}_h}+w_{\backslash i,h,j}^*w_j^{{r'_j}^{(d)}}),
\ g_d(w)=w_{\backslash i,h,j}^*w_i^{s^{(d)}_i-1}w_h^{s^{(d)}_h-1}w_j^{s^{(d)}_j-1}\]
Where $w_{\backslash i,h,j}^*$ is the monomial which composes of the components of 
the variables vector $w\in\RR^d$ except for $i,h,j$-th components.
And $p_{\tau}^{(d)}(\tau=i,h,j)$,$r_{\tau}^{(d)}(\tau=i,h),{r'_j}^{(d)}\in\NN_0$,
$s_i^{(d)},s_h^{(d)},s_j^{(d)}>0$ hold.
In the above conditions, $\pi(f_d,g_d)_{i,j},\pi(f_d,g_d)_{h,j}$can be represented as follows.

\[\pi(f_d,g_d)_{i,j}=\frac{p_i^{(d)}{r'_j}^{(d)}+r_i^{(d)}p_j^{(d)}+r_i^{(d)}{r'_j}^{(d)}}
{s_i^{(d)}{r'_j}^{(d)}+r_i^{(d)}s_j^{(d)}},\ \pi(f_d,g_d)_{h,j}=
\frac{p_h^{(d)}{r'_j}^{(d)}+r_h^{(d)}p_j^{(d)}+r_h^{(d)}{r'_j}^{(d)}}
{s_h^{(d)}{r'_j}^{(d)}+r_h^{(d)}s_j^{(d)}}\]
The denominator and numerator of $\pi(f_d,g_d)_{i,j}$ are 
constants from invariant equation \eqref{eq:blowup_invariant}.
Let $C_{ij}=p_i^{(d)}{r'_j}^{(d)}+r_i^{(d)}p_j^{(d)}+r_i^{(d)}{r'_j}^{(d)}$ and 
$C'_{ij}=s_h^{(d)}{r'_j}^{(d)}+r_h^{(d)}s_j^{(d)}$.
By the calculation of $\pi(f_d,g_d)_{h,j}-\pi(f_d,g_d)_{i,j}$,

\begin{align*}
&\pi(f_d,g_d)_{h,j}-\pi(f_d,g_d)_{i,j}
=\frac{p_h^{(d)}{r'_j}^{(d)}+r_h^{(d)}p_j^{(d)}+r_h^{(d)}{r'_j}^{(d)}}
{s_h^{(d)}{r'_j}^{(d)}+r_h^{(d)}s_j^{(d)}}-
\frac{C_{ij}}{C'_{ij}}\\
&=\left\{
p_h^{(d)}{r'_j}^{(d)}s_i^{(d)}{r'_j}^{(d)}+r_h^{(d)}p_j^{(d)}s_i^{(d)}{r'_j}^{(d)}
+r_h^{(d)}{r'_j}^{(d)}s_i^{(d)}{r'_j}^{(d)}+p_h^{(d)}{r'_j}^{(d)}r_i^{(d)}s_j^{(d)}
+r_h^{(d)}p_j^{(d)}r_i^{(d)}s_j^{(d)}\right.\\
&\left.+r_h^{(d)}{r'_j}^{(d)}r_i^{(d)}s_j^{(d)}-p_i^{(d)}{r'_j}^{(d)}s_h^{(d)}{r'_j}^{(d)}
-r_i^{(d)}p_j^{(d)}s_h^{(d)}{r'_j}^{(d)}-r_i^{(d)}{r'_j}^{(d)}s_h^{(d)}{r'_j}^{(d)}
-p_i^{(d)}{r'_j}^{(d)}r_h^{(d)}s_j^{(d)}\right.\\
&\left.-r_i^{(d)}p_j^{(d)}r_h^{(d)}s_j^{(d)}
-r_i^{(d)}{r'_j}^{(d)}r_h^{(d)}s_j^{(d)}\right\}/(s_h^{(d)}{r'_j}^{(d)}+r_h^{(d)}s_j^{(d)})C'_{ij}\\
&=\frac{
p_h^{(d)}{r'_j}^{(d)}C'_{ij}+r_h^{(d)}p_j^{(d)}s_i^{(d)}{r'_j}^{(d)}
+r_h^{(d)}{r'_j}^{(d)}s_i^{(d)}{r'_j}^{(d)}-C_{ij}s_h^{(d)}{r'_j}^{(d)}
-p_i^{(d)}{r'_j}^{(d)}r_h^{(d)}s_j^{(d)}}
{(s_h^{(d)}{r'_j}^{(d)}+r_h^{(d)}s_j^{(d)})C'_{ij}}\\
&=\frac{{r'_j}^{(d)}\left[p_h^{(d)}C'_{ij}-C_{ij}s_h^{(d)}+r_h^{(d)}
\left\{s_i^{(d)}\left(p_j^{(d)}+{r'_j}^{(d)}\right)-p_i^{(d)}s_j^{(d)}\right\}\right]}
{(s_h^{(d)}{r'_j}^{(d)}+r_h^{(d)}s_j^{(d)})C'_{ij}}
\end{align*}
Where $s_i^{(d)}\left(p_j^{(d)}+{r'_j}^{(d)}\right)-p_i^{(d)}s_j^{(d)}$ is a constant 
with respect to depth $d$ since it is obtained 
by substituting $(a,b,c)=(0,1,0)$ in invariant equation \eqref{eq:blowup_invariant}.
This shows that $\pi(f_d,g_d)_{h,j}-\pi(f_d,g_d)_{i,j}$ approaches 0 monotonically with respect to depth $d$.
Thus, $\max\pi(f_d,g_d)$ and $\min\pi(f_d,g_d)$ is monotonically 
non-increasing and non-decreasing respectively.

From the above discussion, since $(\max\{\max\rho(f_d,g_d),\pi(f_d,g_d)_{i,j}\})_{d=0}^D$ and 
$\max\pi(f_d,g_d)$ are monotonically non-increasing, 
and the union $M_d=\max\{\max\rho(f_d,g_d),\max\pi(f_d,g_d)\}$ is also monotonically non-increasing.
$m_d=\min\{\min\rho(f_d,g_d),\min\pi(f_d,g_d)\}$ is monotonically non-decreasing in the same way.
Thus, the lemma is proved.

\end{proof}
\end{lmm}
(The proof of the theorem \ref{thm:nonNC_d2sop_rlct} continues below.)
Let $\calT'$ be the blowup tree obtained by blowup between variables 
(Algorithm \ref{alg:between_two_vars_with_jacobian}) centered with $\{w\in\RR^d|w_i=w_j=0\}$ 
and its stem $\calS'$.
We consider the relationship between the stem and leaf nodes at the same depth 
and their parent nodes in this tree.
Let the parent nodes $(f_d,g_d)$ at depth $d$ be as follows $(d=0,1,\ldots,d(\calT)-1)$.

\[f_d(w)=w_i^{p_i}w_j^{p_j}
(w_{\backslash i,j}^*w_i^{r_i}+w_{\backslash i,j}^*w_j^{r'_j}),
\ g_d(w)=w_{\backslash i,j}^*w_i^{s_i-1}w_j^{s_j-1}\]
Suppose $r_i\geq r'_j$, we represent the stem $(f_{d+1},g_{d+1})$ at depth $d+1$ as follows.

\[f_{d+1}(w)=w_i^{p_i+p_j+r'_j}w_j^{p_j}
(w_{\backslash i,j}^*w_i^{r_i-r'_j}+w_{\backslash i,j}^*w_j^{r'_j}),
\ g_{d+1}(w)=w_{\backslash i,j}^*w_i^{s_i+s_j-1}w_j^{s_j-1}\]
The leaf $(f'_{d+1},g'_{d+1})$ at depth $d+1$ as follows.

\[f'_{d+1}(w)=w_i^{p_i}w_j^{p_j+p_i+r'_j}
(w_{\backslash i,j}^*w_i^{r_i}w_j^{r_i-r'_j}+w_{\backslash i,j}^*),
\ g'_{d+1}(w)=w_{\backslash i,j}^*w_i^{s_i-1}w_j^{s_j+s_i-1}\]
Comparing index ratios and potential index ratios, the following inequations hold.

\[
\left\{
\begin{array}{rcl}
\max\rho(f'_{d+1},g'_{d+1})&\leq&\max\{\max\rho(f_d,g_d),\max\rho(f_{d+1},g_{d+1})\}\\
\max\pi(f'_{d+1},g'_{d+1})&\leq&\max\{\max\pi(f_d,g_d),\max\pi(f_{d+1},g_{d+1})\}
\end{array}
\right.
\]
From these inequations and Lemma \ref{lmm_all_monotonous},
for root node $(f,g)$ and any leaf node $^\forall(f',g')\in\calL(\calT')$ of 
the blowup binary tree $\calT'$ obtained by 
blowup between variables (Algorithm \ref{alg:between_two_vars_with_jacobian}),
it is shown that the following inequation holds.

\begin{equation}\label{eq:compare_root_and_leaf}
\max\{\max\rho(f',g'),\max\pi(f',g')\}\leq\max\{\max\rho(f,g),\max\pi(f,g)\}
\end{equation}
Since blowup between terms (Algorithm \ref{alg:between_terms}) repeats blowup between variables 
with the leaves of a generated blow-up tree as the roots of a new blowup tree,
equation \eqref{eq:compare_root_and_leaf} also holds for the root node $(f,g)$ and 
any leaf node $^\forall(f',g')\in\calL(\calT)$ 
in the blowup binary tree $\calT$ obtained by blowup between terms.
Finally, we show that there exists a leaf node $(f',g')\in\calL(\calT)$ such that 
the equality of inequation \eqref{eq:compare_root_and_leaf}.
We consider by cases.

First, consider the case $\max\{\max\rho(f,g),\max\pi(f,g)\}=\max\rho(f,g)_i$$(i=1,2,\ldots,d)$.
Since the blow-up tree $\calT$ is binary, 
there exists a child node whose index of $w_i$ does not change at each branching by a blow-up.
Therefore, we continue to select such branches, then we reach a leaf node 
$(f',g')\in\calL(\calT)$ whose index of $w_i$ is the same as $(f,g)$ since tree $\calT$ is finite.
$\max\rho(f',g')=\rho(f',g')_i=\rho(f,g)_i$ holds there, 
so the equality of \eqref{eq:compare_root_and_leaf} holds.

Then, we consider the case $\max\{\max\rho(f,g),\max\pi(f,g)\}=\max\pi(f,g)_{i,j}$
$(i\in I_{(f,g)},j\in J_{(f,g)})$.
The algorithm runs a blowup between variables centered with $\{w\in\RR^d|w_i=w_j=0\}$ in the first while cycle.
Let $\calT'$ be the resulting binary tree in the first while cycle.
Note that $\calT'\subset\calT$ holds.
Let $(f',g')\in\calL(\calT')$ be the deepest leaf of the small binary tree $\calT'$($d((f',g'))=d(\calT')$).
Then either $\pi(f,g)_{i,j}=\rho(f',g')_i$ or $\pi(f,g)_{i,j}=\rho(f',g')_j$ holds.
From symmetry, it is sufficient to consider the case of $\pi(f,g)_{i,j}=\rho(f',g')_i$.
As in the previous discussion, there exists a child node whose index of $w_i$ does not change 
in subsequent branches. We continue to select such branches, 
then we reach a leaf node $(f'',g'')\in\calL(\calT)$ whose index of $w_i$ is equal to $(f',g')$.
$\max\rho(f'',g'')=\max\rho(f'',g'')_i=\max\rho(f',g')_i=\max\pi(f,g)_{i,j}$ holds there, 
thus the equality of \eqref{eq:compare_root_and_leaf} holds.

The above shows that there exists a leaf node $(f',g')\in\calL(\calT)$ such that 
the equality of \eqref{eq:compare_root_and_leaf} holds in all cases.
The theorem is now shown to hold since the following equation holds 
for the maximum pole $-\lambda$ of the zeta function and the blowup tree $\calT$.

\end{proof}
\end{thm}

\subsection{case of general sop polynomial}\label{sec:dn_case}
Finally, we consider the case of general sop polynomials.
In this case, unlike the binomial case, it is difficult to determine the normal crossing.
Therefore, we define the local normal crossing, which is looser than the normal crossing, 
then we decide whether or not to continue the blowup based on the given sop polynomial belonging to the set.
The following properties hold for local normal crossing.

\begin{prp}
Let $f$ be a non-negative polynomial with $d$-variables.
Let the multi-indexes matrix of $f$ be $A_f=\bmat{a_{*1}&a_{*2}&\cdots&a_{*n}}\in\RR^{d\times n}$.
$f\in\calLNC$ if and only if there exist $a_{*i}$ such that $a_{*i}\leq a_{*j}$ holds for any $a_{*j}$.

\begin{proof}
Suppose that $f(w)$ is a local normal crossing.
$f$ can be represented $f(w)=h(w)w^k$, where $h(w)$ is a polynomial 
that is not zero in a neighborhood of the origin and $k$ is a multi-indexes vector.
Since the polynomial $h$ is not zero at the origin, it has a constant term.
Therefore, let $(c_j)_{j=1}^n$ be the real sequence such that $c_j\neq0(^\forall j=1,2,\ldots,n)$ and 
$(a'_{*j})_{j=2}^n(a'_{*j}\in\NN_0^d,j=2,3,\ldots,n)$ be 
nonzero multi-indexes vectors that are different from each other, 
we represent $h(w)$ as $h(w)=c_1+\sum_{j=2}^nc_jw^{a'_{*j}}$.
From the above, we obtain $f(w)=c_1w^k+\sum_{j=2}^nc_jw_j^{a'_{*j}+k}$, 
so the multx-indexes matrix $A_f$ can be represented as $A_f=\bmat{k&a'_2+k&a'_3+k&\cdots&a'_n+k}$.
Therefore, $k\leq a'_j+k$ holds for any $j=2,\ldots,d$.

Conversely, let $A_f=\bmat{a_{*1}&a_{*2}&\cdots&a_{*n}}$ and suppose that 
there exists $a_{*i}$ such that $a_{*i}\leq a_{*j}$ for any $a_{*j}$.
It is ok to fix $i=1$.
Let $a_{*1}=k$ and $a'_{*j}=a_{*j}-k(^\forall j=2,3,\ldots,n)$, 
then we obtain a polynomial $h(w)$ satisfying $f(w)=h(w)w^k$ and $h(0)>0$, 
following the proof of the first part in reverse order.
Thus $f\in\calLNC$.

\end{proof}
\end{prp}

\subsubsection{blowup algorithm for general sop polynomial}\label{sec:dn_sop_alg}
In this subsubsection, we propose the blowup algorithm such that 
the inner polynomial of any leaf node of the blowup tree is 
a local normal crossing and show it halts.

\begin{prp}\label{prp:nc_blowup_halt}
Algorithm \ref{alg:local_nc_blowup} halts for any $f(w)\in\calSP^{d,n},g(w)\in\calM^d$.

\begin{proof}
Let $h=p_1+p_2+\cdots+p_n$ for $(h,k,n'_1,n'_2)$ added to $Q$, it is sufficient to show that 
$p_{\min\{n'_1,n'_2\}}$ is the minimum in set $\{p_1,\ldots,p_{\max\{n'_1,n'_2\}-1}\}$.
Where $p_i(i=1,2,\ldots,n)$ is a monomial and the minimum of monomials set $\{w^{a_1},w^{a_2},\ldots,w^{a_n}\}$ 
means the minimum of the multi-indexes vectors $\{a_1,a_2,\ldots,a_n\}$.
We show the proposition by induction.
First, $f_1=p_1$ is monomial for $(f,g,1,2)$ so obviously $p_1$ is the minimum in the set $\{p_1\}$.
Next, suppose that $p_{\min\{n'_1,n'_2\}}$ is the minimum in the set $\{p_1,\ldots,p_{\max\{n'_1,n'_2\}-1}\}$ 
for $(h,k,n'_1,n'_2)$, we show $p'_{\min\{n'_1,\max\{n'_1,n'_2\}+1\}}=p'_{n'_1}$ is the minimum in the set 
$\{p'_1,\ldots,p'_{\max\{n'_1,\max\{n'_1,n'_2\}+1\}-1}\}$ 
for $(h\circ h'',k',n'_1,\max\{n'_1,n'_2\}+1)$ if $p'_{n'_1}\leq p'_{n'_2}$.
Since $p_{\min\{n'_1,n'_2\}}$ is the minimum in the set $\{p_1,\ldots,p_{\max\{n'_1,n'_2\}-1}\}$,
$p'_{\min\{n'_1,n'_2\}}$ is also the minimum in the set $\{p'_1,\ldots,p'_{\max\{n'_1,n'_2\}-1}\}$ 
after blowing up between terms.
If $\min\{n'_1,n'_2\}=n'_1$, because of $p'_{n'_1}\leq p'_{n'_2}$ by assumption, 
$p'_{n'_1}$ is the minimum in the set $\{p'_1,\ldots,p'_{\max\{n'_1,n'_2\}}\}$.
Next, if $\min\{n'_1,n'_2\}=n'_2$, because of $p'_{n'_1}\leq p'_{n'_2}$ by assumption, 
the chain rule derives that 
$p'_{n'_1}$ is the minimum in the set $\{p'_1,\ldots,p'_{n'_2},\ldots,p'_{\max\{n'_1,n'_2\}-1},p'_{n'_1}\}$.
Furthermore, since $n'_1=\max\{n'_1,n'_2\}$, $p'_{n'_1}$ is the minimum 
in the set $\{p'_1,\ldots,p'_{\max\{n'_1,n'_2\}}\}$.
From the above discussion, we show all the cases.
The same discussion holds for the case of $p'_{n'_2}\leq p'_{n'_1}$.
Thus the proposition is proved.

\end{proof}
\end{prp}

From Proposition \ref{prp:nc_blowup_halt} and Theorem \ref{thm:terms_blowup_halt}
derives the main theorem \ref{thm:main_theorem_blowup_algorithm}.

\subsubsection{RLCT of general sop polynomial}

We then consider the RLCT of the general sop polynomial $\calSP^{d,n}$.
As described in the previous section \ref{sec:dn_sop_alg}, 
it is difficult to derive the exact value of the RLCT for general sop polynomials 
because it is also difficult to judge that they are normal crossings.
In this paper, we derive an upper bound of the RLCT instead of the exact value.
First, we define the concepts necessary for the proof.

\begin{dfn}[index ratio at the origin]
Let $f(w)\in\calSP^{d,n}_{\geq0}\cap\calLNC$ and $g(w)\in\calM^d$, 
and they are represented as follows.

\[f(w)=w^p\left(1+\sum_{j=2}^nw^{r_j}\right)\ ,\ g(w)=w^{s-1}\]
Where $p,r_j\in\NN_0^d(j=2,3,\ldots,n),s\in\RR^d_{>0}$,$r_j\neq0(j=2,3,\ldots,n)$,
$r_j\neq r_{j'}(j\neq j')$, and $p+\sum_{j=2}^nr_j\in\NN^d$.
We define the vector for multi-indexes as follows.

\[\rho(f,g)=\frac{p}{s}=\left(\frac{p_i}{s_i}\right)_{i=1}^d\]
The components of the vector are referred to as \textbf{index ratios at the origin}.

\end{dfn}

The following proposition holds for index ratios at the origin.

\begin{prp}\label{prp:local_nc_upperbound}
Let $f(w)\in\calSP^{d,n}_{\geq0}$ and $g(w)\in\calM^d$.
Let the finite blowup tree whose local coordinates are all local normal crossings and 
the root node is $(f,g)$ be $\calT$.
Where the centers of the blowup are limited to the non-singular set around the axes $C_T\in\calC$.
For the maximum pole $-\lambda$ of the zeta function defined below

\[\zeta(z)=\int_Wf(w)^z|g(w)|dw,\]
the following inequation holds.

\[\frac{1}{\lambda}\geq\underset{(h,k)\in\calL(\calT)}{\max}\left\{\max\rho(h,k)\right\}\]
Where the integral domain $W$ is a sufficiently large compact set.

\begin{proof}
It is ok to replace the integral domain $W$ with $[-D,D]^d$ for a sufficiently large number $D>0$.
For the multi-indexes matrix $A_f,A_h$ of the root $f$ and a leaf $h$, 
there exists a blow matrix $B_h\in\RR^{d\times d}$ such that $A_h=B_h^\top A_f$.
Let $w\in\RR^d$ be the variable at root $f$ and $w'\in\RR^d$ be the variable at leaf $h$.
Then $B_h\log w'=\log w$ holds.
Therefore, the radius of integration after variable transformation for each variable 
can be represented as $B_h^{-1}1_d(\log D)$.
Let $\delta^{(h)}=\exp(B_h^{-1}1_d\log D)\in\RR^d$, 
then the zeta function can be transformed as follows.

\[\zeta(z)=\sum_{(h,k)\in\calL(\calT)}\int_{[\pm\delta^{(h)}]}h(w)^z|k(w)|dw\]
Where $[\pm\delta^{(h)}]=\prod_{j=1}^d[-\delta^{(h)}_j,\delta^{(h)}_j]$ for the integral domain.
Since $h(w)$ is normal crossing in a neighborhood of the origin, 
let $A_h=(a^{(h)}_{ij})\in\RR^{d\times n}$ and define 
$h'(w)=h(w)/\prod_{i=1}^dw_i^{\underset{j}{\min}\left(a^{(h)}_{ij}\right)}$, then 
we get $h'(0)>0$.
The following minute amount $\epsilon_h>0$ is defined as follows.

\[
\epsilon_h=\frac{1}{2}\min\left\{\underset{j}{\min}\left(\delta^{(h)}_j\right),
\underset{h'(w)=0}{\min}\|w\|_\infty\right\}
\]
Since it is obvious that $1_d\epsilon_h<\delta^{(h)}$ holds, 
from the property of the inequality (Proposition \ref{eq:rlct_prp:inequality}), 
let the maximum pole of the following zeta function be $-\lambda'$,

\[
\zeta'(z)=\sum_{(h,k)\in\calL(\calT)}\int_{[-\epsilon_h,\epsilon_h]^d}h(w)^z|k(w)|dw
\]
$\lambda'\geq\lambda$ holds. Since $h(w)$ is normal crossing on $[-\epsilon_h,\epsilon_h]^d$,

\[\frac{1}{\lambda'}=\underset{(h,k)\in\calL(\calT)}{\max}\left\{\rho(h,k)\right\}\]
holds. Thus the proposition is proved.

\end{proof}
\end{prp}
With the above preparation, we prove the following theorem.

\begin{thm}\label{thm:dn_rlct_upper_bound}
Let $f\in\calSP^{d,n}_{\geq0}$,$g\in\calM^d$ and they are represented as follows.

\[f(w)=\sum_{j=1}^nw^{a_{*j}},\ g(w)=w^{s-1}\]
Where $a_{*j}\in\NN_0^d(j=1,2,\ldots,n)$,$s\in\RR^d,s>0$.
we define the following quantity for the multi-indexes.

\begin{equation}\label{eq:dn_smplx_ub}
\frac{1}{\smplx}=\underset{\alpha\in S_{d-1}}{\max}\left\{
\underset{j=1,\ldots,n}{\min}\left(\sum_{h=1}^d\alpha_h\frac{a_{hj}}{s_h}\right)\right\}
\end{equation}
Then, the following inequation holds for the maximum pole $-\lambda$ of 
the zeta function $\zeta(z)=\int_Wf(w)^z|g(w)|dw$.

\[\lambda\leq\smplx\]
Where $S_{d-1}$ is $d-1$-dimensional simplex, that is,

\[S_{d-1}=\left\{\alpha\in\RR^d\middle|\alpha_h\geq0, \sum_{h=1}^d\alpha_h=1\right\}.\]
The integral domain $W$ is a sufficiently large compact set.

\begin{proof}
Let $(f',g')\in\calL(\calT)$ be a leaf of a blowup tree $\calT$ whose root is $(f,g)$.
Where the centers of the blowups of this tree $\calT$ are limited to 
a non-singular set around the coordinate axes $C_T\in\calC$.
For any leaf $(f',g')\in\calL(\calT)$, 
there exists a blow matrix $B=(b_{hi})_{h,i=1}^d\in\NN_0^{d\times d}$ such that $A_{f'}=B^\top A_f$.
Then the following holds.

\[\max\rho(f',g')=\underset{i=1,2,\ldots,d}{\max}\left(\underset{j=1,2,\ldots,n}{\min}
\frac{\sum_{h=1}^db_{ih}a_{hj}}{\sum_{h=1}^db_{ih}s_h}\right)\]
From Proposition \ref{prp:local_nc_upperbound},

\[\frac{1}{\lambda}\geq\underset{(f',g')\in\calL(\calT)}{\max}
\left\{\underset{i=1,2,\ldots,d}{\max}\left(\underset{j=1,2,\ldots,n}{\min}
\frac{\sum_{h=1}^db_{ih}a_{hj}}{\sum_{h=1}^db_{ih}s_h}\right)\right\}\]
holds. In the following, $t_{i,(f',g')}=\underset{j=1,2,\ldots,n}{\min}
\left(\sum_{h=1}^db_{ih}a_{hj}\right)/\left(\sum_{h=1}^db_{ih}s_h\right)$ for the simplicity.
For $t_{i,(f',g')}$

\[\underset{(f',g')\in\calL}{\max}\left(\underset{i=1,\ldots,d}{\max}\ t_{i,(f',g')}\right),
\underset{i=1,\ldots,d}{\max}\left(\underset{(f',g')\in\calL}{\max}\ t_{i,(f',g')}\right)\geq t_{i,(f',g')}\]
holds, so we can reverse the order of taking the maximum as follows.

\begin{equation}\label{eq:swap_max}
\frac{1}{\lambda}\geq\underset{i=1,\ldots,d}{\max}\left\{
\underset{(f',g')\in\calL(\calT)}{\max}\left(\underset{j=1,\ldots,n}{\min}
\frac{\sum_{h=1}^db_{ih}a_{hj}}{\sum_{h=1}^db_{ih}s_h}\right)\right\}
\end{equation}
We now prove the following lemma.

\begin{lmm}\label{lmm:arbitary_tree_given_gcd1}
Let $b=(b_h)_{h=1}^d\in\NN_0^d$ be a vector satisfying $\underset{h=1,\ldots,d}{\gcd}(b_h)=1$.
There is a blowup tree $\calT$ whose leaves are all local normal crossings and such that 
one of the blow matrix $B\in\RR^{d\times d}$ corresponding to its leaves $(f',g')$ has $b$ as a column.

\begin{proof}
First, we construct matrix $B'$ which has $b$ as a column by basic transformation from the identity matrix.
Since $\gcd(b_1,\ldots,b_d)=1$, there is a component $i$ such that $b_i\notin2\ZZ$.
we multiply $i$-th column by $b_i$ for the identity matrix at first.
Then we carry the diagonal components of $h$-th component ($h\neq i$) to 
$i$-th column for $b_h$ times by the basic transformation.
That is, let $P_i(b_i)=\bmat{I_{i-1}&&\\&b_i&\\&&I_{d-i}}$, we construct $B'$ as follows.

\[
B'\!=\!\bmat{e_1&\cdots&e_{i-1}&b&e_{i+1}&\cdots&e_d}\!=\!\left\{
\prod_{h=1,h\neq i}^d\left(R_{hi}^{b_h}\right)\right\}
P_i(b_i)
\]

Next, we construct a blowup tree $\calT$ with leaf $(f',g')\in\calL(\calT)$ 
corresponding to the blowup matrix $B$ whose $i$-th column is the same as $B'$.
There is a matrix $B''$ such that $B''=\prod_{k=1}^nR_{i_k,j_k}$ and $B'=B''P_i(b_i)$.
First of all, variable transform $\log w=P_i(b_i)\log w^{(0)}$.
Let $(f(\exp(P_i(b_i)\log w^{(0)})),g(\exp(P_i(b_i)\log w^{(0)})))
=(f_0,g_0)$ obtained by the variable transformation be the root node of tree $\calT$
We blowup centered with $\{w^{(0)}_{i_1}=w^{(0)}_{j_1}=0\}\in\calC$, 
then we obtain two variable transformations.
Let $(f_1,g_1)$ be the node obtained by variable transformation 
such that $\log w^{(0)}=R_{i_1,j_1}\log w^{(1)}$and $(f'_1,g'_1)$ be the node 
obtained by variable transformation such that $\log w^{(0)}=R_{j_1,i_1}\log {w'}^{(1)}$.
In the same way, for node $(f_k,g_k)$, blowup centered with 
$\{w^{(k)}_{i_{k+1}}=w^{(k)}_{j_{k+1}}=0\}\in\calC$, 
then let $(f_{k+1},g_{k+1})$ be the node obtained by variable transformation such that 
$\log w^{(k)}=R_{i_{k+1},j_{k+1}}\log w^{(k+1)}$ and $(f'_{k+1},g'_{k+1})$ be the node obtained by 
variable transformation $\log w^{(k)}=R_{j_{k+1},i_{k+1}}\log {w'}^{(k+1)}$ recursively.
Let the small blowup tree obtained by this operation be $\calT'$. 
For any multi-indexes matrix of the inner polynomial of the node,

\[A_{f_l}=\left(\prod_{k=1}^lR_{i_k,j_k}\right)^\top P_i(b_i)A_f\ \ \ (l=1,2,\ldots,n)\]
\[A_{f'_l}=R_{j_l,i_l}^\top\left(\prod_{k=1}^{l-1}R_{i_k,j_k}\right)^\top P_i(b_i)A_f\ \ \ (l=1,2,\ldots,n)\]
hold. Note that $A_{f_n}={B'}^\top A_f$ especially.
We define the stem of small blowup binary tree $\calT'$ as $\calS=\{(f_l,g_l)\}_{l=0}^n$.
The leaves $\calL(\calT)$ of $\calT$ are $\calL(\calT')=\{(f'_l,g'_l)\}_{l=0}^n\cup\{(f_n,g_n)\}$.
Let $\calT$ be the blowup tree extended by applying the local normal crossing blowup 
(Algorithm \ref{alg:local_nc_blowup}) to each leaf of this tree $\calT'$.
Since local canonical crossing blowup halts 
(Main Theorem \ref{thm:main_theorem_blowup_algorithm} and Proposition \ref{prp:nc_blowup_halt}), 
$\calL(\calT)\subset\calLNC$ holds for the leaves set of this tree $\calT$.
Finally, we show that this tree $\calT$ has a leaf $(f',g')\in\calL(\calT)$ corresponding to 
the blow matrix $B$ whose $i$-th column is the same as $B'$, that is, 
it has a leaf $(f',g')$ such that $A_{f'}=B^\top A_f$.
The relation $\log w=B'\log w^{(n)}$ holds for variable $w\in\RR^d$ at the root and 
variable $w^{(n)}\in\RR^d$ at local coordinates corresponding to $(f_n,g_n)$.
Consider the following leaves extended from the leaf $(f_n,g_n)\in\calL(\mathcal{T}')$ 
of a small blow-up binary tree $\calT'$.
There exists a branch where the blowup does not change the degree of the $i$-th variable of the polynomial.
Let $(f',g')$ be the leaf at which the branch ends and $B$ be the blowup matrix.
There exists some $C\in\calK$ such that $B=B'C$.
Since $A_{f'}=C^\top{B'}^\top A_f=C^\top A_{f_n}$ and 
the $i$-th row of $A_{f'}$ is equivalent to the $i$-th row of $A_{f_n}$,
the $i$-th column of $C^\top$ is $e_i^\top$($e_i\in\RR^d$ is a standard basis).
That is, $i$-th column of $B$ and that of $B'$ are equal because $i$-th column of $C$ is $e_i$.
Thus, the lemma is proved.

\end{proof}
\end{lmm}
(The proof of the theorem \ref{thm:dn_rlct_upper_bound} continues below.)
Suppose that there exists the non-negative integer vector $(b_h)_{h=1}^d$ such that $\gcd(b)=1$ and 
maximizes $\underset{j}{\min}\left(\sum_{h=1}^db_ha_{hj}\right)/\left(\sum_{h=1}^db_hs_h\right)$.
That is, suppose that there exists the vector which archives the following.

\begin{equation}\label{eq:discrete_form}
\underset{\gcd(b_1,\ldots,b_d)=1}{\sup}\left\{\underset{j=1,\ldots,n}{\min}
\left(\frac{\sum_{h=1}^db_ha_{hj}}{\sum_{h=1}^db_hs_h}\right)\right\}
\end{equation}
Since $b\neq0$ and Lemma \ref{lmm:arbitary_tree_given_gcd1}, 
there exists the blowup tree $\calT^*$ which has 
a leaf whose index ratios at the origin as \eqref{eq:discrete_form}.
$t_{i,(f',g')}\leq\eqref{eq:discrete_form}$ by the definition, 
then substituting $\calT^*$ for $\calT^*$ in \eqref{eq:swap_max}, we obtain the following inequation.

\[\frac{1}{\lambda}\geq\underset{\gcd(b_1,\ldots,b_d)=1}{\sup}\left\{\underset{j=1,\ldots,n}{\min}
\left(\frac{\sum_{h=1}^db_ha_{hj}}{\sum_{h=1}^db_hs_h}\right)\right\}\]

We consider the continuous optimization problem \eqref{eq:dn_smplx_ub} to prove 
the existence of $b\in\NN_0^d$ which achives \eqref{eq:discrete_form}.
In the following, first we show \eqref{eq:discrete_form} is equivalent to \eqref{eq:dn_smplx_ub},
and derive such $b\in\NN_0^d$.

We show \eqref{eq:discrete_form}$\leq$\eqref{eq:dn_smplx_ub} at first.
Let $(b_h)_{h=1}^d$ be the non-negative integer vector such that $\gcd(b)=1$.
Note that this is not a zero vector since $\gcd(b)=1$.
We define the following non-negative vector.

\[\alpha_h=\frac{b_hs_h}{\sum_{h=1}^db_hs_h}\ (^\forall h=1,\ldots,d)\]
The following inequation holds.

\[\frac{\sum_{h=1}^db_ha_{hj}}{\sum_{h=1}^db_hs_h}=\sum_{h=1}^d\alpha_h\frac{a_{hj}}{s_h}
\leq\underset{\alpha\in S_{d-1}}{\max}\left\{
\underset{j=1,\ldots,n}{\min}\left(\sum_{h=1}^d\alpha_h\frac{a_{hj}}{s_h}\right)\right\}
\]
Where $(b_h)_{h=1}^d$ is arbitrary, 
\eqref{eq:discrete_form}$\leq$\eqref{eq:dn_smplx_ub} holds by taking the upper bound of left-hand side.

Next, we show \eqref{eq:dn_smplx_ub}$\leq$\eqref{eq:discrete_form}.
We can formularize \eqref{eq:dn_smplx_ub} as the following linear programming problem (
In the following, we denote it as LPP for short).

\begin{align*}
\mathrm{maximize}:&\beta\\
\mathrm{subject\ to}:&\\
&\alpha_1\mu_{11}+\alpha_2\mu_{21}+\cdots+\alpha_d\mu_{d1}-\beta\geq0\\
&\alpha_1\mu_{12}+\alpha_2\mu_{22}+\cdots+\alpha_d\mu_{d2}-\beta\geq0\\
&\vdots\\
&\alpha_1\mu_{1n}+\alpha_2\mu_{2n}+\cdots+\alpha_d\mu_{dn}-\beta\geq0\\
&\alpha_1+\alpha_2+\cdots+\alpha_d=1\\
&\alpha_1,\alpha_2,\ldots,\alpha_d,\beta\geq0
\end{align*}
Where $M=(\mu_{hj})=(\frac{a_{hj}}{s_h})_{h=1,j=1}^{d,n}\in\RR^{d\times n}$.
The above LPP can be transformed into the standard LPP 
by introducing slack variables $\gamma=(\gamma_1,\ldots,\gamma_n)^\top\in\RR^n$.

\begin{align}
\mathrm{maximize}:&(0_d^\top,1,0_n^\top)(\alpha^\top,\beta,\gamma^\top)^\top\nonumber\\
\mathrm{subject\ to}:&\nonumber\\
&\bmat{-M^\top&1_n&I_n\\1_d^\top&0&0_n^\top}\bmat{\alpha\\ \beta\\ \gamma}=e_{n+1}
\label{eq:regular_form_LP}\\
&\alpha,\beta,\gamma\geq0\nonumber
\end{align}
Where $e_{n+1}=\bmat{0_n\\1}$, the standard basis in $n+1$-dimensional vector space.
This LPP has an optimal (basis) solution since it is obviously bounded and has a feasible solution.
Let $D\in\mathbb{R}^{(n+1)\times(d+1+n)}$ be the coefficient matrix 
of the conditional equation \eqref{eq:regular_form_LP}.
We now consider the principle of the simplex method.
In the simplex method, we divide variables into $n+1$ basis variables and $d$ non-basis variables 
whose values are fixed at $0$. Let $B\in\RR^{(n+1)\times(n+1)}$ be the submatrix of the coefficient matrix $D$, 
that is, $B$ composes of $n+1$ columns selected from $D$.
Let $B$ be reversible.
If the solution $y\in\RR^{n+1}$ to the linear equation $By=\bmat{0_n\\1}$ satisfies $y\geq0$, 
$y$ is a feasible basis solution of the simplex method.
There exists an optimal basis solution that achieves the optimal value of LPP.
Let $D=\bmat{d_{*1}&\cdots&d_{*d+1+n}}$ and $B=\bmat{b_{*1}&\cdots&b_{*n+1}}$.
And let the map between indexes such that $d_{k_h}=b_h(h=1,\ldots,n+1)$ be $k:h\mapsto k_h$.
It is ok to assume that $k$ is monotonically increasing.
When $B$ is reversible, that is , when $B$ is a basis matrix, Cramer's formula gives

\begin{equation}\label{eq:cramel_formula}
y_i=\frac{\det(\bmat{b_{*1}&\cdots&b_{*i-1}&e_{n+1}&b_{*i+1}&\cdots&b_{*n+1}})}{\det(B)}.
\end{equation}
By expanding cofactor for the numerator of \eqref{eq:cramel_formula}, 
since $i$-th column is $e_{n+1}$, we obtain

\[y_i\det(B)=(-1)^{i+n+1}\det(B'_{\backslash i}).\]
Where $B'_{\backslash i}=
\bmat{b'_{*1}&\cdots&b'_{*i-1}&b'_{*i+1}&\cdots&b'_{*n+1}}\in\RR^{n\times n}$ and 
$B'=\bmat{b'_{*1}&\cdots&b'_{*n+1}}\in\RR^{n\times(n+1)}$.
From the definition, $B'$ is submatrix of $\bmat{-M^\top&1_n&I_n}\in\RR^{n\times(d+1+n)}$.
Thus $B'$ can be divided into the columns from $-M^\top$ and the columns from $\bmat{1_n&I_n}$.
Since it is supposed that $B$ is reversible, $1\leq I\leq n+1$ holds.
Then we obtain $-s_{k_h}b'_{*h}=a_{k_h*}^\top\in\NN_0^n(h=1,\ldots,I)$.
Where we decompose the multi-indexes matrix $A_f$ into $\bmat{a_{1*}\\ \vdots\\a_{d*}}\in\RR^{d\times n}$.
From the above discussion, the determinant of $B'_{\backslash i}$ ($i=1,2,\ldots,I$) 
is represented as follows by multi-linearity.

\[\det(B'_{\backslash i})=\frac{1}{\prod_{h=1,h\neq i}^I(-s_{k_h})}
\sum_{\sigma\in\mathfrak{S}_n}\sign(\sigma)\left(\prod_{h=1,h\neq i}^Ia_{k_h,\sigma(k_h)}\right)
\left(\prod_{h=I+1}^Ip_{k_h,\sigma(k_h)}\right)\]
Where $\mathfrak{S}_n$ is $n$-dimensional symmetric group and $p_{k_h,\sigma(k_h)}$ is $0$ or $1$.
Since each element of the summation $\sigma\in\mathfrak{S}_n$ is an integer,

\[\det(B'_{\backslash i})\left(\prod_{h=1,h\neq i}^Is_{k_h}\right)\in\ZZ\ \ \ \mathrm{if}\ i=1,2,\ldots,I.\]
Also $y_i=\alpha_{k_i}$ for $i=1,2,\ldots,I$,

\[\frac{\alpha_{k_i}}{s_{k_i}}\det(B)\left(\prod_{h=1}^Is_{k_h}\right)\in\ZZ\ (^\forall i=1,2,\ldots,I)\]
holds. Thus $\left(\frac{\alpha_{k_h}}{s_{k_h}}\right)_{h=1}^I$ are integer ratios, and 
$\left(\frac{\alpha_h}{s_h}\right)_{h=1}^d$ are also integer ratios with non-basis variables.

The above discussion also holds for any feasible basis solution, which includes the optimal solution.
Let the ratios which satisfy the optimal solution be $(\alpha_h^*)_{h=1}^d$.
Since we can represent it as $\frac{\alpha_h^*}{s_h}=kq_h(q_h\in\NN_0)$ for some $k$,

\[\sum_{h=1}^d\alpha_h^*\frac{x_{hj}}{s_h}=k\sum_{h=1}^dq_hx_{hj}.\]
Where $\sum_{h=1}^d\alpha_h^*=1$, so $k$ satisfies the following equation.

\[\sum_{h=1}^d\alpha_h^*=k\sum_{h=1}^ds_hq_h=1\]
Thus we obtain

\[\sum_{h=1}^d\alpha_h^*\frac{a_{hj}}{s_h}=\frac{\sum_{h=1}^dq_hx_{hj}}{\sum_{h=1}^dq_hs_h}.\]
we can let $(q_h)_{h=1}^d$ be prime each other by dividing the great common divisor.
This is $(b_h)_{i=1}^d$, so \eqref{eq:dn_smplx_ub}$\leq$\eqref{eq:discrete_form} and 
the existence of $b\in\NN_0^d$ which archives \eqref{eq:discrete_form} are proved.

\end{proof}
\end{thm}
In the proof of Theorem \ref{thm:dn_rlct_upper_bound}, 
we constructed a blowup tree that clearly has local coordinates with 
the pole satisfying the claim of the theorem by using local normal crossing blowup 
(Algorithm \ref{alg:local_nc_blowup}) supplementarily.
Some numerical experiments show that the blowup trees 
obtained by using local normal crossing blowup directly for the target polynomial 
also have local coordinates with the pole satisfying the claim of the theorem.
This is expected to be provable.

(The proof of the main theorem \ref{thm:main_theorem_rlct})
Finally, we show that the upper bound \eqref{eq:dn_smplx_ub} coincides with the result 
(Theorem \ref{thm:nonNC_d2sop_rlct}) in the binomial case.
The upper bound \eqref{eq:dn_smplx_ub} is the inverse of the maximum value of $\beta$ 
in the feasible basis solution of the LPP \eqref{eq:regular_form_LP}.
A feasible basis solution requires that at least one of the components of $\alpha$ be positive.
First, consider a feasible basis solution such that $\alpha_h=1$ and $\beta>0$.
In this case, $\beta=\min(\mu_{h1},\mu_{h2})$ and 
$\gamma_{\underset{n=1,2}{\arg\min}(\mu_{hn})}=\max(\mu_{h1},\mu_{h2})-\min(\mu_{h1},\mu_{h2})$ 
is a feasible basis solution.
Next, consider feasible basis solutions such that $\alpha_h,\alpha_{h'}>0$ and $\beta>0$.
It is sufficient to solve the following simultaneous equation.

\[\left\{\begin{array}{rcl}
\alpha_h\mu_{h1}+\alpha_{h'}\mu_{h'1}&=&\alpha_h\mu_{h2}+\alpha_{h'}\mu_{h'2}\\
\alpha_h+\alpha_{h'}&=&1
\end{array}\right.\]
$\alpha_h,\alpha_{h'}>0$ if and only if $(\mu_{h1}-\mu_{h2})(\mu_{h'2}-\mu_{h'1})>0$.
Solving this, we obtain $\alpha_h=(\mu_{h'2}-\mu_{h'1})/(\mu_{h1}-\mu_{h2}+\mu_{h'2}-\mu_{h'1})$.
Thus, the following basis solution

\[\beta=\frac{\mu_{h1}(\mu_{h'2}-\mu_{h'1})+\mu_{h'1}(\mu_{h1}-\mu_{h2})}{\mu_{h1}-\mu_{h2}+\mu_{h'2}-\mu_{h'1}}
=\frac{\mu_{h1}\mu_{h'2}-\mu_{h'1}\mu_{h2}}{\mu_{h1}-\mu_{h2}+\mu_{h'2}-\mu_{h'1}}\]
is feasible. The other feasible basis solution derives $\beta=0$.
From the above discussion, let $\nu_h=\max(\mu_{h1},\mu_{h2})$,$\mu_h=\min(\mu_{h1},\mu_{h2})$ 
for $h=1,2,\ldots,d$, and define the index set $H=\{(h,h')|h\neq h',\nu_h-\mu_h,\nu_{h'}-\mu_{h'}>0\}$,
we obtain

\[\underset{\alpha\in S_{d-1}}{\max}\left\{\underset{j=1,2}{\min}\left(
\sum_{h=1}^2\alpha_h\frac{a_{hj}}{s_h}\right)\right\}=\max\left\{
\underset{h}{\max}(\mu_h),\ \underset{(h,h')\in H}{\max}\left(
\frac{\nu_h\nu_{h'}-\mu_h\mu_{h'}}{\nu_h+\nu_{h'}-\mu_h-\mu_{h'}}
\right)\right\}.\]
Therefore the main theorem \ref{thm:main_theorem_rlct} is proved.

\section{Application and discussion of main theorem}
In this section, we discuss applications of the main theorem, in particular, 
the results of the general sop polynomial (Section \ref{sec:dn_case}).
One is the derivation of the upper bound of the RLCT for general polynomials and 
the other is the derivation of the optimal weights in the weighted blowup.

\subsection{simplex upper bound}\label{sec:simplex_upper_bound}

The result of Theorem \ref{thm:dn_rlct_upper_bound} is independent of the coefficients of the polynomial, 
as the proof of the theorem shows.
Therefore, it can be applied to general polynomials.
In other words, the following result holds.

\begin{cor}[simplex upper bound]\label{cor:simplex_upper_bound}
Let $f\in\RR[w_1,w_2,\ldots,w_d]$ be a polynomial and $g\in\calM^d$ be a monomial, 
and they can be represented as $f(w)=\sum_{j=1}^nc_jw^{a_{*j}}$ and $g(w)=w^{s-1}$.
Where $c_j\neq0(^\forall j)$,$a_{*j}\neq a_{*j'}(j\neq j')$ and $s\in\RR^d_{>0}$.
We define the following quantity for the multi-indexes.

\begin{equation}\label{eq:simplex_upper_bound}
\frac{1}{\smplx}=\underset{\alpha\in S_{d-1}}{\max}\left\{\underset{j=1,2,\ldots,n}{\min}
\left(\sum_{h=1}^d\alpha_h\frac{a_{hj}}{s_h}\right)\right\}
\end{equation}
If right-hand side of \eqref{eq:simplex_upper_bound} is equal to $0$, let $\smplx=\infty$.
The following holds for the maximum pole $-\lambda$ of the zeta function $\zeta(z)=\int_W|f(w)|^z|g(w)|dw$.

\[\lambda\leq\smplx\]
Where the integral domain $W$ is a sufficiently large compact set.
$\smplx$ is referred to as \textbf{simplex upper bound} of RLCT $\lambda$.

\end{cor}

The simplex upper bound $\smplx$ can be obtained by formularizing the following LPP, 
as described in the proof of Theorem \ref{thm:dn_rlct_upper_bound}.
The inverse of the objective function $\beta$ corresponds to $\smplx$.

\begin{eqnarray}
\mathrm{maximize}:&\beta&\nonumber\\
\mathrm{subject\ to}:&&\nonumber\\
&&\alpha_1\frac{a_{11}}{s_1}+\alpha_2\frac{a_{21}}{s_2}+\cdots+\alpha_d\frac{a_{d1}}{s_d}-\beta\geq0\nonumber\\
&&\alpha_1\frac{a_{12}}{s_1}+\alpha_2\frac{a_{22}}{s_2}+\cdots+\alpha_d\frac{a_{d2}}{s_d}-\beta\geq0\nonumber\\
&&\vdots\label{eq:mi_poly_lp}\\
&&\alpha_1\frac{a_{1n}}{s_1}+\alpha_2\frac{a_{2n}}{s_2}+\cdots+\alpha_d\frac{a_{dn}}{s_d}-\beta\geq0\nonumber\\
&&\alpha_1+\alpha_2+\cdots+\alpha_d=1\nonumber\\
&&\alpha_1,\alpha_2,\ldots,\alpha_d,\beta\geq0\nonumber
\end{eqnarray}

The LPP \eqref{eq:mi_poly_lp} is called a 
\textbf{linear programming problem based on multi-indexes of polynomials}.
The computational complexity of obtaining the simplex upper bound depends on the method of solving 
the LPP.
The simplex method is a typical method for solving LPP, 
but there are other methods such as the interior point method \cite{IPM}.
Both of these methods can have huge computational complexity for problems with complex conditions, 
but otherwise, they are known to have relatively tiny computational complexity.

\subsubsection{property of simplex upper bound}

Next, we discuss the properties of the simplex upper bound.

\begin{prp}\label{prp:simplex_infinity_if_nonzero}
let $f\in\RR[w_1,w_2,\ldots,w_d]$.
$f(0)\neq 0$ is equivalent to $\smplx=\infty$.

\begin{proof}
Suppose $f(0)\neq 0$. When we expand $f(w)=\sum_{j=1}^nc_jw^{a_{*j}}$, 
there exists term such that $c_j=f(0),a_{*j}=0\in\NN_0^d$. Thus,

\begin{equation}\label{eq:min_is_zero}
\underset{j=1,2,\ldots,n}{\min}\left(\sum_{h=1}^d\alpha_h\frac{a_{hj}}{s_h}\right)=
\left(\sum_{h=1}^d\alpha_h\frac{0}{s_h}\right)=0
\end{equation}
From the above discussion, $\smplx=\infty$ holds.

Conversely, suppose $\smplx=\infty$.
\eqref{eq:min_is_zero} holds for any $\alpha\in S_{d-1}$, 
so there exists $j'\in\{1,2,\ldots,n\}$ such that $a_{*j'}=0,c_{j'}\neq0$.
Thus,

\[f(0)=\sum_{j=1}^nc_jw^{a_{*j}}=c_{j'}0^0=c_{j'}\neq0\]
holds. The proposition is proved.

\end{proof}
\end{prp}
From Proposition \ref{prp:simplex_infinity_if_nonzero}, 
the polynomial $f$ at the origin must be such that $f(0)=0$ 
to obtain a valid upper bound of RLCT.
Therefore, if $f(0)\neq 0$, we first consider parallel translation.
For simplicity, we assume that $g(w)=1(^\forall w\in W)$.
Suppose that the point $p\in\RR^d$ satisfies $f(p)=0$
Let $f'(w')=f(w'+p)$ be the polynomial obtained by the parallel translation $w'=w+p$
and the maximum pole of zeta function $\zeta'(z)=f'(w')^zdw'$ be $-\lambda'$,

\[\lambda'=\lambda\]
holds.
On the other hand, since $f'(0)=f(p)=0$, 
the simplex upper bound $\smplx$ is bounded by Proposition \ref{prp:simplex_infinity_if_nonzero}.
Therefore, we can obtain an effective upper bound $\lambda=\lambda'\leq\lambda'_{\mathrm{smplx}}$ 
for the RLCT.
We denote this by $\smplx^{(p)}$ hereafter and referred to as simplex upper bound at $p$.
Unless mentioned, $\smplx$ simply means the simplex upper bound $\smplx^{(0)}$ at the origin.
The simplex upper bound differs at the point. The following is an example.

\begin{xmp}
Let $f\in\RR[w_1,w_2]_{\geq0}$ be defined as follows.

\[f(w_1,w_2)=w_1^4w_2^2\]
Consider the simplex upper bound of RLCT at $(p_1,p_2)\in\RR^2$.

\begin{eqnarray*}
f(w_1+p_1,w_2+p_2)&=&(w_1+p_1)^4(w_2+p_2)^2\\
&=&(w_1^4+4p_1w_1^3+6p_1^2w_1^2+4p_1^3w_1+p_1^4)(w_2^2+2p_2w_2+p_2^2)
\end{eqnarray*}
From the above calculation, 
if $p_1p_2\neq0$, $\smplx^{(p_1,p_2)}=\infty$ from Proposition \ref{prp:simplex_infinity_if_nonzero}.
If $p_1\neq0,p_2=0$, then $\smplx^{(p_1,0)}=\frac{1}{2}$.
If $p_1=0$, we obtain $\smplx^{(0,p_2)}=\frac{1}{4}$.
\end{xmp}
Similarly to the parallel translation, the simplex upper bound after linear transformation 
gives the upper bound of the RLCT.
Let the constant matrix $P\in\RR^{d\times d}$ satisfy $\det(P)\neq0$.
The polynomial $f''(w'')$ obtained by the linear transformation $w''=Pw$ is defined as 
$f''(w'')=f(P^{-1}w'')$, and the maximum pole of its zeta function 
$\zeta''(z)=\int f''(w'')^zdw''$ is defined as $-\lambda''$.
Then,

\[\lambda''=\lambda\]
holds.
Let $\smplx''$ be the simplex upper bound of this polynomial $f''$, 
then $\lambda=\lambda''\leq\smplx''$ holds.
Let $p\in\RR^d$ and $P\in\RR^{d\times d}(\det(P)\neq0)$.
The simplex upper bound of polynomials $f'(u)=f(P^{-1}(u+p))$ is denoted by $\smplx^{(P,p)}$.
This changes not continuously but discretely with the change of $p$ or $P$.
Unless noted, $\smplx$ means $\smplx=\smplx^{(I_d,0)}$ hereafter.
The simplex upper bound by linear transformation depends on the matrix.
The following is an example.

\begin{xmp}
Let $f\in\RR[w_1,w_2]_{\geq0}$ be defined as follows.

\[f(w_1,w_2)=w_1^4w_2^2\]
Define $P^{-1}=\bmat{1&1\\1&-1}$. The inverse matrix is $P=-\frac{1}{2}\bmat{-1&-1\\-1&1}$.
Consider the linear transformation $\bmat{w_1\\w_2}=P\bmat{u_1\\u_2}$.

\begin{eqnarray*}
f(u_1+u_2,u_1-u_2)&=&(u_1+u_2)^4(u_1-u_2)^2=(u_1^2+2u_1u_2+u_2^2)(u_1^4-2u_1^2u_2^2+u_2^4)\\
&=&u_1^6-2u_1^4u_2^2+u_1^2u_2^4+2u_1^5u_2-4u_1^3u_2^3+2u_1u_2^5+u_1^4u_2^2-2u_1^2u_2^4+u_2^6\\
&=&u_1^6+2u_1^5u_2-u_1^4u_2^2-4u_1^3u_2^3-u_1^2u_2^4+2u_1u_2^5+u_2^6
\end{eqnarray*}
The simplex upper bound of polynomial above is $\smplx^{(P,0)}=\frac{1}{3}$.

\end{xmp}

The following corollary is obtained from the above discussion.

\begin{cor}
Let $f\in\RR[w_1,w_2,\ldots,w_d]_{\geq0}$.
For the maximum pole $-\lambda$ of the zeta function $\zeta(z)=\int_Wf(w)^zdw$, 
the following holds. Where the integral domain $W$ is a sufficiently large compact set.

\[\lambda\leq\underset{\det(P)\neq0,p\in W}{\min}\ \smplx^{(P,p)}\]
\end{cor}
In general, for two different non-negative polynomials $f,g$ whose RLCT are equal, 
their simplex upper bounds are not necessarily equal.
In other words, the simplex upper bound is not coordinate-invariant.

In addition, the following relation holds between the simplex upper bound and 
the parameter upper bound for nonnegative polynomials.

\begin{prp}
Let $f\in\RR[w_1,w_2,\ldots,w_d]_{\geq0}$ satisfy $f(0)=0$.
and $g\in\calM^d$ can be represented as $g(w)=w^{s-1}$ ($s\in\RR^d_{>0}$).
The following inequation holds for the simplex upper bound.

\[\smplx\leq\frac{\sum_{h=1}^ds_h}{2}\]
\begin{proof}
Let the multi-indexes matrix of $f$ be $A_f=(a_{hj})\in\RR^{d\times n}$.

\[\frac{1}{\smplx}=\underset{\alpha\in S_{d-1}}{\max}\left\{\underset{j=1,2,\ldots,n}{\min}\left(
\sum_{h=1}^d\alpha_h\frac{a_{hj}}{s_h}\right)\right\}\]
Let $\alpha\in S_{d-1}$ achives the maximum of the above equation be $\alpha^*$, we obtain

\[\underset{j}{\min}\left(\sum_{h=1}^d\alpha_h^*\frac{a_{hj}}{s_h}\right)\geq
\underset{j}{\min}\left(\sum_{h=1}^d\frac{s_h}{\sum_{h'=1}^ds_{h'}}\frac{a_{hj}}{s_h}\right)
=\frac{\underset{j}{\min}\left(\sum_{h=1}^da_{hj}\right)}{\sum_{h=1}^ds_h}.
\]
Since $f(0)=0$, $\sum_{h=1}^da_{hj}\geq1$ holds.
Thus, it is sufficient to show $\underset{j}{\min}\sum_{h=1}^da_{hj}\neq 1$.
Suppose $\underset{j}{\min}\sum_{h=1}^da_{hj}=1$ and 
we prove the proposition by causing the contradiction.
There exists $j\in\{1,2,\ldots,n\}$ and unique $h\in\{1,2,\ldots,d\}$ such that

\[a_{ij}=\left\{\begin{array}{lr}1&(i=h)\\0&(i\neq h)\end{array}\right..\]
Fix all components to $0$ except for the $h$-th component, 
and denote obtained polynomial as 
$\hat{f}(w_h)=f(\overbrace{0,\ldots,0}^{h-1},w_h,\overbrace{0,\ldots,0}^{d-h})$.
Then we represent $\hat{f}$ as

\[\hat{f}(w_h)=w_h\left(c_j+\sum_{j'}c_{j'}w_h^{a_{hj'}}\right)\ (c_j,c_j'\neq0,a_{hj'}\in\NN),\]
so $\hat{f}(\epsilon)\hat{f}(-\epsilon)<0$ for a sufficiently small positive number $\epsilon>0$, 
which contradicts $\hat{f}\geq0$.
Therefore the proposition is proved.

\end{proof}
\end{prp}

\begin{cor}
Let true distribution, statistical model, and prior distribution be $q(x),p(x|w)$, and $\varphi(w)$.
We define the mean error function $K(w)$ as follows and 
let $K'(w)$ be a non-negative polynomial whose RLCT $\lambda$ is equal to that of mean error function $K(w)$.

\[K(w)=\int q(x)\log\frac{q(x)}{p(x|w)}dx\]
Let $W\subset\RR^d$ be the parameter set. If there exists an open set $U\subset W$ such that

\[U_0=\left\{w\in U\middle| K(w)=0,\ \varphi(w)>0\right\}\]
is not empty, the following inequation holds for the simplex upper bound $\smplx^{(I_d,p)}$ 
at any $p\in U_0$ of polynomial $K'(w)$.

\begin{equation}\label{eq:rlct_sub_pub}
\lambda\leq\smplx^{(I_d,p)}\leq\frac{d}{2}
\end{equation}
\end{cor}

\subsubsection{comparison with statistical model whose RLCT is revealed}
Equation \eqref{eq:rlct_sub_pub} indicates that the simplex upper bound exists between 
the RLCT and the parameter upper bound.
In order to know the specific location of the simplex upper bound, 
this section compares simplex upper bound with the results of previous studies and 
the parameter upper bound (Theorem \ref{thm:parameter_upper_bound}).
In this study, we compare our results with the following statistical models and singularities.

\begin{itemize}
\item reduced rank regression
\item Poisson mixture
\item Vandermonde matrix type singularities
\end{itemize}

\paragraph{reduced rank regression\cite{rrr}}
In the regression problem of estimating the output $y\in\RR^N$ from the input $x\in\RR^M$, 
the regression function $y=f(x)+\epsilon$($\epsilon$ is normal noise) is called 
\textbf{reduced rank regression} if it is a linear mapping.
Reduced rank regression can be regarded as the neural network whose activation function is an identity function.
It is represented by the following conditional distribution with parameters 
$A\in\RR^{H\times M},B\in\RR^{N\times H}$.

\[p(y|A,B,x)=\frac{1}{(\sqrt{2\pi})^N}\exp\left\{-\frac{\|y-BAx\|_2^2}{2}\right\}\]
$H$ is called the \textbf{rank} in reduced rank regression.
For practical use, models are often constructed with $H<M,N$.
This is thought to reduce unnecessary information and extract only the main information.
The following results are known about the RLCT when learning reduced rank regression in Bayesian learning.

\begin{thm}[RLCT of reduced rank regression]\label{thm:rlct_rrr}
Let the true distribution $q(y|x)$ be the reduced rank regression with rank $r$ and 
statsitical model $p(y|A,B,x)$ be the reduced rank regression with rank $H$.
Where $r\leq H$. Let the prior $\varphi(A,B)$ be $C^\infty$ on the compact support $W$, and 
positive at the true parameter.
Where the true parameter is $w_0=(A_0,B_0)$ such that $q(y|x)=p(y|A_0,B_0,x)$, that is, 
$\varphi(A_0,B_0)>0$ holds.
Then the RLCT is obtained as follows.

\[\lambda=\min\left\{\frac{(N+M)r-r^2+s(N-r)+(M-r-s)(H-r-s)}{2}\middle|0\leq s\leq\min\{M-r,H-r\}\right\}\]
\end{thm}

We use the following polynomial to calculate the simplex upper bound.

\begin{lmm}[equivalent polynomial to mean error function of reduced rank regression]
Consider the same condition as the theorem \ref{thm:rlct_rrr}.
The mean error function

\[K(A,B)=\iint q(x,y)\log\frac{q(x,y)}{p(x,y|A,B)}dxdy\]
and the polynomial

\[K'(C_1,C_2,C_3,A_4,B_4)=\|C_1\|_2^2+\|C_2\|_2^2+\|C_3\|_2^2+\|B_4A_4\|_2^2\]
have the same RLCT. Where $C_1\in\RR^{r\times r}$,$C_2\in\RR^{(N-r)\times r}$,
$C_3\in\RR^{r\times(M-r)}$,$A_4\in\RR^{(H-r)\times(M-r)}$, and $B_4\in\RR^{(N-r)\times(H-r)}$.

\end{lmm}

\begin{figure}[htbp]
\centering
\includegraphics[width=0.8\linewidth]{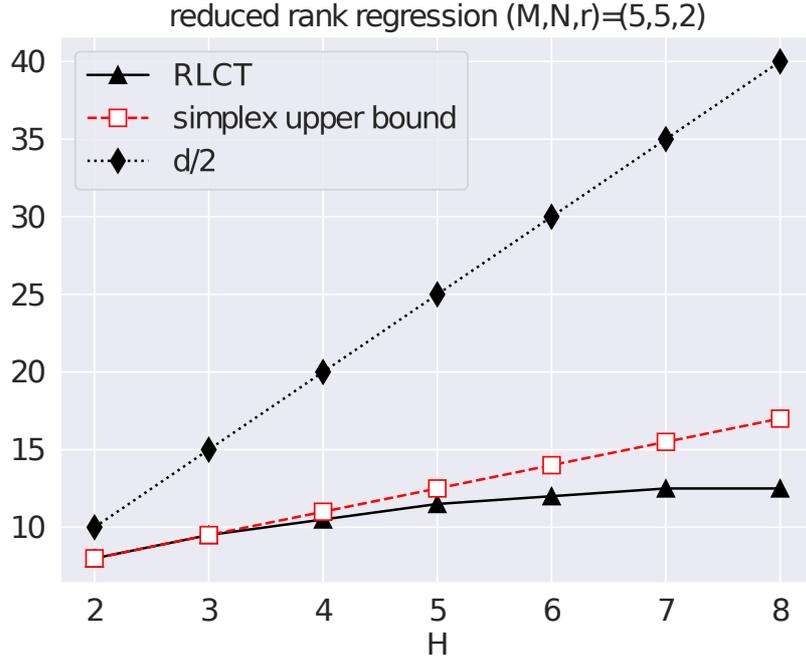}
\caption{
Comparison of the RLCT $\lambda$, the simplex upper bound $\smplx$, 
and the parameter upper bound $d/2$ for reduced rank regression.
$d$ is the number of parameters in the model, where $d=(M+N)H$.
The input dimension $M=5$, the output dimension $N=5$, and 
the mid-dimension of the true distribution $r=2$.
The horizontal axis represents the mid-dimension of the model $H$.}
\label{fig:rrr_comparison1}
\end{figure}

The simplex upper boundary is shown in Figure \ref{fig:rrr_comparison1}.
The simplex upper bound coincides with the RLCT for ranks $H=2,3$.
For $H\geq4$, the simplex upper bound continues to increase with increasing model rank $H$, 
while the RLCT gradually decreases the increasing gap and is flat at $H=7,8$.
The slope of the simplex upper bound is less steep than that of the parameter upper bound $d/2$.

\paragraph{Poisson mixture\cite{pm_arxiv,pm_master}}
The distribution that averages several different $M$-dimensional Poisson distributions 
with probabilistic weights is called \textbf{Poisson mixture}.
Specifically, it is defined as follows with $a\in S_{H-1}$,$b_h\in\RR^M_{>0}(h=1,2,\ldots,H)$.

\[p(x|a,B)=\sum_{h=1}^Ha_h\left(\prod_{m=1}^M\frac{b_{hm}^{x_m}}{x_m!}e^{-b_{hm}}\right)\]
$H$ is called \textbf{the number of components} of the Poisson mixture.
The following results are known for the RLCT when learning Poisson mixture in Bayesian learning.

\begin{thm}[RLCT of Poisson mixture]\label{thm:rlct_pm}
Let the true distribution $q(x)$ be the Poisson mixture with the number of components $r$ and 
the statistical model $p(x|a,B)$ be the Poisson mixture with the number of components $H$.
Let the prior be positive on the compact support $W$.
The RLCT is obtained as follows.

\[\lambda=\frac{3r+H-2}{4}\ \ \ \mathrm{if}\ M=1\]
\[\lambda=\frac{Mr+H-1}{2}\ \ \ \mathrm{if}\ M\geq2\]
\end{thm}

We use the following polynomial to calculate the simplex upper bound.

\begin{lmm}[equivalent polynomial to the mean error function of Poisson mixture]
Consider the same conditions as the theorem \ref{thm:rlct_pm}.
The mean error function

\[K(a,B)=\sum_{x\in\NN_0^M} q(x)\log\frac{q(x)}{p(x|a,B)}\]
and the polynomial

\begin{equation}\label{eq:equivalent_poly_pm}
K'(a,B)=\sum_{x\in\{0,1,\ldots,H+r-1\}^M}\left(\sum_{h=1}^Ha_hb_h^x-\sum_{h=1}^ra_h^*b_h^{*x}\right)^2
\end{equation}
have the same RLCT.
Where $a^*\in S_{r-1},B^*=\bmat{b_1^*&b_2^*&\cdots&b_r^*}\in\RR^{M\times r}_{>0}$ are 
the parameters that satisfy $q(x)=p(x|a^*,B^*)$.

\end{lmm}
Unlike the reduced rank case, the polynomial $K'(a,B)$, 
which is equivalent to the mean error function defined in \eqref{eq:equivalent_poly_pm}, 
needs to be set the parameters for the true distribution.
Since the RLCT is different for $M=1$ and otherwise, 
we compare them in each case.

\begin{figure}[htbp]
\centering
\includegraphics[width=0.8\linewidth]{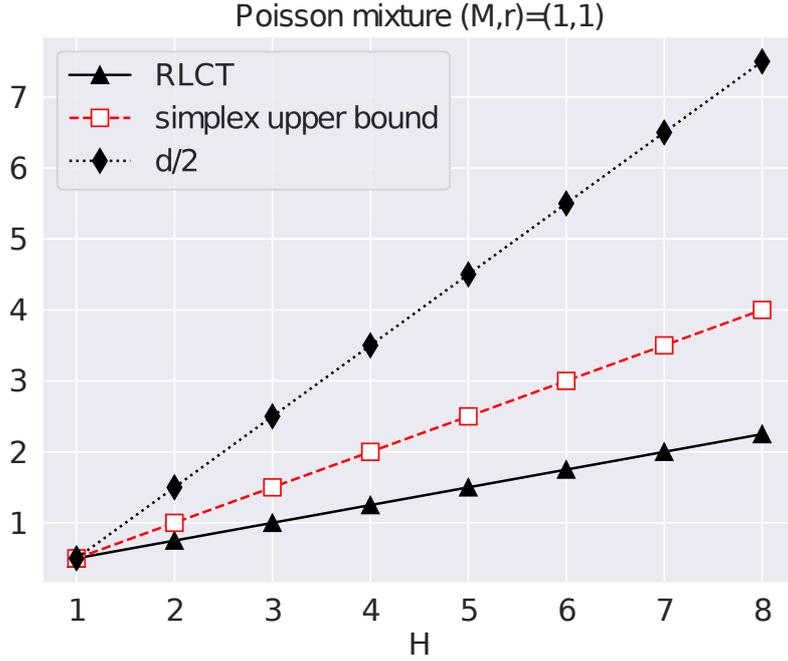}
\caption{Comparison of the RLCT $\lambda$, the simplex upper bound $\smplx$, 
and the parameter upper bound $d/2$ of the Poisson mixture when the data dimension is $M=1$.
$d$ is the number of parameters in the model, where $d=(M+1)H-1$.
The number of components of the true is set as $r=1$.
The horizontal axis represents the number of components in the model $H$.}
\label{fig:pm_comparison1}
\end{figure}

For the case of the data dimension $M=1$, 
consider as an example the case where the number of components of the true $r=1$.
The specific parameters are set as $a^*_1=1$,$b^*_{11}=1$.
The simplex upper bound must satisfy the conditions for Proposition \ref{prp:local_nc_upperbound} to be valid.
From the proof, we assume that all parameters are free variables, 
but the parameter on the simplex $a\in S_{H-1}$ does not satisfy this assumption, 
so a dimension must be reduced.
Specifically, let $a_1=1-\sum_{h=2}^Ha_h$ then deal with$a'=(a_h)_{h=2}^H$ as the free parameters.
Furthermore, since $K'(a',B)$ is not zero at the origin, we do the parallel translation to 
the coordinate whose origin is set as $(a'_2,\ldots,a'_H)^\top=e_{H-1}$,$(b_{11},\ldots,b_{1H})^\top=e_H$.
The simplex upper bound of the polynomial obtained above is shown in Figure \ref{fig:pm_comparison1}.
It can be seen that the RLCT, the simplex upper bound, and the parameter upper bound coincide when $H=1$, 
and they all increase linearly with increasing $H$.
In particular, with respect to the slope, the slope of the simplex upper bound is twice that of RLCT, 
and that of the parameter upper bound is twice that of the simplex upper bound.

\begin{figure}[htbp]
\centering
\includegraphics[width=0.8\linewidth]{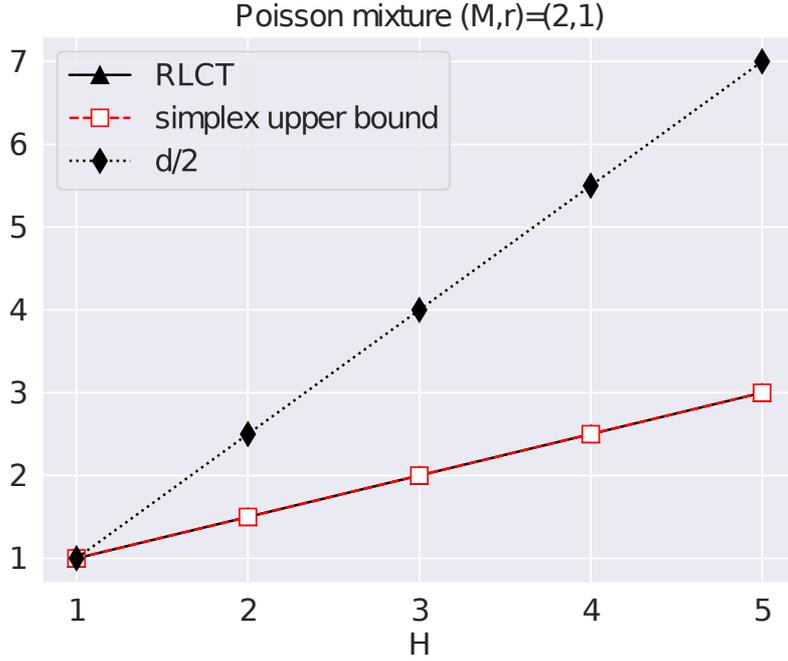}
\caption{Comparison of the RLCT $\lambda$, the simplex upper bound $\smplx$, 
and the parameter upper bound $d/2$ of the Poisson mixture when the data dimension is $M=2$.
$d$ is the number of parameters in the model, where $d=(M+1)H-1$.
The number of components of the true is $r=1$.
The horizontal axis represents the number of components in the model $H$.}
\label{fig:pm_comparison2}
\end{figure}

Next, we consider the case where the data dimension is $M=2$ and 
the number of components of the true is $r=1$ as an example.
We set $a^*_1=1$,$b^*_{11}=b^*_{21}=1$ as the specific parameter settings.
In order to reduce the apparent dimension, substituting $a_1=1-\sum_{h=2}^Ha_h$.
Also, $K'(a',B)$ is not $0$ at the origin, we do the parallel translation to the coordinate whose origin are
$(a'_2,\ldots,a'_H)^\top=e_{H-1}$,$(b_{m1},\ldots,b_{mH})^\top=e_H(m=1,2)$.
The simplex upper bound of the polynomial obtained above is shown in Figure \ref{fig:pm_comparison2}.
Unlike the case $M=1$, 
the simplex upper bound coincides with the RLCT perfectly.

\paragraph{Vandermode matrix type singularities \cite{vandermonde}}
VanderMonde matrix type singularities have the same algebraic structure as 
the mean error function of a three-layer neural network whose activation function 
is the hyperbolic tangent function and the normal mixture.
The definitions are given below.

\begin{dfn}[Vandermonde matrix type singularities]\label{dfn:vandermonde_matrix_type_singularities}
Let $M,H,N,Q\in\NN$,$m,r\in\NN_0$. We define the matrix $A$ as follows.

\[A=\bmat{a_{11}&\cdots&a_{1H}&a^*_{1,H+1}&\cdots&a^*_{1,H+r}\\
\vdots&\ddots&\vdots&\vdots&\ddots&\vdots\\
a_{M1}&\cdots&a_{MH}&a^*_{M,H+1}&\cdots&a^*_{M,H+r}}\in\RR^{M\times(H+r)}\]
For $L=(l_1,\ldots,l_N)\in\NN_0^N$, define the vector $B_L$ as follows.

\[B_L=\bmat{\prod_{j=1}^Nb_{1j}^{l_j}&\cdots&\prod_{j=1}^Nb_{Hj}^{l_j}&
\prod_{j=1}^N{b^*_{H+1,j}}^{l_j}&\cdots&\prod_{j=1}^N{b^*_{H+r,j}}^{l_j}}^\top\]
The matrix $B$ consists of the above vectors.

\[B=\bmat{B_L}_{|L|=Qn+m,0\leq n\leq H}\in\RR^{(H+r)\times(N^*)}\]
Where $N^*=\left|\left\{(l_1,\ldots,l_N)\in\NN_0^N\middle|
\sum_{j=1}^Nl_j=Qn+m,\ 0\leq n\leq H\right\}\right|$.
Then, the singularities of the ideal generated by the components of the matrix $AB$ is
referred to as \textbf{Vandermode matrix type singularities}.
Note that $a_{mh},b_{hn}$ are variables and $a^*_{m.H+h},b^*_{H+h,n}$ are constants.

\end{dfn}

Define the function $\sigma(x)=\sum_{n=0}^Hx^{Qn+m}$ for $x\in\RR$.
The component of $AB$ in the definition is obtained by Talor expansion of the following function.

\[f_m(a_{m*},B)=\sum_{h=1}^Ha_{mh}\sigma\left(\sum_{j=1}^Nb_{hn}\right)+
\sum_{h=H+1}^{H+r}a^*_{mh}\sigma\left(\sum_{j=1}^Nb^*_{hn}\right)\ \ \ (m=1,\ldots,M)\]

The following results are known for the RLCT of $\|AB\|^2$ when $r=0$ and $m=1$.

\begin{thm}[RLCT of Vandermode matrix type singularities]
Let the matrix $A,B$ be defined in the Definition \ref{dfn:vandermonde_matrix_type_singularities}.
let $r=0,m=1$.
For the maximum pole $-\lambda$ of the eta function defined below

\[\zeta(z)=\int_W\|AB\|_2^2dAdB,\]
The following hold, where $W$ is a sufficiently large compact set.

\begin{enumerate}
\item for the case of $H=1$, $\lambda=\min\left\{\frac{M}{2},\frac{N}{2}\right\}$
\item for the case of $H=2$, $\lambda=\min\left\{\frac{\beta N+(2-\beta)M}{2}\middle|\beta=0,1,2\right\}
	\cup\left\{\frac{2N+Q(N-1+M)}{2Q+2}\right\}$
\item for the case of $H=3$, $\lambda=\min\left\{\frac{\beta N+(3-\beta)M}{2}\middle|\beta=0,1,2,3\right\}$\\
	$\cup\left\{\frac{\beta N+(3-\beta)M+Q(\alpha(N+\alpha-\beta)+(3-\alpha)M)}{2(Q+1)}\middle|
	\alpha=1,\ldots,\beta-1,\beta=2,3\right\}\cup
	\left\{\frac{3N+Q(3N-3+3M)}{2(2Q+1)}\right\}$
\item for the case of $H=4$, $\lambda=\min\left\{\frac{\beta N+(4-\beta)M}{2}\middle|\beta=0,1,2,3,4\right\}$\\
	$\cup\left\{\frac{\beta N+(4-\beta)M+Q(\alpha(N+\alpha-\beta)+(4-\alpha)M)}{2(Q+1)}\middle|
	\alpha=1,\ldots,\beta-1,\beta=2,3,4\right\}$\\
	$\cup\left\{\frac{4N+Q(\alpha N-\alpha-1+(8-\alpha)M)}{2(2Q+1)}\middle|\alpha=2,3,4\right\}
	\cup\left\{\frac{4N+Q(5N-5+3M)}{2(2Q+1)}\right\}$\\
	$\cup\left\{\frac{3N+M+Q(\alpha N-\alpha+(8-\alpha)M)}{2(2Q+1)}\middle|\alpha=2,3\right\}
	\cup\left\{\frac{4N+Q(6N-6+6M)}{2(3Q+1)}\right\}$
\item for the case of $N=1$, let $k=\max\{i\in\ZZ|2H\geq M(i(i-1)Q+2i)\}$, then
	$\lambda=\frac{MQk(k+1)+2H}{4(1+kQ)}$
\end{enumerate}

\end{thm}

For $N\geq2$, the RLCT $\lambda$ is revealed for $H\leq4$ and 
for $N=1$, the RLCT is revealed for all $H$.
We compare them in each case.
The polynomial for the simplex upper bound is $\|AB\|_2^2$ itself.

\begin{figure}[htbp]
\centering
\includegraphics[width=0.8\linewidth]{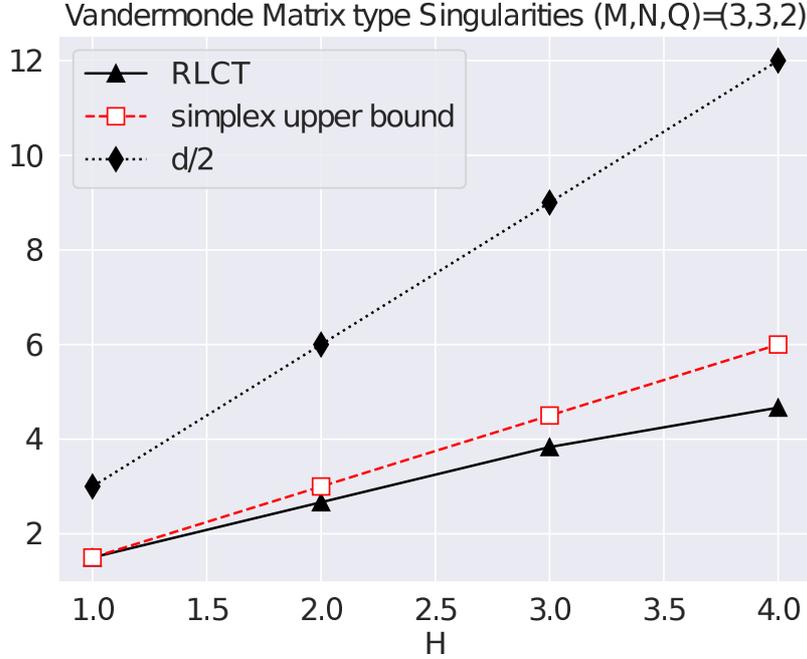}
\caption{Comparison of the RLCT $\lambda$, the simplex upper bound $\smplx$, 
and the parameter upper bound $d/2$ of the Vandermonde matrix type singularities when $N=3$.
$d$ is the number of variables in the polynomial, where $d=(M+N)H$.
The specific structural setting is $(M,Q)=(3,2)$.
The horizontal axis represents the mid-dimension $H$.}
\label{fig:vmts_comparison1}
\end{figure}

The simplex upper bound for $N\geq2$ is shown in Figure \ref{fig:vmts_comparison1}.
The simplex upper bound coincides with the RLCT only when the model rank $H=1$.
As the mid-dimension $H$ increases, the simplex upper bound increases linearly, 
though the slope of the RLCT becomes slower.
It is also found that the simplex upper bound is half of the parameter upper bound.

\begin{figure}[htbp]
\centering
\includegraphics[width=0.8\linewidth]{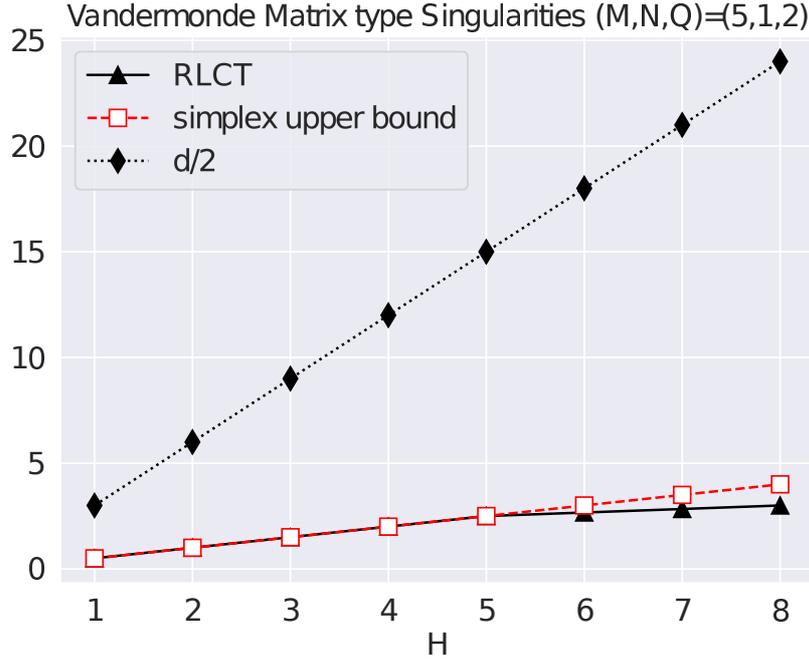}
\caption{Comparison of the RLCT $\lambda$, the simplex upper bound $\smplx$, 
and the parameter upper bound $d/2$ of a Vandermonde matrix type polynomial when $N=1$.
$d$ is the number of variables in the polynomial, where $d=(M+N)H$.
The specific structural setting is $(M,Q)=(5,2)$.
The horizontal axis represents the mid-dimension $H$.}
\label{fig:vmts_comparison2}
\end{figure}

The simplex upper bound for $N=1$ is shown in Figure \ref{fig:vmts_comparison2}.
The simplex upper bound coincides with the RLCT when the model rank $H\leq 5$.
For $H\geq 6$, as the mid-dimension $H$ increases, 
the simplex upper bound increases with the same slope as $H\leq 5$ though 
the RLCT is constant.
It is also found that the slope of the simplex upper bound is lower than that of the parameter upper bound.

\subsection{Application for weighted blowup}\label{sec:apply_for_weighted_blowup}
In this section, we show that the LPP \eqref{eq:mi_poly_lp} can be applied to weighted blowups.
First, we define a weighted blowup\cite{weighted}.
Let $q=(q_1,q_2,\ldots,q_d)^\top\in\NN_0^d$ be an integer vector with at least two positive components.
Let $I=\{i|q_i\neq0\}$. We define a map $g_i:V_i\rightarrow U$ as follows.

\[
%\begin{equation}\label{eq:weighted_blowup_var_trans}
g_i=\left\{\begin{array}{rcl}
x_1&=&x_1^{(i)}x_i^{q_1}\\
x_2&=&x_2^{(i)}x_i^{q_2}\\
&\vdots&\\
x_{i-1}&=&x_{i-1}^{(i)}x_i^{q_{i-1}}\\
x_i&=&x_i^{q_i}\\
x_{i+1}&=&x_{i+1}^{(i)}x_i^{q_{i+1}}\\
&\vdots&\\
x_d&=&x_d^{(i)}x_i^{q_d}
\end{array}
\right.
%\end{equation}
\]
We can past the destinations together by using $g_i(i\in I)$ and 
make a non-singular manifold $V=\bigcup_{i\in I}V_i$.
$(V,\{g_i\}_{i\in I})$ is called a \textbf{weighted blowup} with weight $q$.
An example is as follows.

\begin{xmp}\label{xmp:weighted_blowup}
Let $H(w)=w_1^2+w_2^4+w_3^6$.
We consider the weighted blowup with weight $q=(6,3,2)^\top$ for $H$.

\[g_1:(w_1,w_2,w_3)\mapsto(w_{11}^6,w_{11}^3w_{21},w_{11}^2w_{31})\]
\[g_2:(w_1,w_2,w_3)\mapsto(w_{22}^6w_{12},w_{22}^3,w_{22}^2w_{32})\]
\[g_3:(w_1,w_2,w_3)\mapsto(w_{33}^6w_{13},w_{33}^3w_{23},w_{33}^2)\]
Thus we obtain the following on each local coordinate.

\[H(g_1(w))=w_{11}^{12}+w_{11}^{12}w_{21}^4+w_{11}^{12}w_{31}^6=w_{11}^{12}(1+w_{21}^4+w_{31}^6)\]
\[H(g_2(w))=w_{22}^{12}w_{12}^2+w_{22}^{12}+w_{22}^{12}w_{32}^6=w_{22}^{12}(w_{12}^2+1+w_{32}^6)\]
\[H(g_3(w))=w_{33}^{12}w_{13}^2+w_{33}^{12}w_{23}^4+w_{33}^{12}=w_{33}^{12}(w_{13}^2+w_{23}^4+1)\]
\end{xmp}
The problem in a weighted blowup is what weight to give to a given polynomial $f$.
If a simple polynomial is given, as in Example \ref{xmp:weighted_blowup}, 
we can intuitively consider the weights to be given, 
but if a complex polynomial is given, it is difficult to intuitively consider the weights.
The question of what weights should be given is different for each set of objectives.
In this section, we consider optimal weights in the following sense.

\begin{dfn}
Let $f\in\RR[w_1,w_2,\ldots,w_d]_{\geq0}$ and $h\in\calM^d$.
They are represented as follows.

\[f(w)=\sum_{j=1}^nc_jw^{a_{*j}}\ ,\ h(w)=w^{s-1_d}\]
Where $c_j\neq0$,$a_{*i}\in\NN_0^d,a_{*j}\neq a_{*j'}(j\neq j')$,$s\in\RR^d_{>0}$.
Suppose that the following is obtained polynomial on $w_i$-chart by blowing up with weight $q\in\NN_0^d$.

\[f(g_i(w))=\sum_{j=1}^nc_jw_i^{q^\top a_{*j}}w_{\backslash i}^{a_{*_{\backslash i}j}}=
w_i^{\underset{j}{\min}(q^\top a_{*j})}\left(
\sum_{j=1}^nc_jw_i^{q^\top a_{*j}-\underset{j}{\min}
(q^\top a_{*j})}w_{\backslash i}^{a_{*_{\backslash i}j}}\right)\]
\[h(g_i(w))|g_i|=w_i^{q^\top(s-1_d)}w_{\backslash i}^{s_{\backslash i}-1_{d-1}}w_i^{q^\top1_d-1}=
w_i^{q^\top s-1}w_{\backslash i}^{s_{\backslash i}-1_{d-1}}\]
Then the quantity 
$\mu_{(f,h)}(q)=(\underset{j=1,2,\ldots,n}{\min}q^\top a_{*j})/(q^\top s)$ for multi-indexes 
is referred as to \textbf{minimum index ratios} with weight $q$ for polynomial $(f,h)$.

\end{dfn}
The minimum indexes ratio $\mu_{(f,g)}(q)$ corresponds to the indexes ratio 
at the origin of $(f(g_i(w)),h(g_i(w))(i\in I)$.
The minimum index ratio is closely related to the simplex upper bound that 
gives the lower bound of the inverse of the RLCT.
Consider the problem to find $q$ such that the minimum index ratio is maximized.
This paper can provide the answer to this problem. The following proposition holds.

\begin{prp}\label{prp:optimize_alpha}
Let $f\in\RR[w_1,w_2,\ldots,w_d]_{\geq0}$ and $h\in\calM^d$.
$h(w)$ is represented as $h(w)=w^{s-1}(s\in\RR^d_{>0})$.
The weight $q\in\NN_0^d(q^\top 1_d\geq2)$ that 
maximizes the minimum index ratio $\mu_{(f,h)}(q)$ for polynomial $(f,g)$ 
is the integer multiple of $(\alpha^*_h/s_h)_{h=1}^d\in\QQ^d_{\geq0}$, 
where $\alpha^*\in S_{d-1}$ is the optimal solution of LPP \eqref{eq:mi_poly_lp}.
We also obtain the following for minimum index ratios.

\[\mu_{(f,h)}(q)\leq\smplx\]
The equality holds for $q_h=k\alpha^*_h/s_h(h=1,2,\ldots,d)$, where $k\in\NN$.

\begin{proof}
The proposition follows from a part of the proof of the Theorem \ref{thm:dn_rlct_upper_bound}.

\end{proof}
\end{prp}

\section{Conclusions}

This paper reveals the following for non-negative sop polynomials.

\begin{itemize}
\item For nonnegative sum-of-products binomials, 
we propose an algorithm that forms a blowup tree with normal crossings in all local coordinates.
We also derived the exact value of the RLCT by using the blowup tree formed by the proposed algorithm.
(Sections \ref{sec:22_case},\ref{sec:d2_case})
\item For nonnegative sum-of-products polynomials, 
we define a local normal crossing, which is more relaxed than a normal crossing.
Using the algorithm for the case of binomials, we propose an algorithm to form a blowup tree 
with local normal crossings in all local coordinates and show that it halts.
Using the blowup tree formed by the algorithm, we derived the upper bound of the RLCT.
(Section \ref{sec:dn_case})
\item We confirm that the upper bound of the RLCT of sop polynomials is also applicable to 
general polynomials, that is, the inverse of RLCT, or simplex upper bound, 
can be obtained by solving the LPP \eqref{eq:mi_poly_lp} based on multi-indexes of polynomials.
Furthermore, we prove that the simplex upper bound is less than or equal to the parameter upper bound, 
and confirm that the simplex upper bound is truly tighter than the parameter upper bound 
for some polynomials or mean error functions by some experiments.
We also show that the optimal weight in the weighted blowup can be calculated from 
the optimal solution of the LPP \eqref{eq:mi_poly_lp}.
(Sections \ref{sec:simplex_upper_bound},\ref{sec:apply_for_weighted_blowup})
\end{itemize}

\section{Acknowledgement}

Finally, I would like to thank all the people who have helped me in conducting this study.
I would like to thank my advisor, Professor Sumio Watanabe, for his careful answers to my questions about blowups, which helped me to deepen my understanding of the theory of blowups.
I am also grateful to Atsuyoshi Muta and Yoshiki Mikami, 
who have already graduated from the Watanabe Laboratory, 
for their answers to my questions and suggestions on the proof of the main theorem of their master theses, 
which are previous studies. I would like to express my sincere gratitude to them.

\begin{comment}
\appendix
\section{}
\end{comment}

\end{document}